\documentclass[12pt ]{article}

\newtheorem{theorem}{Theorem}[section]
\newtheorem{lemma}{Lemma}[section]

\newtheorem{remark}{Remark}[section]
\newcommand{\eqnsection}{
   \renewcommand{\theequation}{\thesection.\arabic{equation}}
   \makeatletter
   \csname @addtoreset\endcsname{equation}{section}
   \makeatother}
\def \ov{\overline}

\def \be{\begin{equation}}
\def \ee{\end{equation}}
\def \bt{\begin{theorem}}
\def \et{\end{theorem}} 
\def \bea{\begin{eqnarray}}
\def \eea{\end{eqnarray}}
\def \bas{\begin{eqnarray*}}
\def \eas{\end{eqnarray*}}
\def \bl{\begin{lemma}} 
\def \el{\end{lemma}}
\def \al{\alpha}
\def \bb{\beta}
\def \ga{\gamma}

\def \de{\delta}
\def \De{\Delta}
\def \ep{\epsilon}

\def \la{\lambda}

\def \si{\sigma}

\def \th{\theta}

\def\stl{\stackrel{\mathcal{L}}{=}}
\def \ff{\infty}
\def \wh{\widehat}
\def \wt{\widetilde}

\def \rar{\rightarrow}

\def \cd{\,\cdot\,}

\def \AA{{\cal A}}
\def \BB{{\cal B}}

\def \II{{\cal I}}

\def \MM{{\cal M}}
\def \NN{{\cal N}}

\def \PP{{\cal P}}

\def \RR{{\cal R}}

\def \ZZ{{\cal Z}}

\def \p{p_{t}(x,y)}

\def \({\left(}
\def \){\right)}

\def \lc{\left\{}
\def \rc{\right\}}

 \def \nn{\nonumber}
\def \Proof{\noindent{\bf Proof $\,$ }}

\def \bc{\begin{center} }
\def \ec{\end{center} }
\def \bs{\begin{slide} }
\def \es{\end{slide} }

\def\square{{\vcenter{\vbox{\hrule height.3pt
        \hbox{\vrule width.3pt height5pt \kern5pt
           \vrule width.3pt}
        \hrule height.3pt}}}}
\def\qed{{\hfill $\square$ \bigskip}}

\eqnsection
\begin{document}

\def\wh{\widehat}
\def\ol{\overline}

\title{A CLT  for the $L^{2}$ norm of increments  of   local times of L\'evy processes as time goes to infinity}

\author{   Michael B. Marcus 
 \hspace{ .2in}  Jay Rosen\thanks
 {The research of both authors was  supported, in part, by grants from the National Science
Foundation and PSC-CUNY.}}

%\date{November 20, 2005}

\maketitle

\bibliographystyle{amsplain}

\begin{abstract} 
Let $X=\{X_{t},t\in R_{+}\}$   be a symmetric L\'{e}vy process 
with  local time   $\{L^{ x }_{ t}\,;\,(x,t)\in R^{ 1}\times  R^{  1}_{ +}\}$. 
When the L\'{e}vy exponent $\psi(\la)$  is regularly varying at zero with   index $1<\bb\leq 2$, and satisfies some additional regularity conditions, 
\bea
&& { \int_{-\ff}^{\ff} ( L^{ x+1}_{t}- L^{ x}_{ t})^{ 2}\,dx- E\(\int_{-\ff}^{\ff} ( L^{ x+1}_{t}- L^{ x}_{ t})^{ 2}\,dx\)\over t\sqrt{\psi^{-1}(1/t)}}\label{r5.0tweaksabs}\\
&&\hspace{1  in}
\stackrel{\mathcal{L}}{\Longrightarrow}(8c_{\psi,1 })^{1/2}\(\int_{-\ff}^{\ff} \(L_{\bb,1}^{x}\)^{2}\,dx\)^{1/2}\,\eta\nn
\eea
as $t\rar\ff$, where    $L_{\bb,1}=\{L^{ x }_{\bb, 1}\,;\, x \in R^{ 1} \}$ denotes the local time, at time 1, of   a symmetric stable process with index $\bb$,   
$\eta$ is a normal random variable with mean zero and variance one that is independent of $L _{ \bb,1}$, and $c_{\psi,1}$ is a known constant that depends on $\psi$.

\end{abstract}



\section{Introduction}

The earliest  result we know about the asymptotic behavior   in time, of  increments of  local times  in the spatial variable,  is due to Dobrushin, \cite{D}. Let $ \{S_{ n} \,;\,n=0,1,2,\ldots\}$ be a simple random walk on
$Z^{ 1}$ and let  $ \ell^{ x }_{ n}=\sum_{j=1}^{n}\,1_{\{S_{j}=x\}} $
denote its the local time.    Dobrushin shows that
\be
{\ell_{n}^{1} -\ell_{n}^{0}\over n^{1/4}}\stackrel{\mathcal{L}}{\Longrightarrow}  (2|Z|)^{1/2}\,\,\eta \label{dob}\ee
as $n\rar \ff$, where $Z$ and $\eta$ are independent normal random variables with mean zero and variance one. Two aspects of this result are relevant to  this paper. One is that, since $\ell_{n}^{0}$ grows like $n^{1/2}$, see for example \cite[(10.1), (9.13)]{Revesz}, (\ref{dob}) is a result about fluctuations. The other is that the right-hand side of (\ref{dob})  is not a standard normal random variable, but is the product of a standard normal random variable and an independent random variable.
Extensions of (\ref{dob}) to the local time of Brownian motion and other processes can be found in  R\'{e}v\'{e}sz, \cite[(11.10), (12.17), (12.19)]{Revesz}, Marcus and Rosen, \cite{MR2,MR3}, Rosen, \cite{Rosen} and Yor, \cite{Yor}.

 Our   motivation for considering  the increments of local times is our interest in    the  Hamiltonian  for the  critical attractive random polymer in one dimension,  
\cite{HK,HKK},
\be
H_{ n}= \sum_{x\in Z^{ 1}}\(\ell_{ n}^{ x+1}-\ell_{ n}^{ x}\)^{ 2}\label{1.12qq},
\ee
where  $\ell_{n}^{x}$ is  as  defined above. This is the  square of the $\ell^{2}$ norm of the increments of the local time at time $n$.

We began our study of expressions like (\ref{1.12qq}) in  \cite{CLR}, with X.  Chen and W. Li, by considering the continuous version of this problem for Brownian motion. 
Let   $ \{L^{ x }_{ t}\,;\,(x,t)\in R^{ 1}\times  R^{  1}_{ +}\}$ denote the local times of Brownian motion. We show that  
 \begin{equation} { \int_{-\ff} ^{\ff}( L^{ x+1}_{t}- L^{ x}_{ t})^{ 2}\,dx- 4t\over t^{ 3/4}}
\stackrel{\mathcal{L}}{\Longrightarrow}\( 64 / 3 \)^{ 1/2}\(\int_{-\ff} ^{\ff} \(L_{1}^{x}\)^{2}\,dx\)^{1/2}\,\,\eta\label{5.0tweak}
\end{equation}
as $t\rar\ff$, where   $\eta$ is a normal random variable with mean zero and variance one that is independent of $\{L^{x}_{1},x\in  R^{1} \}$.

The  proof in  \cite{CLR} makes extensive use of the scaling property of Brownian motion. A different proof in \cite{Rosen1} uses   stochastic integrals and the 
 theorem of Papanicolaou, Stroock, and Varadhan, \cite[Chapter XIII]{RY}. Neither of these approaches can   be used to extend (\ref{5.0tweak}) to general L\'evy processes.  We do that in this paper using the method of moments. 
 
  Let $X=\{X_{t},t\in R_{+}\}$ \label{page1} be a   symmetric L\'evy process with characteristic function 
\begin{equation}
E\(e^{i\la X_{t}}\)=e^{-\psi(\la)t}\label{st.1},
\end{equation}
and local time which we continue to denote by   $\{L^{ x }_{ t}\,;\,(x,t)\in R^{ 1}\times  R^{  1}_{ +}\}$.  
The  behavior of a suitably scaled version of  $\int ( L^{ x+1}_{t}- L^{ x}_{ t})^{ 2}\,dx$ as $t$ goes to infinity depends primarily  on the behavior of $\psi(\la)$ as $\la$ goes to $0$. This is not   surprising, since  large time properties of $X$ such as transience and recurrence depend on 
 the behavior of $\psi(\la)$ as $\la$ goes to $0$; (see \cite[Chapter 1, Theorem 17]{Bertoin}, which shows, in particular, that the processes we consider are recurrent.)
 
We assume that  $\psi(\la)$ satisfies the following conditions:
 \be
\hspace{-.9 in}1. \hspace{.3 in}\psi(\la)  \mbox{ is regularly varying at 0 with   index $1<\bb\leq 2$} ; \label{regcond}
 \ee
 \begin{equation}
  \hspace{-3 in}2.\hspace{.3 in}\int_{-\ff}^{\ff} {1 \over 1+\psi(\la)}\,d\la<\ff;\label{regcond2}
 \end{equation}
\hspace{.025 in}3. $\psi$ is twice differentiable almost everywhere,  and   there exist  constants $D_{1},D_{2}<\ff$ such that for $0<\la\le 1$ 
 \be
\la|\psi'(\la)|\le D_{1}\psi(\la)\quad\mbox{and}\quad    \la^{2}  |\psi''(\la)|\le D_{2}\psi(\la) \label{88.m}
   \ee
 and
   \begin{equation}
   \int_{1}^{\ff}\frac{|\psi'(\la)| }{\psi^{2}(\la)}\,d\la<\ff,  \quad    \int_{1}^{\ff}\frac{|\psi'(\la)| ^{2}}{\psi^{2}(\la)}\,d\la<\ff,\quad  \int_{1}^{\ff}\frac{|\psi''(\la)|}{\psi (\la)}\,d\la<\ff.\label{1.12}
   \end{equation}
 (Condition 1. is substantive. Condition 2. is the necessary and sufficient condition for a symmetric L\'evy process to have a local time. The criteria in Condition 3.  are  rather weak.)

\medskip	 We prove the following theorem:
\begin{theorem}\label{theo-clt2r} Let $\{L^{ x }_{ t}\,;\,(x,t)\in R^{ 1}\times  R^{  1}_{ +}\}$ be the local time  of   a symmetric  L\'{e}vy process $X$,
with  L\'{e}vy exponent
 $\psi(\la)$, that is regularly varying at zero with index $1<\bb\le 2$ and satisfies (\ref{regcond2})--(\ref{1.12}). Then   
\bea
&& { \int  _{-\ff}^{\ff} ( L^{ x+1}_{t}- L^{ x}_{ t})^{ 2}\,dx- E\(\int _{-\ff}^{\ff} ( L^{ x+1}_{t}- L^{ x}_{ t})^{ 2}\,dx\)\over t\sqrt{\psi^{-1}(1/t)}}\label{r5.0tweaks}\\
&&\hspace{1.5 in}
\stackrel{\mathcal{L}}{\Longrightarrow}(8c_{\psi,1 })^{1/2}\(\int _{-\ff}^{\ff} \(L_{\bb,1}^{x}\)^{2}\,dx\)^{1/2}\,\eta\nn
\eea
as $t\rar\ff$, where  $L_{\bb,1}=\{L^{ x }_{\bb, 1}\,;\, x \in R^{ 1} \}$   is the local time,  at time 1, of   a symmetric stable process of index $\bb$,  $\eta$ and $L _{\bb, 1}$ are independent, and 
\begin{equation}
  c_{\psi ,1}  =   {16  \over \pi}\int_{0}^{\ff} { \sin^{4} p  /2  \over \psi ^{2} (p)}\,dp.\label{1.16aq}
\ee  

\end{theorem}

  (Since $\psi$ is regularly varying at zero, it is asymptotic to a monotonic function at zero. We define   $\psi^{-1} $ as the inverse of this  function.)
  
  It follows from Lemma 5.1, in this paper, that 
  \be
   E\(\int _{-\ff}^{\ff} ( L^{ x+1}_{t}- L^{ x}_{ t})^{ 2}\,dx\)=4c_{\psi,0}t+o\( t\sqrt{\psi^{-1}(1/t)}\),\label{1.11}
   \ee  where
\be
c_{\psi ,0}: =   {2   \over \pi}\int_{0}^{\ff} { \sin^{2} (p /2)  \over \psi  (p)}\,dp\label{1.16a}.
   \end{equation} 
   Therefore, we can replace the mean in (\ref{r5.0tweaks}) by $4c_{\psi,0}t$. 
   
Note that by Lemma \ref{lem-ilt},   $ \int  _{-\ff}^{\ff} ( L^{ x}_{t})^{ 2}\,dx$ grows like $t^{2}\psi^{-1}(1/t)$, therefore  (\ref{r5.0tweaks}) is also a fluctuation result.
 
In Remark \ref{rem-2.1}, we evaluate the constants and make the necessary changes to verify  that when $X$ is Brownian motion, (\ref{r5.0tweaks}) along with (\ref{1.11}), is the same as (\ref{5.0tweak}).

\medskip	 The methods developed to prove Theorem \ref{theo-clt2r}  also give  a Central Limit Theorem for   the  Hamiltonian of the  critical attractive random polymer in one dimension, which is usually written as
\begin{equation} H_{ n}=2\sum_{ i,j=1}^{ n}1_{ \{S_{ i}=S_{  j} \}}- \sum_{ i,j=1 }^{ n}1_{
\{|S_{ i}-S_{ j}|=1 \}}.\label{rp5c.4}
\end{equation}
 It is easy to see that this is the same as (\ref{1.12qq}).    In \cite{cp} we   use the methods of this paper   to show that \begin{equation}  {H_{n} - 4n\over n^{ 3/4}}
\stackrel{\mathcal{L}}{\Longrightarrow}\( 12 \)^{ 1/2}\(\int_{-\ff} ^{\ff}
 \(L_{ 1}^{x}\)^{2}\,dx\)^{1/2}\,\,\eta\label{5.0tweakrw}
\end{equation}
as $n\rar\ff$, where       $\{L^{x}_{ 1},x\in  R^{1} \}$ is the local time of   Brownian motion at time 1. A similar formula holds for symmetric random walks with  variance 1.  In fact in \cite{cp} we  obtain central limit theorems like (\ref{5.0tweakrw}) for a large class of random walks in the domain of attraction of    stable processes of index $1<\bb\le 2$.   (The conditions are the ones in (\ref{regcond})--(\ref{1.12}) with obvious modifications.)

\medskip	Another  application of the techniques used to prove Theorem \ref{theo-clt2r} is suggested by applying the scaling relationship for the local times of $\bb$ stable processes, 
\begin{equation}
\{ L^{ x}_{\bb,  t/\de ^{\bb}}\,;\,( x,t)\in R^{ 1}\times R^{ 1}_{ +}\}
\stackrel{\mathcal{L}}{=}\{\de ^{ -(\bb-1)} L^{\de x}_{\bb,  t}\,;\,( x,t)\in R^{ 1}\times
R^{ 1}_{ +}\},\label{scales}
\ee
(see   e.g. \cite[Lemma 10.5.2]{book}), to (\ref{r5.0tweaks}). This   gives
 \bea&& { \int_{-\ff}^{\ff} ( L^{ x+h}_{\bb, 1}- L^{ x}_{\bb,  1})^{ 2}\,dx-4  c_{|p|^{\bb},0} h^{\bb-1} \over h^{(2\bb-1)/2}}\label{5.0weaksqq}\\
 &&\qquad\qquad
\stackrel{\mathcal{L}}{\Longrightarrow} (8c_{|p|^{\bb},1})^{1/2}\(\int _{-\ff}^{\ff}\(L_{\bb, 1}^{x}\)^{2}\,dx\)^{1/2}\,\eta \nn,
\eea
as $h\rar 0$. 

  One can not use scaling to obtain limits like (\ref{5.0weaksqq}) for more general  L\'evy processes. In \cite{MR} we do this by modifying the methods developed  in this paper.  We show that when   $\psi(\la)$ is regularly varying at infinity, with index $1<\bb\le 2$, and satisfies some additional regularity conditions, \bea
&& \sqrt{h\psi^{2}(1/h)} \lc \int_{-\ff}^{\ff} ( L^{ x+h}_{1}- L^{ x}_{ 1})^{ 2}\,dx- E\( \int_{-\ff}^{\ff} (  L^{ x+h}_{1}- L^{ x}_{ 1})^{ 2}\,dx\) \rc\nn\\
&&\hspace{1.5 in}
\stackrel{\mathcal{L}}{\Longrightarrow}  (8 c_{{|p|^{\bb},1}})^{1/2}\(\int_{-\ff}^{\ff} \(L_{1}^{x}\)^{2}\,dx\)^{1/2}\,\eta\label{1.2},
\eea
as $h\rar 0$.  It is not surprising that this limit   depends  on the behavior of $\psi(\la)$ as $\la$ goes to infinity, since  
 the behavior of $\psi(\la)$ as $\la$ goes to infinity controls the small jumps of $X$.

\medskip	 In Section \ref{sec-CLT} we   give the proof of Theorem \ref{theo-clt2r} based on four lemmas. The remainder of the paper is devoted to the proof of these lemmas. The critical ingredient in the proof  of  Theorem \ref{theo-clt2r} is Lemma \ref{lem-multiple} which gives moments for the  $L^{2}$ norm of increments of local times of L\'evy processes satisfying (\ref{regcond})--(\ref{1.12}). 
  In Section \ref{sec-indf} we state  Lemmas \ref{lem-vproprvtj}--\ref{lem-2.5} which give properties of $p_{t}(x)$, $\De^{1}p_{t}(x)$ and $\De^{1}\De^{-1}p_{t}(x)$, where $p_{t}(x)$ is the transition probability density of   these L\'evy processes. The lemmas are used in Section \ref{sec-3} in the proof of Lemma \ref{lem-multiple}. They are proved in Section \ref{sec-lemproofs}.  In the proof of  Theorem \ref{theo-clt2r} we need good asymptotic estimates of   the mean and variance of $ \int_{-\ff}^{\ff} ( L^{ x+h}_{1}- L^{ x}_{1})^{2}\,dx$. 
These are stated in   Section \ref{sec-CLT} and proved in Section \ref{sec-9}.   In  Section \ref{sec-Kac}, a short appendix, we derive   the form of the  Kac Moment Formula used in this paper. 

We anticipate that   the   estimates in Section \ref{sec-indf}, of properties of the probability densities of a wide class of L\'evy processes, will be useful in the study  of  L\'evy processes and their local times.

  \section{Proof of Theorem \ref{theo-clt2r}}\label{sec-CLT}
  
The   proof of Theorem \ref{theo-clt2r} is long and difficult.  In order to make it easier to follow we first present   the main steps of the proof heuristically.  We then restate them, precisely, in a series of lemmas and show how Theorem \ref{theo-clt2r} follows from these lemmas.   These lemmas are proved in Sections 
\ref{sec-indf}  through \ref{sec-9}.

Let
\begin{equation}
I_{j,k,t}:=\int ( L^{ x+1}_{t}- L^{ x}_{t})\circ\th_{jt}\,\, (   L^{ x+1}_{t}-  L^{ x}_{ t})\circ  \th_{kt} \,\, \,dx\label{m.1}
\end{equation}
and
\begin{equation}
  \al_{j,k,t  } := \int  L^{ x}_{t }\circ \th_{jt}\,\,  L^{ x}_{t}\circ  \th_{kt}\,\,  \,dx.\label{m.1a}
\end{equation} 
(An integral sign without limits is to be read as $\int_{-\ff}^{\ff}\,$.)

For any integer $l$ set   
\begin{equation}
\wt I_{l,t} :=
 \sum_{\stackrel{j,k=0}{j< k}}^{l-1} I_{j,k,t/l}.\label{e1.6}
\end{equation}
Using the additivity property of   local times we can write
\begin{equation}
L^{ x}_{ t} =\sum_{j=0}^{l-1}L^{ x}_{t/l}\circ \th_{jt/l}\label{e1.3},
\end{equation}
so that 
\begin{equation}
 \int ( L^{ x+1}_{t}- L^{ x}_{ t})^{2} \, \,dx =
 \sum_{j,k=0}^{l-1} I_{j,k,t/l}=2\wt I_{l,t}+\sum_{ j=0 }^{l-1}I_{j,j,t/l}.\label{e1.4}
\end{equation}
Consequently 
\begin{eqnarray}
 &&  \int ( L^{ x+1}_{t}- L^{ x}_{ t})^{ 2}\,dx- E\(\int ( L^{ x+1}_{t}- L^{ x}_{ t})^{ 2}\,dx\)
\label{e2.1m}\\
 &  &\qquad=2\wt I_{l,t}+   \lc \sum_{ j=0 }^{l-1}   I_{j,j,t/l}      -  E\(\int ( L^{ x+1}_{t}- L^{ x}_{ t})^{ 2}\,dx\)\rc.\nonumber
\end{eqnarray}
Similarly we set 
\begin{equation}
\wt \al_{l,t } = \sum_{\stackrel{j,k=0}{j< k}}^{l-1}\al_{j,k,t/l}
,\label{e1.7}
\end{equation}
and  write
\begin{equation}
 \al_{t}:= \int  \(L^{ x}_{t} \)^{2}  \,  \,dx=  \sum_{j,k=0}^{l-1}\al_{j,k,t/l}=2\wt \al_{l,t }+\sum_{ j=0 }^{l-1}\al_{j,j,t/l}.\label{e1.5}
\end{equation}

The main steps in the proof of Theorem \ref{theo-clt2r} are to show that:
\begin{enumerate}
\item 
 The `off-diagonal' terms  $\wt I_{l,t}$ and $ \sqrt{\wt \al_{l,t }} $ are comparable asymptoticly as $t\to\ff$.
\item
 The diagonal term  $\sum_{ j=0 }^{l-1}\al_{j,j,t/l}$ is   negligible, as $t\to\ff$, compared to the terms in 1.
\item The diagonal term  $\sum_{ j=0 }^{l-1}I_{j,j,t/l}$ is such that
\begin{equation}
 \sum_{ j=0 }^{l-1}I_{j,j,t/l}- E\(\int ( L^{ x+1}_{t}- L^{ x}_{ t})^{ 2}\,dx\)
 \ee
 \label{hue.3}
 is   negligible, as $t\to\ff$, compared to the terms in 1.
\end{enumerate}
We now explain the precise meaning of these statements, and show how they imply Theorem \ref{theo-clt2r}.

\medskip	The precise meaning of step 1.  is given by the following lemma. (Lemmas \ref{lem-hue.1}--\ref{lem-hue.3} are proved in Section \ref{sec-5}.)
  \begin{lemma} \label{lem-hue.1} Under the hypotheses of Theorem \ref{theo-clt2r},  for each $m$, with $l=l(t)=[\log t]^{q}$, for any $q>0$,  
\be  
 E\(\( \wt I_{l,t} \)^{m}\)\nn 
   =\left\{\begin{array}{ll}\label{e1.8ax}
 \displaystyle {( 2n)!\over 2^{ n}n!}\(4c_{\psi,1}\)^{ n} E\lc\(\wt  \al_{l,t} \)^{
n}\rc+o((t^{2}\psi^{-1} (1/t))^{n} ) &\!\! m=2n\\\\ 
O((t^{2}\psi^{-1} (1/t))^{m/2}t^{-\ep})  &\!\!  \mbox{otherwise.}
\end{array}
\right. 
\ee 
 \end{lemma}
 This lemma is the crux of this paper. We note that even though
  the summands $I_{j,k,t/l}$ of $\wt I_{l,t}$ are not independent, the fact that $j\neq k$ 
 provides enough structure for a proof. 
  
\medskip	The precise meaning of step 2. is given by the next lemma. 
  \begin{lemma} \label{lem-hue.2} Under the hypotheses of Theorem \ref{theo-clt2r},  for each $n$, with $l=l(t)=[\log t]^{q}$, for any $q>0$, 
\begin{equation}
\lim_{t\to\ff}   \frac{|E(2\wt  \al_{l,t})^{n}-E(  \al_{ t})^{n}|}{\(t^{2}\psi^{-1}(1/t)\)^{n}}=0.\label{2.11}
   \end{equation}
 \end{lemma} 
 
  Lastly, the precise meaning of step 3. is given by the next lemma.
  
  \begin{lemma} \label{lem-hue.3} Under the hypotheses of Theorem \ref{theo-clt2r},   with $l=l(t)=[\log t]^{q}$, for  $q$ sufficiently large, 
\begin{equation}
 \lim_{t\to\ff}{  \sum_{ j=0 }^{l-1}    \( I_{j,j,t/l} -  E\(\int ( L^{ x+1}_{t}- L^{ x}_{ t})^{ 2}\,dx\)\)\over t(\psi^{-1}(1/t))^{1/2}}=0.
   \end{equation}
   in $L^{2}$.
 \end{lemma}

   We also need to know the limiting behavior of the moments of $\al_{t}$.  It is given by   
the  next lemma which is proved in Section \ref{sec-ilt}.
   \bl\label{lem-ilt}Under the hypotheses of Theorem \ref{theo-clt2r},  for each $n$,  
   \begin{equation}
   \lim_{t\rar \ff}  E\lc\({\al_{t} \over t^{2}\psi^{-1} (1/t)}\)^{
n}\rc=E\lc\(\al_{\bb, 1}\)^{
n}\rc \label{e1.57c}.
   \end{equation}
   \el
 
 \medskip	\noindent{\bf Proof of   Theorem \ref{theo-clt2r} }
In (\ref{e1.8ax}) replace $ \wt I_{l,t}$ by $2 \wt I_{l,t}$ and $\(4c_{\psi,1}\)^{ n} E\lc\(\wt  \al_{l,t} \)^{n}\rc$ by $\(8c_{\psi,1}\)^{ n} E\lc\(2\wt  \al_{l,t} \)^{n}\rc$. Then set
\begin{equation}
   E\lc\(2\wt  \al_{l,t} \)^{n}\rc=E(\al_{t})^{n} +E\lc\(2\wt  \al_{l,t} \)^{n}\rc -E(\al_{t})^{n}. 
   \end{equation}
Then use Lemmas \ref{lem-hue.2}  and \ref{lem-ilt} to see that 
  for each integer $m$ 
\begin{eqnarray} &&
\lim_{ t\rar \ff}E\(\(  {2\wt I_{l,t} \over t\sqrt{\psi^{-1} (1/t)}}   \)^{m}\)\nn\\ 
&&\hspace{ 1in}  =\left\{\begin{array}{ll} \displaystyle
{( 2n)!\over 2^{ n}n!}\( 8c_{\psi,1}\)^{ n} E\lc\(\al_{\bb, 1}\)^{
n}\rc &\mbox{ if }m=2n\\\\
0&\mbox{ otherwise.}
\end{array}
\right.
\label{e1.8ajj}
\end{eqnarray}

Note that the right-hand side of (\ref{e1.8ajj})  is the $2n$--th moment of  $ \(8c_{\psi,1}\)^{ 1/2}\newline	\sqrt{\al_{\bb, 1}}\,\eta$ when $\al_{\bb, 1}$ and $\eta$ are independent. Furthermore,
it follows from  \cite[(6.12)]{CLRR}    that 
 \begin{equation}
 E \(\al_{\bb, 1}\)^{ n} \leq C^{ n}( (2n)!)^{ 1/(2\bb)}.\label{rb.1}
 \end{equation}
Consequently, since $\sqrt{(2n)!}\le 2^{n}n!$
\begin{equation}
   E\( \(8c_{\psi,1}\)^{ 1/2}\sqrt{\al_{\bb, 1}}\,\eta\)^{m}\le C^{m}( m!)^{ (\bb+1)/(2\bb )}.
   \end{equation}
 This implies that   $ \(8c_{\psi,1}\)^{ 1/2}  \sqrt{\al_{\bb, 1}}\,\eta$
  is determined by its moments; (see \cite[p. 227-228]{Feller}).
    Therefore, by the method of moments, \cite[Theorem 30.2]{B}\label{meth}), it follows from  (\ref{e1.8ajj})
     that
 \begin{eqnarray}
&& \lim_{t\rar \ff}  {2\wt I_{l,t} \over t\sqrt{\psi^{-1} (1/t)}}  
\label{e1.8}\stackrel{\mathcal{L}}{\Longrightarrow}\(8c_{\psi,1}\)^{ 1/2}\sqrt{\al_{\bb,1}}\,\,\eta.
\end{eqnarray}
  Theorem \ref{theo-clt2r} then follows from (\ref{e2.1m}), (\ref{e1.8}) and  Lemma   \ref{lem-hue.3}.

 \begin{remark} \label{rem-2.1}
 {\rm
 For Brownian motion $\psi(p)=p^{2}/2$. We have
  \bea
c_{p^{2}/2,0}&  = &  {4  \over \pi}\int_{0}^{\ff} { \sin^{2} (p /2)  \over p^{2}}\,dp=1;\\
c_{p^{2}/2,1}&=&{64  \over \pi}\int_{0}^{\ff} { \sin^{4} p  /2  \over p ^{4}  }\,dp=\frac{8}{3}.
   \eea
 Also
$
 L_{2,1}^{x}\stl  (1/2) L_{2 }^{x}\stl (1/\sqrt2) L_{1 }^{x/\sqrt2}
$
 where $\{L_{t}^{x}\}$ denotes the local time of Brownian motion, so that
 \begin{equation}
   \(\int _{-\ff}^{\ff} \(L_{2,1}^{x}\)^{2}\,dx\)^{1/2}=\frac{1}{2^{1/4}}   \(\int _{-\ff}^{\ff} \(L_{ 1}^{x}\)^{2}\,dx\)^{1/2}.
   \end{equation}
   Since  $\sqrt{\psi^{-1}(1/t)}=  2^{1/4}/t^{1/4} $, we see that (\ref{r5.0tweaks}) and (\ref{1.11}) imply (\ref{5.0tweak}).
 }\end{remark}

 \section{Estimates for the probability densities of certain L\'evy processes  }\label{sec-indf}

 Let $p_{s}(x)$ denote the density of the symmetric  L\'{e}vy process $X$ with L\'{e}vy exponent $\psi (\la)$ as described in (\ref{st.1}).  Let $\De_{ x}^{\ga}$  denote
the finite difference operator on the variable $x$, i.e.
\begin{equation}
\De_{x}^{ \ga}\,f(x)=f(x+\ga)-f(x).\label{pot.3w}
\end{equation}
We write  $\De^{ \ga}$ for $\De_{x}^{ \ga}$ when the variable $x$ is clear.

Let
 \bea
u(x,t)&:=& \int_{0}^{t}\, \,p_{s }(x)\,ds \\
 v_{\ga}(x,t)&:=&\int_{0}^{t}\, |\De ^{\ga}\,p_{s }(x)|\,ds\label{2.3} \\
w_{\ga}(x,t)&:=&\int_{0}^{t}\,|\De^{\ga}\De ^{-\ga} \,p_{s }(x)|\,ds 
\eea
We also write $ v(x,t)$ for  $v_{1}(x,t)$ and   $ w(x,t)$ for  $w_{1}(x,t)$.

\medskip	
The following lemmas   provide the main estimates we use in this paper. Their  proofs are given in Section \ref{sec-lemproofs}.

\begin{lemma}\label{lem-vproprvtj} Under the hypotheses of 
Theorem  \ref{theo-clt2r},  for all $t$ sufficiently large  
    \bea
\sup_{x\in R^{1}}u(x,t)&\le&  Ct \psi^{-1}(1/t)  \label{jr.1j};\\
 \sup_{x\in R^{1}}v(x,t)&\le& C\log t   \label{2.4j};\\
 \sup_{x\in R^{1}}w (x,t)&\le&  C    \label{jrst1.3yj},
 \eea 
and
 \bea
   \int u(x,t)\,dx&= & t\label{ee1};\\
   \int v (x,t)\,dx&\le&  C   t\(\psi^{-1}(1/  t)\)\log t\label{151ee1};\\
       \int w(x,t)\,dx&\le& C (\log t)^{2} \label{132ee1}; \\ 
    \int w^{2} (x,t)\,dx&\le& C \log t \label{132ee2}; \\
     \int_{|x|\ge u} w^{2} (x,t)\,dx&\le& C  {(\log t)^{2}\over u }. \label{132ee3} 
\eea
  
\end{lemma}

  \begin{lemma}\label{lem-2.4}  Under the hypotheses of 
Theorem  \ref{theo-clt2r},  for all $t$ sufficiently large
  \begin{equation}
   |\De^{1}p_{t}(0)|\le C \(\psi^{-1}(1/  t)\)^{3},\label{mac.1}
   \end{equation}
and  
\begin{equation}
\int_{0}^{2t}\int_{0}^{2t} |\De^{1}p_{r+s }(0)|\,dr\,ds \leq C\(t^{2}\(\psi^{-1}(1/t)\)^{3}+L(t)+ 1\),
\label{wo.4a}
 \end{equation}
 where $L(t)$ is a slowly varying function at infinity.
 \end{lemma}

\begin{lemma}\label{lem-2.5}  Under the hypotheses of 
Theorem  \ref{theo-clt2r},

 \bea 
  \int_{0}^{\ff}\De^{ 1} \,\,p_{s }(0)\,ds & =& - \, c_{\psi, 0};\label{2.13}\\
\int  \(\int_{0}^{\ff}\De^{ 1}\De^{ -1}\,\,p_{s }(x)\,ds\)^{2}\,dx&=&   c_{\psi, 1};\label{h3.1}
\eea
  and   
\begin{equation} 
\int  \(\,\int_{0}^{t }\,\De^{ 1}\De^{ -1} \,p_{s }(x)\,ds\)^{2}\,dx=   c_{\psi, 1} +O(t^{-1/3}),\label{h3.2}
\end{equation}
as $t\to\ff$.
\end{lemma}

\section{Moments of  increments of local times.}\label{sec-3}

  The goal of this section is to establish the following lemma, in which    $I_{j,k,t}$ and $ \al_{j,k,t  } $ are  defined in (\ref{m.1}) and (\ref{m.1a}).

\bl\label{lem-multiple}
 Let 
$m_{j,k}$, $0\leq j<k \leq K$ be positive integers with  $\sum_{ j,k=0, j<k }^{K}\newline  m_{j,k}=m$. 
 If all the   integers $m_{j,k}$ are even, then  for some $\ep>0$
\begin{eqnarray} && 
E\(\prod_{\stackrel{j,k=0}{j< k}}^{K} \(I_{j,k,t} \)^{m_{j,k}}\)
\label{m.2}\\ 
&&\qquad
=  \prod_{\stackrel{j,k=0}{j< k}}^{K}{( 2n_{j,k})!\over 2^{ n_{j,k}}(n_{j,k}!)}\( 4c_{\psi ,1}\)^{ n_{j,k}} E\(\prod_{\stackrel{j,k=0}{j< k}}^{K}\(\al_{j,k,t}\)^{
n_{j,k}}\) +O\(t^{(2-1/\bb)n-\ep}\),
\nn
\end{eqnarray}
  where  $n_{j,k}=m_{j,k}/2$.  

If any of the $m_{j,k}$
are odd, then
\begin{equation}
E\(\prod_{\stackrel{j,k=0}{j< k}}^{K} \(I_{j,k,t} \)^{m_{j,k}}\)=   O\(t^{(2-1/\bb)\frac{m}{2}-\ep}\).\label{m.2a}
   \end{equation} 
   In (\ref{m.2}) and (\ref{m.2a}) the error terms  may depend on $m$, but not on the 
individual terms $m_{j,k} $.  
\el
 
 \Proof 
We can write
\begin{eqnarray} \lefteqn{
E\(\prod_{\stackrel{j,k=0}{j< k}}^{K} \(I_{j,k,t} \)^{m_{j,k}}\)\label{m.3}}\\
&&=E\(\prod_{\stackrel{j,k=0}{j< k}}^{K} \prod_{i=1}^{m_{j,k}}\(\int ( \De^{1}_{x_{j,k,i}} L^{ x_{j,k,i}}_{t}\circ\th_{jt})\,\, (  \De^{1}_{x_{j,k,i}}L^{ x_{j,k,i}}_{ t}\circ  \th_{kt}) \,\, \,dx_{j,k,i} \) \)\nn\\
&&=\int  \lc \prod_{\stackrel{j,k=0}{j< k}}^{K} \prod_{i=1}^{m_{j,k}}\De^{1,j}_{x_{j,k,i}}  \De^{1,k}_{x_{j,k,i}} \rc E\(\prod_{\stackrel{j,k=0}{j< k}}^{K} \prod_{i=1}^{m_{j,k}}\(  (  L^{ x_{j,k,i}}_{t}\circ\th_{jt})\,\, (  L^{ x_{j,k,i}}_{ t}\circ  \th_{kt})  \) \)\nonumber\\
&&\hspace{3.5 in}\prod_{\stackrel{j,k=0}{j< k}}^{K} \prod_{i=1}^{m_{j,k}}\,dx_{j,k,i}\nn,
\end{eqnarray}
where the notation $\De^{1,j}_{x_{j,k,i}}$   indicates that we apply the difference operator $\De^{1}_{x_{j,k,i}}$ to    $ L^{ x_{j,k,i}}_{t}\circ\th_{jt}$.  Note that there are $2m$ applications of the difference operator $\De^{1}$.   

 Consider 
\begin{equation}
   E\(\prod_{\stackrel{j,k=0}{j< k}}^{K} \prod_{i=1}^{m_{j,k}}\(  (  L^{ x_{j,k,i}}_{t}\circ\th_{jt})\,\, (  L^{ x_{j,k,i}}_{ t}\circ  \th_{kt})  \) \).
   \end{equation}
 We collect  all the factors containing $\th_{lt}$ and write
\begin{eqnarray}
&&E\(\prod_{\stackrel{j,k=0}{j< k}}^{K} \prod_{i=1}^{m_{j,k}}\(  (  L^{ x_{j,k,i}}_{t}\circ\th_{jt})\,\, (  L^{ x_{j,k,i}}_{ t}\circ  \th_{kt})  \) \)
\label{reo.1}\\
&&\qquad=E\(\prod_{l=0}^{K}\lc \(\prod_{j=0}^{l-1} \prod_{i=1}^{m_{j,l}}   L^{ x_{j,l,i}}_{t}\) 
 \(\prod_{k=l+1}^{K} \prod_{i=1}^{m_{l,k}}L^{ x_{l,k,i}}_{t}\)   \rc \circ  \th_{lt}\)   \nonumber\\
&&\qquad=E\(\prod_{l=0}^{K}H_{l} \circ  \th_{lt}\)   \nonumber,
\end{eqnarray}
where
\begin{equation}
H_{l}=\(\prod_{j=0}^{l-1} \prod_{i=1}^{m_{j,l}}    L^{ x_{j,l,i}}_{t}\) 
 \(\prod_{k=l+1}^{K} \prod_{i=1}^{m_{l,k}}    L^{ x_{l,k,i}}_{t}\).\label{reo.2}
\end{equation}
By the Markov property
\begin{equation}
E\(\prod_{l=0}^{K}H_{l} \circ  \th_{lt}\)=E\(H_{0} \,E^{X_{t}}\(\prod_{l=1}^{K}H_{l} \circ  \th_{(l-1)t}\)\)  .\label{reo.3}
\end{equation}
Let 
\be
m_{l}=\sum_{k=l+1}^{K}\,m_{l,k}+\sum_{j=0}^{l-1}\,m_{j,l},\qquad l=0,\ldots,K-1,\label{m4.9}
\ee
and note that  $m_{l}$ is the number of local time factors in $H_{l}$.  

Let 
\begin{equation}
   f(y)=E^{y}\(\prod_{l=1}^{K}H_{l} \circ  \th_{(l-1)t}\).
   \end{equation}
   It follows from  Kac's Moment Formula, Theorem \ref{KMT},  that  for any $z\in R^{1}$
\begin{eqnarray}
\lefteqn{E^{z}\(\prod_{l=0}^{K}H_{l} \circ  \th_{lt}\)}\\
&=& E^{z}\(H_{0} \,f(X_{t})\)\nn
\label{reo.4} \\
&=&   \sum_{\pi_{0}}\int_{\{\sum_{q=1}^{m_{0}}r_{0,q}\leq t\}}
p_{r_{0,1}}(x_{\pi_{0}(1)}-z) \prod_{q=2}^{m_{0}}p_{r_{0,q}}(x_{\pi_{0}(q)}-x_{\pi_{0}(q-1)}) \nonumber\\
&&\hspace{1 in}\(\int p_{(t-\sum_{q=1}^{m_{0}}r_{0,q}) }(y-x_{\pi_{0}(m_{0})})f(y)\,dy\)\prod_{q=1}^{m_{0}}\,dr_{0,q} \nonumber,
\end{eqnarray}
where the sum runs over all 
  bijections $\pi_{0} $ from $[1,m_{0}]$ to
\begin{equation}
I_{0}=\bigcup_{k=1}^{K}\{(0,k,i),\,1\leq i\leq m_{0,k}\}.\label{reo.5}
\end{equation}
Clearly, $I_{0}$ is the set of subscripts  of  the terms $x_{\cd} $ appearing in the local time factors in $H_{0}$. 

By the Markov property 
\bea
   f(y)&=&E^{y}\(H_{1}E^{X_{2t}}\(\prod_{l=2}^{K}H_{l} \circ  \th_{(l-2)t}\)\)\label{4.13}\\
   &:=&E^{y}\(H_{1}g(X_{2t}) \).\nn
   \eea
 Therefore, by (\ref{reo.3})--(\ref{4.13}), for any $z'\in R^{1}$ 
 \bea
\lefteqn{E^{z'}\(\prod_{l=0}^{K}H_{l} \circ  \th_{lt}\)\label{reo.6}}\\
&&=E^{z'}\(H_{0}E^{X_{t}}\(H_{1} \,g(X_{2t})\)\)
\nn\\
&&=  \sum_{\pi_{0}}\int_{\{\sum_{q=1}^{m_{0}}r_{0,q}\leq t\}}
p_{r_{0,1}}(x_{\pi_{0}(1)}-z') \prod_{q=2}^{m_{0}}p_{r_{0,q}}(x_{\pi_{0}(q)}-x_{\pi_{0}(q-1)}) \nonumber\\
&&\hspace{.5 in}\(\int p_{(t-\sum_{q=1}^{m_{0}}r_{0,q}) }(y-x_{\pi_{0}(m_{0})})
E^{y}\(H_{1} \,g(X_{2t})\)\,dy\)\prod_{q=1}^{m_{0}}\,dr_{0,q} \nonumber\\
&&=  \sum_{\pi_{0}}\int_{\{\sum_{q=1}^{m_{0}}r_{0,q}\leq t\}}
p_{r_{0,1}}(x_{\pi_{0}(1)}-z') \prod_{q=2}^{m_{0}}p_{r_{0,q}}(x_{\pi_{0}(q)}-x_{\pi_{0}(q-1)}) \nonumber\\
&&\hspace{.5 in}  p_{(t-\sum_{q=1}^{m_{0}}r_{0,q}) }(y-x_{\pi_{0}(m_{0})})\nn\\
&&\hspace{.3 in} \sum_{\pi_{1}}\int_{\{\sum_{q=1}^{m_{1}}r_{1,q}\leq t\}}
p_{r_{1,1}}(x_{\pi_{1}(1)}-y) \prod_{q=2}^{m_{1}}p_{r_{1,q}}(x_{\pi_{1}(q)}-x_{\pi_{1}(q-1)}) \nonumber\\
&&\hspace{.3 in}\(\int p_{(t-\sum_{q=1}^{m_{1}}r_{1,q}) }(y'-x_{\pi_{1}(m_{1})})g(y')\,dy'\)\prod_{q=1}^{m_{1}}\,dr_{1,q} \nonumber
\,dy\prod_{q=1}^{m_{0}}\,dr_{0,q} \nonumber
\end{eqnarray}
where the second sum   runs over all 
  bijections $\pi_{1} $ from $[1,m_{1}]$ to
\begin{equation}
I_{1}=\{(0,1,i),\,1\leq i\leq m_{0,1}\}\bigcup_{k=2}^{K}\{(1,k,i),\,1\leq i\leq m_{1,k}\}\label{reo.7}
\end{equation}

As above, $I_{1}$ is the set of subscripts  of  the terms $x_{\cd} $ appearing in the   local time factors in $H_{1}$. 

  We now use the Chapman-Kolmogorov equation to integrate with respect to $y$ to get
\begin{eqnarray}
\lefteqn{E^{z'}\(H_{0}E^{X_{t}}\(H_{1} \,g(X_{t})\)\)
\label{reo.7a}}\\
&&=  \sum_{\pi_{0},\pi_{1}}\int_{\{\sum_{q=1}^{m_{0}}r_{0,q}\leq t\}}
p_{r_{0,1}}(x_{\pi_{0}(1)}-z') \prod_{q=2}^{m_{0}}p_{r_{0,q}}(x_{\pi_{0}(q)}-x_{\pi_{0}(q-1)}) \nonumber\\ 
&&\hspace{.3 in}  \int_{\{\sum_{q=1}^{m_{1}}r_{1,q}\leq t\}}
p_{(t-\sum_{q=1}^{m_{0}}r_{0,q})+r_{1,1}}(x_{\pi_{1}(1)}-x_{\pi_{0}(m_{0})}) \nonumber\\
&&\hspace{2 in}\prod_{q=2}^{m_{1}}p_{r_{1,q}}(x_{\pi_{1}(q)}-x_{\pi_{1}(q-1)}) \nonumber\\
&&\hspace{.3 in}\(\int p_{(t-\sum_{q=1}^{m_{1}}r_{1,q}) }(y'-x_{\pi_{1}(m_{1})})g(y')\,dy'\)\prod_{q=1}^{m_{1}}\,dr_{1,q} \nonumber
 \prod_{q=1}^{m_{0}}\,dr_{0,q}. \nonumber
\end{eqnarray}

Iterating this procedure, and recalling (\ref{reo.1}) we see that
\begin{eqnarray}
\lefteqn{E\(\prod_{\stackrel{j,k=0}{j< k}}^{K} \prod_{i=1}^{m_{j,k}}\( (  L^{ x_{j,k,i}}_{t}\circ\th_{jt})\,\, ( L^{ x_{j,k,i}}_{ t}\circ  \th_{kt})   \) \)
\label{m.4}}\\
&&=  \sum_{\pi_{0},\ldots, \pi_{K}}  \prod_{l=0}^{K}\int_{\{\sum_{q=1}^{m_{l}}r_{l,q}\leq t\}}
p_{(t-\sum_{q=1}^{m_{l-1}}r_{l-1,q})+r_{l,1}}(x_{\pi_{l}(1)}-x_{\pi_{l-1}(m_{l-1})})  \nonumber\\
&&\hspace{1.8 in}\prod_{q=2}^{m_{l}}p_{r_{l,q}}(x_{\pi_{l}(q)}-x_{\pi_{l}(q-1)})\prod_{q=1}^{m_{l}}\,dr_{l,q},\nonumber
\end{eqnarray}
  where $\pi_{-1}(m_{-1}):=0$ and   $1-\sum_{q=1}^{m_{ -1}}r_{ -1,q}:=0$. In (\ref{m.4})
  the sum runs over all 
$\pi_{0},\ldots, \pi_{K}$ such that each $\pi_{l}$ is a bijection from $[1,m_{l}]$ to  
\begin{equation}
I_{l}=\bigcup_{j=0}^{l-1}\{(j,l,i),\,1\leq i\leq m_{j,l}\}\bigcup_{k=l+1}^{K}\{(l,k,i),\,1\leq i\leq m_{l,k}\}.\label{m.5}
\end{equation}
As in the observations about $I_{0}$ and $I_{1}$, we see  that  $I_{l}$ is the set of subscripts of  the terms $x_{\cd} $ terms appearing in the   local time factors in   $H_{l}$. Since there are   $2m$ local time factors  we have that   $\sum_{l=0}^{K}m_{l}=2m$.

\medskip	We  now use (\ref{m.4}) in (\ref{m.3}) and continue to develop an expression for the left-hand side of   (\ref{m.3}). Let $\mathcal{B}$ to denote the set of  $(K+1)$--tuples, $\pi=(\pi_{0},\ldots, \pi_{K})$, of bijections described in (\ref{m.5}). Clearly 
\be
|\mathcal{B}|=\prod_{l=0}^{K}m_{l}!\leq (2m)!.
\ee
  Also, similarly to the way we obtain the first equality in (\ref{reo.1}), we see that  
\begin{equation}
\prod_{\stackrel{j,k=0}{j< k}}^{K} \prod_{i=1}^{m_{j,k}}\De^{1,j}_{x_{j,k,i}}  \De^{1,k}_{x_{j,k,i}} = \prod_{l=0}^{K}\prod_{q=1}^{m_{l}}\De^{1,l}_{x_{\pi_{l}(q)}}\label{m.5a}.
\end{equation}
Consequently   
\be  
E\(\prod_{\stackrel{j,k=0}{j< k}}^{K} \(I_{j,k,t} \)^{m_{j,k}}\)
 = \sum_{\pi_{0},\ldots, \pi_{K}}\int \wt\mathcal{T}_{t}( x;\,\pi  )\prod_{j,k,i}\,dx_{j,k,i} \label{m.6a}
\ee  
where  we take the product over  $ \{  0\leq j<k\leq K,\, 1\leq i\leq m_{j,k}\}$,   $\pi\in \mathcal{B}$ 
and
\begin{eqnarray}
\lefteqn{ \wt\mathcal{T}_{t}( x;\,\pi  )
\label{m.4a}}\\
&&=    \prod_{l=0}^{K}\prod_{q=1}^{m_{l}}\De^{1}_{x_{\pi_{l}(q)}}\int_{\{\sum_{q=1}^{m_{l}}r_{l,q}\leq t\}}
p_{(t-\sum_{q=1}^{m_{l-1}}r_{l-1,q})+r_{l,1}}(x_{\pi_{l}(1)}-x_{\pi_{l-1}(m_{l-1})})  \nonumber\\
&&\hspace{1.8 in}\prod_{q=2}^{m_{l}}p_{r_{l,q}}(x_{\pi_{l}(q)}-x_{\pi_{l}(q-1)})\prod_{q=1}^{m_{l}}\,dr_{l,q}.\nonumber
\end{eqnarray}

\medskip	We continue to rewrite the right-hand side of   (\ref{m.6a}).

 In (\ref{m.4a}),   each   difference operators, say $\De^{1}_{u} $ is applied to the product of two terms, say $p_{\cd}(u-a)\,p_{\cd}(u-b)$, using the product rule for difference operators  we see that
\bea 
   && \De^{1}_{u}\{p_{\cd}(u-a)p_{\cd}(u-b)\}\label{3.24q}\\&&\qquad= \De^{1}_{u}\,p_{\cd}(u-a) p_{\cd}(u+1-b)+p_{\cd}(u-a) \De^{1}_{u}\, p_{\cd}(u-b)\nn 
 \eea
Consider an example of how the term $\De^{1}_{a}\De^{1}_{u}p_{\cd}(u-a)$ may appear. It could be  by the application
\begin{equation}
  \De^{1}_{a}  \(\De^{1}_{u}\,p_{\cd}(u-a) p_{\cd}(v-a)\),\label{3.25}
   \end{equation}
in which we take account of the two terms to which $\De^{1}_{a} $ is applied. Using the product rule in (\ref{3.24q}) we see that (\ref{3.25})
\begin{equation}
   =  \(\De^{1}_{a}   \De^{1}_{u}\, p_{\cd}(u-a) \)p_{\cd}(v -(a+1))+ \De^{1}_{u}\,p_{\cd}(u-a) \De^{1}_{a} p_{\cd}(v-a).\label{3.26}
   \end{equation}
Consider one more example
\bea 
  &&  \De^{1}_{a}  \(\De^{1}_{u}\,p_{\cd}(u-a)\De^{1}_{v}\, p_{\cd}(v-a)\)\label{3.27}\\
  &&\qquad\nn = \(\De^{1}_{a}   \De^{1}_{u}\, p_{\cd}(u-a) \)\De^{1}_{v}\, p_{\cd}(v -(a+1))\nn\\
  &&\hspace{1in}\nn + \De^{1}_{u}\,p_{\cd}(u-a) \De^{1}_{a}\De^{1}_{v} p_{\cd}(v-a).
 \eea
Note that in both examples the arguments of probability densities  with two difference operators applied to it does not contain a $1$. This is true in general because the difference 
formula, (\ref{3.24q}), does not add a  $  1$ to the argument of a term to which a difference operator is applied. Otherwise we may have a $\pm 1$ added to the arguments of probability densities to which one  difference operator is applied, as in (\ref{3.27}), or to the arguments of probability densities to which no difference operator is applied, as in (\ref{3.26}).

  Based on the argument of the preceding paragraph we write (\ref{m.4a}) in the form
\be  
 E\(\prod_{\stackrel{j,k=0}{j< k}}^{K} \(I_{j,k,t} \)^{m_{j,k}}\)=  \sum_{a}\sum_{\pi_{0},\ldots, \pi_{K}}\int \mathcal{T}'_{t}( x;\,\pi ,a )\prod_{j,k,i}\,dx_{j,k,i},\label{4.25}
\ee  
where  
\begin{eqnarray}
\lefteqn{\mathcal{T}'_{t}( x;\,\pi ,a )=   \prod_{l=0}^{K}\int_{\RR_{l}}
\(\(\De^{ 1}_{ x_{\pi_{l}(1)}}\)^{a_{ 1}(l,1)}
\(\De^{ 1}_{ x_{\pi_{l-1}(m_{l-1})}}\)^{a_{ 2}(l,1)}\right.
\label{m.7}}\\
&&\hspace{1.5 in}
\left.\,p^{\sharp}_{(t-\sum_{q=1}^{m_{l-1}}r_{l-1,q})+r_{l,1}}(x_{\pi_{l}(1)}-x_{\pi_{l-1}(m_{l-1})})\)  \nonumber\\
&&\hspace{.1 in}\prod_{q=2}^{m_{l}}\(\(\De^{ 1}_{ x_{ \pi_{l}( q)}}\)^{a_{ 1}(l,q)}
\(\De^{ 1}_{ x_{ \pi_{l}( q-1)}}\)^{a_{ 2}(l,q)}\,p^{\sharp}_{r_{l,q}}(x_{\pi_{l}(q)}-x_{\pi_{l}(q-1)})\)
\prod_{q=1}^{m_{l}}\,dr_{l,q}.\nonumber
\end{eqnarray}

In (\ref{m.7})
  $\RR_{l}=\{\sum_{q=1}^{m_{l}}r_{l,q}\leq t\}$. In (\ref{4.25})
   the first sum is taken  over    all  
\be
a=(a_{ 1},a_{ 2})\,:\,\{(l,q),\, 0\leq l\leq K,\, 1\leq q\leq m_{l}\}\mapsto \{ 0,1\}\times \{ 0,1\}
\ee
  with the
restriction that for each triple $j,k,i$, there are exactly two  factors  of the form $\De^{
1}_{ x_{j,k,i}}$, each of which is applied to one of the  terms  $p^{\sharp}_{r_{\cdot}}(\cdot)$
that contains $x_{j,k,i}$ in its argument. This   condition can be stated more formally by saying that 
  for each $l$ and $q=1,\ldots, m_{l}-1$, if $\pi_{l}(q)=(j,k,i)$, then   $\{a_{ 1}(l,q), a_{ 2}(l,q+1) \}=\{0,1\}$ and if $q=m_{l}$ then   $\{a_{ 1}(l,m_{l}), a_{ 2}(l+1,1) \}=\{0,1\}$.    (Note that  when we write $\{a_{ 1}(l,q), a_{ 2}(l,q+1) \}=\{0,1\}$ we mean as two  sets, so, according to what $a$ is, we may have  $a_{ 1}(l,q)=1$ and $ a_{ 2}(l,q+1)=0$ or $a_{ 1}(l,q)=0$ and $a_{ 2}(l,q+1)=1$ and similarly for $\{a_{ 1}(l,m_{l}), a_{ 2}(l+1,1) \}$.)
  Also, in (\ref{m.7}) we define $(\De_{x_{i}}^{1})^{0}=1 $ and $(\De_{0}^{1}) =1 $.

  In (\ref{m.7}),  because of  (\ref{3.24q}),   $p^{\sharp}_{r_{\cdot}}(z)$ can take any of the values $p_{r_{\cdot}}(z)$, $p_{r_{\cdot}}(z+1)$ or $p_{r_{\cdot}}(z-1)$. (We   must  consider all three possibilities.) Finally, it is important to emphasize that in (\ref{m.7}) each of  the difference operators is applied to only one of the  terms  $p^{\sharp}_{r_{\cdot}}(\cdot)$.

\medskip	 Rather than (\ref{4.25}), we first analyze 
\begin{equation}
\sum_{a}\sum_{\pi_{0},\ldots, \pi_{K}}\int \mathcal{T}_{t}( x;\,\pi ,a )\prod_{j,k,i}\,dx_{j,k,i}\label{m.8},
\end{equation}
where  
\begin{eqnarray}
&&\mathcal{T}_{t}( x;\,\pi ,a )=   \prod_{l=0}^{K}\int_{\RR_{l}}
\(\(\De^{ 1}_{ x_{\pi_{l}(1)}}\)^{a_{ 1}(l,1)}
\(\De^{ 1}_{ x_{\pi_{l-1}(m_{l-1})}}\)^{a_{ 2}(l,1)}\right.
\label{m.9}\\
&&\hspace{1.5 in}
\left.\,p_{(t-\sum_{q=1}^{m_{l-1}}r_{l-1,q})+r_{l,1}}(x_{\pi_{l}(1)}-x_{\pi_{l-1}(m_{l-1})})\)  \nonumber\\
&&\hspace{.1 in}\prod_{q=2}^{m_{l}}\(\(\De^{ 1}_{ x_{ \pi_{l}( q)}}\)^{a_{ 1}(l,q)}
\(\De^{ 1}_{ x_{ \pi_{l}( q-1)}}\)^{a_{ 2}(l,q)}\,p_{r_{l,q}}(x_{\pi_{l}(q)}-x_{\pi_{l}(q-1)})\)\prod_{q=1}^{m_{l}}\,dr_{l,q}.\nonumber
\end{eqnarray}
The difference between   $\mathcal{T} _{t}( x;\,\pi ,a )$ and $\mathcal{T}'_{t}( x;\,\pi ,a )$ is that in the former we replace $p^{\sharp}$ by $p$. It is easier to analyze (\ref{m.8}) than (\ref{4.25}). At the conclusion of this proof  we show  that both (\ref{m.8}) than (\ref{4.25}) have the same asymptotic limit as $t$ goes to infinity. 

 \medskip	We first obtain (\ref{m.2}). Let  $m=2n$, since $m_{j,k}=2n_{j,k}$,  $m_{l}=2n_{l}$ for some integer $n_{l}$. (Recall (\ref{m4.9})).  To begin we consider the case in which  $a=e$, where
 \be
 e(l, 2q)=( 1,1)\quad \mbox{and}\quad e(l, 2q-1)=( 0,0) \qquad \forall q.
 \ee

  When $a =e$ we have
\begin{eqnarray}
&&\mathcal{T}_{t}( x;\,\pi ,e )=   \prod_{l=0}^{K}\int_{\RR_{l}}p_{(t-\sum_{q=1}^{m_{l-1}}r_{l-1,q})+r_{l,1}}(x_{\pi_{l}(1)}-x_{\pi_{l-1}(m_{l-1})}) 
\nn\\
&&
\hspace{1.2 in}\prod_{q=2}^{n_{l}} p_{r_{l,2q-1}}(x_{\pi_{l}(2q-1)}-x_{\pi_{l}(2q-2)})    \label{m.10}\\
&&\hspace{1.2 in}\prod_{q=1}^{n_{l}} \De^{ 1} 
\De^{-1} \,p_{r_{l,2q}}(x_{\pi_{l}(2q)}-x_{\pi_{l}(2q-1)})\prod_{q=1}^{m_{l}}\,dr_{l,q}.\nn
\end{eqnarray}
  Here we use the following notation: $\De^{1} p (u-v) =p (u-v+1)-p (u-v)$, i.e., when $\De^{1}$ has no subscript, the difference operator is applied to the whole argument of the function. In this notation,
\begin{equation}
  \De_{u}^{1}  \De_{v}^{1} p (u-v) =  \De^{1}  \De^{-1} p (u-v) .\label{multigraph}
   \end{equation}

  Consider the multigraph $G_{\pi }$ with vertices   
$\{(j,k,i),\,0\leq j<k\leq K,\, 1\leq i\leq m_{j,k}\}$.  Assign an edge 
between the vertices $\pi_{l } (2q-1)$ and $ \pi_{l}(2q)$  for each $ 0\leq l\leq K$ and $1\leq q\leq n_{l}$.  Each vertex is connected to two edges.  To see this suppose that $  \pi_{l}(2q)=\{(j,k,i)\}$,  with $j=l$ and    $ k=l'\neq l$, then there is a unique $q'$ such that   $  \pi_{l'}(2q')$ or $\pi_{l'}(2q'-1)$  is equal to  $ \{(j,k,i)\}$.  Therefore all the vertices lie in some cycle. Assume that there are $S$ cycles. We denote them by $C_{s}$, $s=1,\ldots, S$. 
Clearly, it is possible to have cycles of order two, in which case two vertices are connected by two edges. 

  It is important to note that the graph $G_{\pi}$ does not assign edges between $\pi_{l } (2q)$ and $ \pi_{l}(2q+1)$, although these vertices may be connected by   the edge assigned between $\pi_{l '} (2q'-1)$ and $ \pi_{l'}(2q')$ for some $l'$ and $q'$.

\medskip	  We proceed to estimate (\ref{m.9})  by breaking the calculation into two cases: when $a=e$ and all the cycles of $G_{\pi}$ are of order two; when $a=e$ and not all the cycles of $G_{\pi}$ are of order two or when $a\ne e$.

 \subsection{${\bf a =e}$, with all cycles of order two}\label{assinged}

Let $\mathcal{P}=\{(\ga_{2v-1},\ga_{2v})\,,\,1\leq v\leq n\}$ be a pairing of the $m$ vertices   
\[\{(j,k,i),\,0\leq j<k\leq K,\, 1\leq i\leq m_{j,k}\}\] of $G_{\pi }$,  that satisfies the following special property: whenever $(j,k,i)$ and $(j',k',i')$ are paired together,   $j=j'$  and $k=k'$.    Equivalently,
\begin{equation}
\mathcal{P}=\bigcup_{\stackrel{j,k=0}{j< k}}^{K} \mathcal{P}_{j,k}\label{pair}
\end{equation}
where each $\mathcal{P}_{j,k}$ is a pairing of the $m_{j,k}$ vertices
\[\{(j,k,i),\, 1\leq i\leq m_{j,k}\}.\] 
We refer to such a pairing $\mathcal{P}$ as a special pairing and  denote the set of special pairings by $\mathcal{S}$.

Given a special  pairing $\mathcal{P}\in \mathcal{S}$, let $\pi$   be   such that for each $ 0\leq l\leq K$ and $1\leq q\leq n_{l}$, 

\be
\{\pi_{l}(2q-1), \pi_{l}(2q)\}=\{\ga_{2v-1},\ga_{2v}\}\label{4.34}
\ee
 for some, necessarily unique,   $ 1\leq v\leq n_{l}$. In this case we say that $\pi$ is compatible with the pairing $\mathcal{P}$ and write this  as $ \pi \sim \mathcal{P}$. (Recall that  when we write $\{\pi_{l}(2q-1), \pi_{l}(2q)\}=\{\ga_{2v-1},\ga_{2v}\}$, we mean as two  sets, so, according to what $\pi_{l}$ is, we may have  $\pi_{l}(2q-1)=\ga_{2v-1}$ and $ \pi_{l}(2q)=\ga_{2v}$ or $\pi_{l}(2q-1)=\ga_{2v}$ and $ \pi_{l}(2q)=\ga_{2v-1}$.)  Clearly
\begin{equation}
| \mathcal{S} |\leq {(2n)! \over 2^{n}n!}\label{ord.1}
\end{equation}
the number of pairings of $m=2n$ objects.

	Let $\pi\in \mathcal{B}$   be  such that $G_{\pi}$ consists of cycles of order two. It is easy to see that $ \pi \sim \mathcal{P}$ for some $\mathcal{P}\in \mathcal{S}$. To see this note that if $\{(j,k,i),(j',k',i')\}$  form a cycle of order two, there must exist   $l$ and $l'$ with $l\neq l'$ and $q$ and $q'$ such that both
$\{(j,k,i),(j',k',i')\}=\{\pi_{l} (2q-1), \pi_{l} (2q)\}$ and $\{(j,k,i),(j',k',i')\}=\{\pi_{l'} (2q'-1), \pi_{l'} (2q')\}$.   This implies that   $j=j'$, $k=k'$ and 
$\{j,k\}=\{l,l'\}$.  Furthermore, by (\ref{4.34}) we have
\begin{equation}
   \{\pi_{l}(2q-1), \pi_{l}(2q)\}=\{\pi_{l'} (2q'-1), \pi_{l'} (2q')\}=\{\ga_{2v-1},\ga_{2v}\}\label{4.36}
   \end{equation}

  When $\pi\sim \PP$ and all cycles are of order two we can write
\bea
  && \prod_{l=0}^{K}\prod_{q=1}^{n_{l}}  \De^{ 1} 
\De^{-1} \,p_{r_{l,2q}}(x_{\pi_{l}(2q)}-x_{\pi_{l}(2q-1)})\label{4.37} \\
&&\qquad=\prod_{v=1}^{n}\De^{1}\De^{-1} \,p_{r_{2\nu}}( x_{\ga_{2v}}-x_{\ga_{2v-1}})\De^{1}\De^{-1} \,p_{r'_{2\nu}}( x_{\ga_{2v}}-x_{\ga_{2v-1}})\nn,
   \eea
where $r_{2\nu}$ and $r'_{2\nu}$ are the rearranged indices $r_{l,2q}$ and $r_{l',2q'}$. We also use the fact that $\sum_{l=0}^{K}n_{l}=2n$.

For use in (\ref{91.4})   below we note that   
\bea 
\lefteqn{
 \int_{0}^{t}\!\! \int_{0}^{t}| \De^{1}\De^{-1} \,p_{r_{2\nu}}( x_{\ga_{2v}}-x_{\ga_{2v-1}})|\,|\De^{1}\De^{-1} \,p_{r'_{2\nu}}( x_{\ga_{2v}}-x_{\ga_{2v-1}})|\,dr_{2\nu}\,dr'_{2\nu}\nn}\\ && =\(\int_{0}^{t}|\De^{1}\De^{-1} \,p_{r }( x_{\ga_{2v}}-x_{\ga_{2v-1}})|\,dr\)^{2}= w^{2}(x_{\ga_{2v}}-x_{\ga_{2v-1}},t),\,\,\,\label{4.37j}\,\,\,
 \eea
  (see (\ref{jrst1.3y}).)  

\medskip	 We want to estimate the integrals in (\ref{m.8}). However, it is difficult   to integrate $\mathcal{T}_{t}( x;\,\pi ,e )$ directly, because   the variables, 
\bea
&&\{x_{\pi_{l}(1)}-x_{\pi_{l-1}(m_{l-1})},\,x_{\pi_{l}(2q-1)}-x_{\pi_{l}(2q-2)},\,x_{\pi_{l}(2q)}-x_{\pi_{l}(2q-1)};\nn \\&& \hspace{1.5in}
\,l\in [0,K],\, q\in [1,n_{l}]\}\nn,
\eea
are not independent.   We begin the estimation by showing that over much of the   domain of integration, the integral is negligible, asymptotically, as $t\to \ff$. To begin, we write  
\begin{equation}
1=\prod_{ v=1}^{ n}\(1_{\{|x_{\ga_{2v}}-x_{\ga_{2v-1}}|\leq t^{(\bb-1)/(4\bb)}\}}+
1_{\{|x_{\ga_{2v}}-x_{\ga_{2v-1}}|\geq t^{(\bb-1)/(4\bb)}\}}\)\label{expm.1}
\end{equation}
and expand it as a sum of $2^{n}$ terms and use it to  write
\begin{eqnarray}
&&
\int \mathcal{T}_{t}( x;\,\pi,e)\prod_{j,k,i}\,dx_{j,k,i}\label{91.3a}\\
&&\qquad =\int \prod_{ v=1}^{ n}\(1_{\{|x_{\ga_{2v}}-x_{\ga_{2v-1}}|\leq t^{(\bb-1)/(4\bb)}\}}\)\mathcal{T}_{t}( x;\,\pi,e)\prod_{j,k,i}\,dx_{j,k,i}+E_{1,t}.\nonumber 
\end{eqnarray}
We now show that 
\begin{equation}
E_{1,t}=O\(t^{-(\bb-1)/(5\bb)}\(t^{2}\psi^{-1}(1/t)\)^{n}\).\label{91.3}
\end{equation}
Note that every term in $E_{1,t}$   can be written in the form 
\begin{equation}
B_{t}(\pi,e, D ):=\int \prod_{v=1}^{n} 1_{ D_{v} } 
\mathcal{T}_{t}( x;\,\pi,e)\prod_{j,k,i}\,dx_{j,k,i}\label{f9.33}
\end{equation}
where each   $D_{v}$ is either $\{|x_{\ga_{2v}}-x_{\ga_{2v-1}}|\leq t^{(\bb-1)/(4\bb)}\}$ or $\{|x_{\ga_{2v}}-x_{\ga_{2v-1}}|\geq t^{(\bb-1)/(4\bb)}\}$, and at least one of the $D_{v}$ is of the second type.  

 Consider (\ref{f9.33}) and the representation of $\mathcal{T}_{t}( x;\,\pi,e)$ in (\ref{m.10}). We  
 take absolute values in the integrand in  (\ref{m.10})   and
  take all the integrals with $r_{\cdot}$  between 0 and t  and use (\ref{4.37j})   
   to get   
\begin{eqnarray} 
 | B_{t}(\pi,e, D )|&\leq &
 \int \prod_{v=1}^{n} 1_{ D_{v} }  w^{2}( x_{\ga_{2v}}-x_{\ga_{2v-1}},t)  \prod_{l=0}^{K}u(x_{\pi_{l}(1)}-x_{\pi_{l-1}(m_{l-1})},t)   \nonumber\\
&& \qquad
  \prod_{q=2}^{n_{l}} u(x_{\pi_{l}(2q-1)}-x_{\pi_{l}(2q-2)},t) \prod_{j,k,i}\,dx_{j,k,i} \label{91.4}.
\end{eqnarray}
 We now take 
\be
\{x_{\ga_{2v}}-x_{\ga_{2v-1}},\,v=1,\ldots,n \}
\ee
 and an additional $n$ 
variables from  the $2n $ arguments of the $u$ terms,
\be \cup_{l=0}^{K} \{x_{\pi_{l}(1)}-x_{\pi_{l-1}(m_{l-1})}, x_{\pi_{l}(2q-1)}-x_{\pi_{l}(2q-2)},\,q=2,\ldots,n_{l} \}\label{kik}\ee
so that the chosen $2n$ variables generate   the space spanned by the 2n variables $\{x_{j,k,i} \}$. 
There are $n$ variables in (\ref{kik}) that are not used.  We bound  the functions $u$  of these variables by their sup norm, which by (\ref{jr.1}) is bounded by $Ct\psi^{-1}(1/t)$. Then we make a change of variables and get that  
\bea 
  |B_{t}(\pi,e, D )|&\le& C\(t\psi^{-1}(1/t)\)^{n}\int \prod_{v=1}^{n} 1_{ D_{v} }  w^{2}( y_{v},t)\prod_{v=n+1}^{2n}u(y_{v},t))\prod_{v=1}^{2n}\,dy_{v}
\nn\\
&\le&C\(t^{2}\psi^{-1}(1/t)\)^{n}\int \prod_{v=1}^{n} 1_{ D_{v} }  w^{2}( y_{v},t)) \prod_{v=1}^{n}\,dy_{v}\nn,\\
&=&O\(t^{-(\bb-1)/(5\bb)}\(t^{2}\psi^{-1}(1/t)\)^{n}\).\label{91.oo}
\eea
  Here we use  (\ref{ee1})
to see that the integral of a $u$ term is $t$.
 Then we use (\ref{132ee2}) and  (\ref{132ee3}) to obtain (\ref{91.3}).   (Note that it is because at least one of the $D_{v}$ is of the second type that we can use (\ref{132ee3}).)

\medskip	We now study  
\begin{equation} \hspace{.4 in}
\int \prod_{ v=1}^{ n}\(1_{\{|x_{\ga_{2v}}-x_{\ga_{2v-1}}|\leq t^{(\bb-1)/(4\bb)}\}}\)\mathcal{T}_{t}( x;\,\pi,e)\prod_{j,k,i}\,dx_{j,k,i}. \label{f91.36}
\end{equation}
  We identify the  relationships in (\ref{4.36})  by setting $v=\si_{l} (q) $  so that
\begin{equation}
\{\pi_{l}(2q-1), \pi_{l}(2q)\}=\{\ga_{2\si_{l} (q)-1},\ga_{2\si_{l} (q)}\},\label{3.51}
   \end{equation}
   for each $ 0\leq l\leq K$ and $1\leq q\leq n_{l}$. We use both (\ref{4.36}) and (\ref{3.51}) in what follows.

  We now make a change of variables that, eventually, enables us to make the arguments of the $u$ terms and the $w$ terms independent.
 For $q\geq 2$ we write 
\begin{eqnarray}
\lefteqn{p_{r_{l,2q-1}}(x_{\pi_{l}(2q-1)}-x_{\pi_{l}(2q-2)})
  \label{f91.37}}\\
&&=p_{r_{l,2q-1}}(x_{\ga_{2\si_{l} (q)-1}}-x_{\ga_{2\si_{l} (q-1)-1}})+\De^{h_{l,q}} p_{r_{l,2q-1}}(x_{\ga_{2\si_{l} (q)-1}}-x_{\ga_{2\si_{l} (q-1)-1}}), \nn
\end{eqnarray}
where $h_{l,q}=(x_{\pi_{l}(2q-1)}-x_{\ga_{2\si_{l} (q)-1}})+(x_{\ga_{2\si_{l} (q-1)-1}}-x_{\pi_{l}(2q-2)})$.  When $q=1$  we can make a similar decomposition
\begin{eqnarray}
&&p_{(t-\sum_{q=1}^{m_{l-1}}r_{l-1,q})+r_{l,1}}(x_{\pi_{l}(1)}-x_{\pi_{l-1}(m_{l-1})})
\label{f91.37.1}\\
&&\qquad= p_{(t-\sum_{q=1}^{m_{l-1}}r_{l-1,q})+r_{l,1}}(x_{\ga_{2\si_{l} (1)-1}}-x_{\ga_{2\si_{l-1} (n_{l-1})-1}})  \nonumber\\
&&\qquad\quad+\De^{h_{l,1}} p_{(1-\sum_{q=1}^{m_{l-1}}r_{l-1,q})+r_{l,1}}(x_{\ga_{2\si_{l} (1)-1}}-x_{\ga_{2\si_{l-1} (n_{l-1})-1}}) , \nonumber
\end{eqnarray}
where $h_{l,1}=(x_{\pi_{l}(1)}-x_{\ga_{2\si_{l} (1)-1}})+ (x_{\ga_{2\si_{l-1} (n_{l-1})-1}}-x_{\pi_{l-1}(m_{l-1})})$.
Note that because of the presence of  the term   $\prod_{ v=1}^{ n}\(1_{\{|x_{\ga_{2v}}-x_{\ga_{2v-1}}|\leq t^{(\bb-1)/(4\bb)}\}}\)$ in the integral in (\ref{f91.36})  we need only be concerned with values of  $|h_{l,q}|\leq 2t^{(\bb-1)/(4\bb)}$,   $ 0\leq l\leq K$ and   $1\leq q\leq n_{l}$.

  For $q=1,\ldots,n_{l}$, $l=0\ldots,K$, we substitute   (\ref{f91.37}) and (\ref{f91.37.1})   into the term $\mathcal{T}_{t}( x;\,\pi ,e )$ in  (\ref{f91.36}), (see also (\ref{m.10})),   and expand the products so that we can write 
(\ref{f91.36}) as a sum of $2^{\sum_{l=0}^{K}n_{l} }$  terms, which we write as  
\begin{eqnarray}
&&\int \prod_{ v=1}^{ n}\(1_{\{|x_{\ga_{2v}}-x_{\ga_{2v-1}}|\leq t^{(\bb-1)/(4\bb)}\}}\)\mathcal{T}_{t}( x;\,\pi,e)\prod_{j,k,i}\,dx_{j,k,i}
\label{91.6}\\
&&\qquad=\int \prod_{ v=1}^{ n}\(1_{\{|x_{\ga_{2v}}-x_{\ga_{2v-1}}|\leq t^{(\bb-1)/(4\bb)}\}}\)\mathcal{T}_{t,1}( x;\,\pi,e)\prod_{j,k,i}\,dx_{j,k,i} + E_{2,t},\nonumber
\end{eqnarray}
  where
\begin{eqnarray}
&&\mathcal{T}_{t,1}( x;\,\pi ,e )=   \prod_{l=0}^{K}\int_{\RR_{l}}p_{(t-\sum_{q=1}^{m_{l-1}}r_{l-1,q})+r_{l,1}}(x_{\ga_{2\si_{l} (1)-1}}-x_{\ga_{2\si_{l-1} (n_{l-1})-1}}) 
\nn\\
&&
\hspace{1.4 in}\prod_{q=2}^{n_{l}} p_{r_{l,2q-1}}(x_{\ga_{2\si_{l} (q)-1}}-x_{\ga_{2\si_{l} (q-1)-1}})  \label{m.11}\\
&&\hspace{1.6 in} \prod_{q=1}^{n_{l}} \De^{ 1} 
\De^{-1} \,p_{r_{l,2q}}(x_{\pi_{l}(2q)}-x_{\pi_{l}(2q-1)})\prod_{q=1}^{m_{l}}\,dr_{l,q}.\nn
\end{eqnarray} 
 Using  (\ref{4.37}) we can rewrite this as  
 \begin{eqnarray}
&&\mathcal{T}_{t,1}( x;\,\pi ,e )\label{m.11c}\\
&&\qquad
=   \int_{\RR_{0}\times\cdots\times \RR_{K}}\(\prod_{l=0}^{K} p_{(t-\sum_{q=1}^{m_{l-1}}r_{l-1,q})+r_{l,1}}(x_{\ga_{2\si_{l} (1)-1}}-x_{\ga_{2\si_{l-1} (n_{l-1})-1}})\right. 
\nn\\
&&\left. 
\hspace{2 in}\prod_{q=2}^{n_{l}} p_{r_{l,2q-1}}(x_{\ga_{2\si_{l} (q)-1}}-x_{\ga_{2\si_{l} (q-1)-1}})\)  \nn\\ 
&&\hspace{.1 in}\( \prod_{v=1}^{n}\De^{1}\De^{-1} \,p_{r_{2\nu}}( x_{\ga_{2v}}-x_{\ga_{2v-1}})\De^{1}\De^{-1} \,p_{r'_{2\nu}}( x_{\ga_{2v}}-x_{\ga_{2v-1}}) \)\nn\\
&&\hspace{3.5in}\prod_{l=0}^{K}\prod_{q=1}^{m_{l}}\,dr_{l,q}\nn,
\end{eqnarray}
where $r_{2\nu}$ and $r'_{2\nu}$ are the rearranged indices $r_{l,2q}$ and $r_{l',2q'}$. 

  The usefulness of the representations in (\ref{f91.37}) and (\ref{f91.37.1}) is now apparent. Since the variables $x_{\ga_{2v}},\,v=1,\ldots, n$, occur only in the last line of (\ref{m.11c}),    we make the change of variables $ x_{\ga_{2v}}-x_{\ga_{2v-1}}\to  x_{\ga_{2v}} $ and $x_{\ga_{2v-1}}\to x_{\ga_{2v-1}}$ and get that 
\begin{eqnarray}
\lefteqn{\int  \mathcal{T}_{t,1}( x;\,\pi,e)\prod_{j,k,i}\,dx_{j,k,i}
\label{m.11d}}\\
&&
= \int  \int_{\RR_{0}\times\cdots\times \RR_{K}}
\(\prod_{l=0}^{K} p_{(t-\sum_{q=1}^{m_{l-1}}r_{l-1,q})+r_{l,1}}(x_{\ga_{2\si_{l} (1)-1}}-x_{\ga_{2\si_{l-1} (n_{l-1})-1}}) \right.\nn\\
&&\hspace{1.9in}\left. \prod_{q=2}^{n_{l}} p_{r_{l,2q-1}}(x_{\ga_{2\si_{l} (q)-1}}-x_{\ga_{2\si_{l} (q-1)-1}})\)  \nn\\
&&\qquad \( \prod_{v=1}^{n}\De^{1}\De^{-1} \,p_{r_{2\nu}}( x_{\ga_{2v}})\De^{1}\De^{-1} \,p_{r'_{2\nu}}( x_{\ga_{2v}}) \)\prod_{l=0}^{K}\prod_{q=1}^{m_{l}}\,dr_{l,q}
\prod_{j,k,i}\,dx_{j,k,i}.\nonumber
\end{eqnarray}
  Now, since the variables $x_{\ga_{2v}},\,v=1,\ldots, n$ occur only in the last line of (\ref{m.11d}) and the variables $x_{\ga_{2v-1}},\,v=1,\ldots, n$ occur only in the second and third lines of (\ref{m.11d}), we can  write (\ref{m.11d}) as
\begin{eqnarray}
\lefteqn{\int \mathcal{T}_{t,1}( x;\,\pi,e)\prod_{j,k,i}\,dx_{j,k,i}
\label{m.11e}}\\
&&
=   \int_{\RR_{0}\times\cdots\times \RR_{K}} \int\(\prod_{l=0}^{K} p_{(t-\sum_{q=1}^{m_{l-1}}r_{l-1,q})+r_{l,1}}(x_{\ga_{2\si_{l} (1)-1}}-x_{\ga_{2\si_{l-1} (n_{l-1})-1}})\right. 
\nn\\
&&\left. 
\hspace{1.5 in}\prod_{q=2}^{n_{l}} p_{r_{l,2q-1}}(x_{\ga_{2\si_{l} (q)-1}}-x_{\ga_{2\si_{l} (q-1)-1}})\) 
 \prod_{v=1}^{n} \,dx_{\ga_{2v-1}}\nn\\
&&\hspace{.5 in}\( \prod_{v=1}^{n}\int \De^{1}\De^{-1} \,p_{r_{2\nu}}( x_{\ga_{2v}})\De^{1}\De^{-1} \,p_{r'_{2\nu}}( x_{\ga_{2v}}) \,dx_{\ga_{2v}}\)\prod_{l=0}^{K}\prod_{q=1}^{m_{l}}\,dr_{l,q}.\nn
\end{eqnarray}
 Note that we also use Fubini's Theorem, which is    justified since the    absolute value of the integrand is integrable, (as we point out in the argument preceding (\ref{91.4})).   (In the rest of this section use Fubini's Theorem frequently for  integrals like (\ref{m.11e})  without repeating the explanation about why it is justified.) 
  
    We now show that   
 \be
 E_{2,t}=O\(t^{-(\bb-1)/(3\bb)}\(t^{2}\psi^{-1}(1/t)\)^{n}\).\label{3.24}
 \ee
 To see this note that   the terms in $E_{2,t}$ are of the form    
\begin{eqnarray}
&& 
\int  \prod_{ v=1}^{ n}\(1_{\{|x_{\ga_{2v}}-x_{\ga_{2v-1}}|\leq t^{(\bb-1)/(4\bb)}\}}\) \label{f9.40}\\
&&\hspace{.5in} \prod_{l=0}^{K}\int_{\RR_{l}}\wt p_{(t-\sum_{q=1}^{m_{l-1}}r_{l-1,q})+r_{l,1}}(x_{\ga_{2\si_{l} (1)-1}}-x_{\ga_{2\si_{l-1} (n_{l-1})-1}}) 
\nn\\
&&
\hspace{.7 in}\prod_{q=2}^{n_{l}}\wt  p_{r_{l,2q-1}}(x_{\ga_{2\si_{l} (q)-1}}-x_{\ga_{2\si_{l} (q-1)-1}})    \nonumber\\
&&\hspace{.9 in}\prod_{q=1}^{n_{l}} \De^{ 1} 
\De^{-1} \,p_{r_{l,2q}}(x_{\pi_{l}(2q)}-x_{\pi_{l}(2q-1)})\prod_{q=1}^{m_{l}}\,dr_{l,q}\prod_{j,k,i}\,dx_{j,k,i}\nn,
\end{eqnarray}
where $\wt p_{r_{l,2q-1}}$ is either $  p_{r_{l,2q-1}}$ or $\De^{h_{l,q}} p_{r_{l,2q-1}}$. Furthermore,   at least one of the terms  $\wt p_{r_{l,2q-1}}$  is of the form      $\De^{h_{l,q}} p_{r_{l,2q-1}}$.  

As in the transition from (\ref{f9.33}) to (\ref{91.4}) we  bound the absolute value of (\ref{f9.40})    by  
\begin{eqnarray}
\lefteqn{\int \prod_{ v=1}^{ n}\(1_{\{|x_{\ga_{2v}}-x_{\ga_{2v-1}}|\leq t^{(\bb-1)/(4\bb)}\}}\)w^{2}( x_{\ga_{2v}}-x_{\ga_{2v-1}},t) 
\label{m.12}}\\
&&
\prod_{l=0}^{K}\wt u(x_{\ga_{2\si_{l} (1)-1}}-x_{\ga_{2\si_{l-1} (n_{l-1})-1}},t) \prod_{q=2}^{n_{l}} \wt u(x_{\ga_{2\si_{l} (q)-1}}-x_{\ga_{2\si_{l} (q-1)-1}},t) \nn\\
&&\hspace{4in}  \prod_{j,k,i}\,dx_{j,k,i},\nn
\end{eqnarray} where each $\wt u(\cd,t)$ is either of the form $u(\cd,t)$ or $v_{h_{l,q}}(\cd,t)$; (see (\ref{2.3})).

\medskip  We need to introduce the following notation and estimates. The next lemma is proved in Section \ref{sec-lemproofs}.  Let 
\begin{equation}
   v _{*}(x,t):=\(\log t\wedge{ t\psi^{-1}(1/t)\over |x|}\wedge t  {1+  \log^{+} x\over x^{2}}\)\label{vstar},
   \end{equation}
\begin{lemma}\label{lem-3.2} Under the hypotheses of 
Theorem  \ref{theo-clt2r},  for all $t$ sufficiently large,
\bea 
 v_{h_{l,q}}(x,t) &\le& C h_{l,q}^{2}  v _{*}(x,t)\label{v_{*}}\\
\sup_{x\in R^{1}}v_{*}(x,t)&\le&   \log t \label{v_{**}}\\                 
   \int\,v_{*}(x,t)\,dx&\le &  C   t\(\psi^{-1}(1/  t)\)\log t.     \label{v_{***}}
 \eea

 \el

We  have $J$ terms of the type $v_{h_{l,q}}(\cd,t)$, for some   $J\geq 1$. It follows from       (\ref{v_{*}})  and   and the fact that $|h_{l,q}|\leq 2t^{(\bb-1)/(4\bb)}$, that  we can  
 bound the  integral in (\ref{m.12}) by  
\begin{eqnarray}
\lefteqn{
 Ct^{J(\bb-1)/(2\bb)}\int \prod_{ v=1}^{ n} w^{2}( x_{\ga_{2v}}-x_{\ga_{2v-1}},t) 
	\label{m.13}}\\
&&   \hspace{.1 in}\prod_{l=0}^{K}\wt u(x_{\ga_{2\si_{l} (1)-1}}-x_{\ga_{2\si_{l-1} (n_{l-1})-1}},t) \prod_{q=2}^{n_{l}} \wt u(x_{\ga_{2\si_{l} (q)-1}}-x_{\ga_{2\si_{l} (q-1)-1}},t) \nn\\
&&\hspace{4in} \prod_{j,k,i}\,dx_{j,k,i}, \nonumber
\end{eqnarray}
where  $\wt u(\cd,t) $ is either $ u(\cd,t) $ or $v_{*} (\cd,t)$,   and we have precisely $J$ of the latter.

 Since   the variables $x_{\ga_{2\nu}}$, $\nu=1,\ldots, n$, occur only in the $w$ terms in (\ref{m.13}) and the variables $x_{\ga_{2v-1}},\,v=1,\ldots, n$ occur only in the $\wt  u$ terms in (\ref{m.13})  , (refer to the change of variables arguments in (\ref{m.11d}) and (\ref{m.11e})), we can  write (\ref{m.13}) as
 \begin{eqnarray}
\lefteqn{ Ct^{t^{J(\bb-1)/(2\bb)}}\int \(\prod_{l=0}^{K}\wt u(x_{\ga_{2\si_{l} (1)-1}}-x_{\ga_{2\si_{l-1} (n_{l-1})-1}},t) \right.
	\label{m.13ww}}\\
&&  \left. \hspace{.1 in}\prod_{q=2}^{n_{l}} \wt u(x_{\ga_{2\si_{l} (q)-1}}-x_{\ga_{2\si_{l} (q-1)-1}},t) \)\prod_{v=1}^{n}\,dx_{\ga_{2v-1}}\prod_{v=1}^{n}w^{2}(x_{\ga_{2v}},t )\prod_{v=1}^{n}\,dx_{\ga_{2v}}
\nonumber\\
&&\leq Ct^{t^{J(\bb-1)/(2\bb)}}  \(\log t\)^{n} \int \(\prod_{l=0}^{K}\wt u(x_{\ga_{2\si_{l} (1)-1}}-x_{\ga_{2\si_{l-1} (n_{l-1})-1}},t)  
	\right.\nn\\
&&  \left. \hspace{1.5 in}\prod_{q=2}^{n_{l}} \wt u(x_{\ga_{2\si_{l} (q)-1}}-x_{\ga_{2\si_{l} (q-1)-1}},t) \)\prod_{v=1}^{n}\,dx_{\ga_{2v-1}}
 \nn 
\end{eqnarray}
where the last inequality uses (\ref{132ee2}).

As we have been doing we extract a linearly independent set of variables from the arguments of the $\wt u$ terms. The other $\wt u$ terms we bound by their supremum. Then we make a change of variables and integrate the remaining $\wt u$ terms.  

Compare  (\ref{jr.1j}) with (\ref{v_{**}}). Replacing the sup of a $u$ term by the sup of a $v_{*}$ term reduces the upper bound by a factor of $1/(t\psi^{-1}(1/t))$, (neglecting the factor of $\log t$ which is irrelevant.)   
On the other hand, considering 
  (\ref{ee1}) and (\ref{v_{***}}), we see that replacing the integral  of a $u$ term by the integral of a $v_{*}$ term reduces the upper bound by a   factor of $\psi^{-1}(1/t)$, 
  (again neglecting the factor of $\log t$.)   
   Counting the initial   factor of 
$ t^{J(\bb-1)/(2\bb)}$ we have a  reduction with is at least
\begin{equation}
   \( t^{(\bb-1)/(2\bb)}(t\psi^{-1}(1/t))^{-1}\)^{J}=o\(\(t^ {- (\bb-1)/(3\bb)}\)^{J} \)
   \end{equation}
for all $\ep>0$.
Since $J\geq 1$, we get (\ref{3.24}).

\medskip 	  Analogous to (\ref{91.3a}) we note that  
\begin{eqnarray}
&&
\int \mathcal{T}_{t,1}( x;\,\pi,e)\prod_{j,k,i}\,dx_{j,k,i}\label{91.3ab}\\
&&\qquad   =\int \prod_{ v=1}^{ n}\(1_{\{|x_{\ga_{2v}}-x_{\ga_{2v-1}}|\leq  t^{(\bb-1)/(4\bb)}\}}\)\mathcal{T}_{t,1}( x;\,\pi,e)\prod_{j,k,i}\,dx_{j,k,i}+\wt E_{1,t},\nonumber 
\end{eqnarray}
where $\wt E_{1,t}=O\(t^{-(\bb-1)/(5\bb)}\(t^{2}\psi^{-1}(1/t)\)^{n}\)$. The proof of (\ref{91.3ab}) is the same as the proof of  (\ref{91.3}).

\medskip	
Since $\psi$ is regularly varying with index $\bb>1$ we see that there exists an   $\ep(\bb):=\ep>0$ such that
\begin{equation}
   E_{1,t}+   E_{2,t}+ \wt E_{2,t}=O\(t^{(2-1/\bb)n-\ep}\).
   \end{equation}
Therefore, it follows from (\ref{91.3a}),  (\ref{91.6})   and (\ref{91.3ab})
  that 
\begin{eqnarray}
&&\int  \mathcal{T}_{t}( x;\,\pi,e)\prod_{j,k,i}\,dx_{j,k,i}
\label{m.13q}\\
&&\qquad=\int \mathcal{T}_{t,1}( x;\,\pi,e)\prod_{j,k,i}\,dx_{j,k,i}+O\(t^{(2-1/\bb)n-\ep}\).\nonumber
\end{eqnarray}

\medskip	  We now obtain a sharp estimate , (asymptotically as $t\to\ff$), of the second integral in (\ref{m.13q}) that leads to the (\ref{m.2}).
Let $\wt\RR_{l}(s )=\{\sum_{q=1}^{n_{l}}r_{l,2q-1}\newline \leq t-s \}$ and $\wt  \si_{l} (q):=\ga_{2\si_{l} (q)-1}$.   We define
\begin{eqnarray} \lefteqn{
F_{t} (\wt  \si,s_{0},\ldots, s_{K}) 
\label{m.14}}\\
&& =\int\(\int_{\wt\RR_{0}(s_{0})\times\cdots \times\wt\RR_{K}(s_{K})} \prod_{ l=0}^{ K}\, p_{(t-\sum_{q=1}^{n_{l-1}}r_{l-1,2q-1}-s_{l-1})+r_{l,1}}\right.
\nn\\
&&\left.
\hspace{.2 in}(x_{\wt  \si_{l} (1) }-x_{\wt  \si_{l-1} (n_{l-1}) }) \prod_{q=2}^{n_{l}} p_{r_{l,2q-1}}(x_{\wt  \si_{l} (q) }-x_{\wt  \si_{l} (q-1) })
\prod_{q=1}^{n_{l}}\,dr_{l,2q-1}\)  \,d x    \nonumber,
\end{eqnarray}
where $(t-\sum_{q=1}^{n_{-1}}r_{-1,2q-1}-s_{-1}):=0$ and $ \wt  \si_{ -1} (n_{ -1}):=0$. 
Here the generic term $dx$ indicates integration with respect to all the variables $x_{\cdot}$ that appear in the integrand.

Since   $\wt  \si_{l} (q)=\ga_{2\si_{l} (q)-1}$ we can also write (\ref{m.14}) as
\begin{eqnarray} \lefteqn{
F_{t} (\wt\si,s_{0},\ldots, s_{K})\label{m.14a}}\\
&& =\int\(\int_{\wt\RR_{0}(s_{0})\times\cdots \times\wt\RR_{K}(s_{K})} \prod_{ l=0}^{ K}\, \right.
p_{(t-\sum_{q=1}^{n_{l-1}}r_{l-1,2q-1}-s_{l-1})+r_{l,1}}
\nn\\
&& \hspace{.3in}
(x_{\ga_{2\si_{l} (1)-1}}-x_{\ga_{2\si_{l-1} (n_{l-1})-1}}) \prod_{q=2}^{n_{l}} p_{r_{l,2q-1}}(x_{\ga_{2\si_{l} (q)-1}}-x_{\ga_{2\si_{l} (q-1)-1}})
  \nonumber\\
&&\hspace{3
in}\left.\prod_{q=1}^{n_{l}}\,dr_{l,2q-1}\)  \,dx, \nn
\end{eqnarray} 
with  $x_{\ga_{2\si_{-1} (n_{-1})-1}}:=0$.

 \medskip	Consider (\ref{m.14a}). By extending the time integration we have 
 \begin{eqnarray} \lefteqn{
F_{t} (\wt\si,s_{0},\ldots, s_{K})\label{m.14av}}\\
&& \leq \int   \prod_{ l=0}^{ K}\,  
u(x_{\ga_{2\si_{l} (1)-1}}-x_{\ga_{2\si_{l-1} (n_{l-1})-1}})
\nn\\
&& \hspace{1in} 
( \prod_{q=2}^{n_{l}} u_{r_{l,2q-1}}(x_{\ga_{2\si_{l} (q)-1}}-x_{\ga_{2\si_{l} (q-1)-1}})
  \,dx. \nn
\end{eqnarray}
 Note that there are $n$ different $x_{\cd}$ variables, each one of which appears twice. Therefore, by an argument similar to the one   in the paragraph containing (\ref{91.oo}), we see that 
   \begin{equation}
   F_{t} (\wt\si,s_{0},\ldots, s_{K}) \le C\(t^{2}\psi^{-1}(1/t)\)^{n},\label{4.60}
   \end{equation}
 for some constant depending only on   $m=2n$.

\medskip	
  Let $\wh\RR_{l} =\{\sum_{q=1}^{n_{l}}r_{l,2q}\leq t\}$,   $l=0,\ldots,K$.  We break up the integration over $ \RR_{0}\times\cdots \times \RR_{K}$ in  (\ref{m.11e}) into integration over $\wt\RR_{0}\times\cdots \times\wt\RR_{K} $ and $\wh\RR_{0}\times\cdots \times\wh\RR_{K}$; (see 
  (\ref{m.14a})).  If one carefully examines  the time indices in (\ref{m.9}) and (\ref{m.14}) and uses Fubini's Theorem, one sees that  
  \begin{eqnarray} 
&& \int \mathcal{T}_{t}( x;\,\pi,e)\prod_{j,k,i}\,dx_{j,k,i}\label{m.15}\\ 
&&\qquad=
 \int_{\wh\RR_{0}\times\cdots \times\wh\RR_{K}}    
F_{t} (\wt\si,\sum_{q=1}^{n_{0}}r_{0,2q}\,,\ldots, \sum_{q=1}^{n_{K}}r_{K,2q})  \nn\\ 
  &&
 \hspace{.6 in} \prod_{i=1}^{n} \(\int \(  \De^{ 1}\De^{-1}
\,p_{r_{i} }(x ) \)  \(\De^{ 1}\De^{-1}
\,p_{r'_{i} }(x )\) \,dx\) \prod_{i=1}^{n}\,dr_{i}\,dr'_{i}.\nn
\end{eqnarray}
The variables $\{r_{i},r'_{i}\,|\,i=1,\ldots,n\}$ are simply a relabeling of the variables $\{r_{l,2q}\,|\,0\leq l\leq K,\,1\leq q\leq n_{l}\}$.
(The exact form of this relabeling  does not matter in what follows.) Here, as always, we set $p_{r}(x)=0$, if $r\leq 0$.

\medskip	 By   Parseval's Theorem  
\begin{eqnarray}
&&\int \(\De^{ 1}\De^{-1}\,p_{r  }(x)\)
\(\De^{ 1}\De^{-1}\,p_{r' }(x)\) \,dx
\label{f9.2}\\
&&\qquad ={1 \over 2\pi}\int |2-e^{ip }-e^{-ip }  |^{2}e^{-r\psi (p) } e^{-r'\psi (p) }\,dp  \nonumber\\
&&  \qquad ={16 \over  \pi}\int_{0}^{\ff} \sin^{4}(p/2)  e^{-r\psi (p) } e^{-r'\psi (p) }\,dp\geq 0\nn.
\end{eqnarray}
 Using this,   (\ref{4.60}) and Fubini's Theorem, we see  that 
\begin{eqnarray}
\lefteqn{ \int_{\(\wh\RR_{0}\times\cdots \times\wh\RR_{K}\)\cap ([0,\sqrt{t}]^{2n})^{c}}    
F_{t} (\wt\si,\sum_{q=1}^{n_{0}}r_{0,2q}\,,\ldots, \sum_{q=1}^{n_{K}}r_{K,2q})
\label{parm.1}}\\
&&\qquad   \prod_{i=1}^{n} \(\int \(  \De^{ 1}\De^{-1}
\,p_{r_{i} }(x ) \)  \(\De^{ 1}\De^{-1}
\,p_{r'_{i} }(x )\) \,dx\)\prod_{i=1}^{n}\,dr_{i}\,dr'_{i}\nn\\
&& \leq C\(t^{2}\psi^{-1}(1/t)\)^{n}\nonumber\\
&&\hspace{.4 in} \int_{ ([0,\sqrt{t}]^{2n})^{c}}    
   \prod_{i=1}^{n} \(\int \(  \De^{ 1}\De^{-1}
\,p_{r_{i} }(x ) \)  \(\De^{ 1}\De^{-1}
\,p_{r'_{i} }(x )\) \,dx\)\prod_{i=1}^{n}\,dr_{i}\,dr'_{i}\nn\\
&& \leq C \(t^{2}\psi^{-1}(1/t)\)^{n} \(\int \(\int ( \De^{ 1}\De^{-1}
\,p_{r  }(x ) )  \,dr \)^{2}  \,dx \)^{n-1}\nn\\
&&  \qquad\quad\int \bigg\{\int_{0}^{\ff}\int_{\sqrt t}^{\ff} \(  \De^{ 1}\De^{-1}
\,p_{r_{i} }(x ) \)  \(\De^{ 1}\De^{-1}
\,p_{r'_{i} }(x )\)\,dr_{i}\,dr'_{i} \bigg\}\,dx
\nn\\
&&  \le C  \(t^{2}\psi^{-1}(1/t) \)^{n} \nonumber\\
&& \hspace{.4 in}\int \bigg\{\int_{0}^{\ff}\int_{\sqrt t}^{\ff} \(  \De^{ 1}\De^{-1}
\,p_{r_{i} }(x ) \)  \(\De^{ 1}\De^{-1}
\,p_{r'_{i} }(x )\)\,dr_{i}\,dr'_{i} \bigg\}\,dx,
\nn
\eea by (\ref{132ee2}).
  By (\ref{h3.1}) and (\ref{h3.2}) the integral in the final line of  (\ref{parm.1})
\bea 
   \le  c_{\psi,1} - \int  \(\int_{0}^{\sqrt{t}}\De^{ 1}\De^{-1}\,\,p_{s }(x)\,ds\)^{2}\,dx\le O\(t^{-1/6}\).
 \eea
Therefore  the first integral in (\ref{parm.1})
is  
   $O(t^{(2-1/\bb)n-\ep}),
$
for some $\ep>0$.

 \medskip	
Since $\(\wh\RR_{0}\times\cdots \times\wh\RR_{K}\)\supseteq [0,\sqrt{t}]^{2n}$, for $2n\sqrt{t}\leq t$, it follows from (\ref{m.15})  and the preceding sentence, that
 \begin{eqnarray} 
&&\lefteqn{ \int \mathcal{T}_{t}( x;\,\pi,e)\prod_{j,k,i}\,dx_{j,k,i}\label{m.15s}}\\ 
&&\qquad=
 \int_{[0,\sqrt{t}]^{2n}}    
F_{t} (\wt\si,\sum_{q=1}^{n_{0}}r_{0,2q}\,,\ldots, \sum_{q=1}^{n_{K}}r_{K,2q})  \prod_{i=1}^{n} \bigg(\int \(  \De^{ 1}\De^{-1}
\,p_{r_{i} }(x ) \)   \nn\\ 
&&\hspace{.7in}
  \(\De^{ 1}\De^{-1}
\,p_{r'_{i} }(x )\) \,dx\bigg)\prod_{ l=0}^{ K}\prod_{q=1}^{n_{l}}\,dr_{l,2q}+O(t^{(2-1/\bb)n-\ep}) \nn
\end{eqnarray}
 
We use the next lemma which is proved in  Subsection \ref{sec-m4}. 

\begin{lemma}\label{lem-m4}  Under the hypotheses of Theorem \ref{theo-clt2r}, for any fixed $m$ and  $s_{0},\ldots, s_{K}\leq m\sqrt{t}$  and $1<\bb\le 2$,  there exists an $\ep>0$ such that for all $t> 0$, sufficiently large, 
\begin{equation}
|F_{t} (\wt\si,s_{0},\ldots, s_{K})-F_{t} (\wt\si,0,\ldots, 0)|\leq C\(t^{2}\psi^{-1}(1/t)\)^{ n-\ep}.\label{m4}
\end{equation}
  \end{lemma} 
  
 \noindent{\bf Proof of Lemma \ref{lem-multiple} continued: }  It follows from  (\ref{m.15s}) and Lemmas  \ref{lem-m4} and    \ref{lem-2.5}, that
  \begin{eqnarray} 
&&\lefteqn{ \int \mathcal{T}_{t}( x;\,\pi,e)\prod_{j,k,i}\,dx_{j,k,i}
\label{m.15sw}}\\ 
&&\qquad=F_{t} (\wt\si, 0,\ldots,0) 
 \int_{[0,\sqrt{t}]^{2n}}    
 \prod_{i=1}^{n} \bigg(\int \(  \De^{ 1}\De^{-1}
\,p_{r_{i} }(x ) \)   \nn\\ 
&&\hspace{.7in}
  \(\De^{ 1}\De^{-1}
\,p_{r'_{i} }(x )\) \,dx\bigg)\prod_{ l=0}^{ K}\prod_{q=1}^{n_{l}}\,dr_{l,2q}+O\(t^{(2-1/\bb)n-\ep}\)\nn\\
&&\qquad=\(c_{\psi,1}\)^{n}
  F_{t} (\wt\si,0,\ldots, 0)+O\(t^{(2-1/\bb)n-\ep}\),\nn
\end{eqnarray}
  for some $\ep>0$.

\medskip	  Consider the mappings  $\wt  \si_{l}$  that are used in (\ref{m.14}).   Recall that $\si_{l} (q) $ is defined by the relationship  $\{\pi_{l}(2q-1), \pi_{l}(2q)\}=\{\ga_{2\si_{l} (q)-1},\ga_{2\si_{l} (q)}\}$.   Therefore,    since $\wt  \si_{l}(q)=\ga_{2\si_{l} (q)-1}$
  we can have that either  $\wt\si_{l} (q) =\pi_{l}(2q-1)$ or
$\wt\si_{l} (q) =\pi_{l}(2q)$.  
 However, since the terms $\wt  \si_{l}(q)$ are subscripts of the terms $x$, all of which are integrated, it is more convenient to define   $\wt  \si_{l}$ differently. 

 Recall that $\mathcal{P}$,  (see (\ref{pair})), is a union of pairings $\mathcal{P}_{j,k}$
 of the $m_{j,k}$ vertices
\[\{(j,k,i),\, 1\leq i\leq m_{j,k}\}.\] 
Each $\mathcal{P}_{j,k}$ consists of $n_{j,k}$ pairs, that  can ordered arbitrarily.  Consider one such ordering. If 
$\{\ga_{2\si_{l} (q)-1},\newline \ga_{2\si_{l} (q)}\}$ is the $i$-th pair in $\mathcal{P}_{j,k}$, we set 
$\wt  \si_{l}(q)=(j,k,i)$. (Necessarily, $l$ will be either $j$ or $k$, as we point out in the paragraph containing (\ref{4.36})). Thus, each $\wt  \si_{l}$ is a bijection from $[1,n_{l}]$ to
\begin{equation}
\wt I_{l}=\bigcup_{k=l+1}^{K}\{(l,k,i),\,1\leq i\leq n_{l,k}\}\bigcup_{j=0}^{l-1}\{(j,l,i),\,1\leq i\leq n_{j,l}\}.\label{m.5kp}
\end{equation}
Let $\wt \mathcal{B}$  denote the set of $K+1$ tuples, $\wt\si=(\wt\si_{0},\ldots, \wt\si_{K})$ of such bijections.  Note that with this definition of $\wt  \si_{l}(q)$,  (\ref{m.14})  remains  unchanged since we have simply renamed the variables of integration.

 \label{page33}By interchanging the elements in any of the $2n$ pairs $\{\pi_{l}(2q-1),\pi_{l}(2q)\}$ 
 we obtain a new $\pi' \sim \mathcal{P}$. In fact we obtain $2^{2n}$ different permutations $\pi$, in this way,  all of which are compatible with $\mathcal{P}$, and all of which   give   the same $\wt\si$ in (\ref{m.14}). Furthermore, by permuting the 
 pairs $\{\pi_{l}(2q-1),\pi_{l}(2q)\}$, $1\leq q\leq n_{l}$, for each $l$, we get all the possible permutation $\wt\pi\sim \mathcal{P}$,  and these give   all possible  mappings $\wt\si\in \wt \mathcal{B}$. Note that $|\wt\mathcal{B}|= \prod_{l=0}^{K}n_{l}!\leq (2n)!$.

We now use  the notation introduced in the paragraph containing   (\ref{m.5kp}),  and the fact that there are $2^{2n}$ permutations that are compatible with $\PP$, to see that      \bea
      && 
 \sum_{\pi\sim\mathcal{P}} \int \mathcal{T}_{t}( x;\,\pi,e)\prod_{j,k,i}\,dx_{j,k,i}\label{m6a}\\
 &&\qquad=\(4c_{\psi,1}\)^{n}
  \sum_{\wt\si\in\wt\mathcal{B}}F_{t} (\wt\si,0,\ldots, 0)+O\(t^{(2-1/\bb)n-\ep}\).\nn
  \eea
   Since   $|\wt \mathcal{B}|\leq (2n)!$,  we see that the error term only depends on $m$.
 Consider (\ref{m6a}) and the definition of $F_{t} (\wt\si,0,\ldots, 0)$  in (\ref{m.14}) and use
  (\ref{m.4}),  with $m_{j,k}$ replaced by $n_{j,k}$, to see that
        \bea
      && 
 \sum_{\pi\sim\mathcal{P}} \int \mathcal{T}_{t}( x;\,\pi,e)\prod_{j,k,i}\,dx_{j,k,i}\label{m6b}\\
 &&\qquad=\(4c_{\psi,1}\)^{n}
  E\(\prod_{\stackrel{j,k=0}{j< k}}^{K}\(\al_{j,k,t}\)^{
n_{j,k}}\) +O\(t^{(2-1/\bb)n-\ep}\);\nn
  \eea
  ($\al_{j,k,t}$ is defined in (\ref{m.1a})).   

  Recall the definition of $\mathcal{S}$, the set of special pairings, given in the first paragraph of this subsection.
  Since   there are ${( 2n_{j,k})!\over 2^{ n_{j,k}}n_{j,k}!}$ pairings of the $2n_{j,k}$
elements $\{1,\ldots, m_{j,k}\}$, (recall that    $m_{j,k}=2n_{j,k}$), we see that when we sum over all the special pairings we get
 \begin{eqnarray} && 
\sum_{\mathcal{P}\in \mathcal{S}}\sum_{\pi\sim \mathcal{P}}\int \mathcal{T}_{t}( x;\,\pi,e)\prod_{j,k,i}\,dx_{j,k,i}
\label{m.25}\\ 
&&\quad
=  \prod_{\stackrel{j,k=0}{j< k}}^{K}{( 2n_{j,k})!\over 2^{ n_{j,k}}n_{j,k}!}\(4c_{\psi,1}\)^{ n_{j,k}} E\lc\prod_{\stackrel{j,k=0}{j< k}}^{K}\(\al_{j,k,t}\)^{
n_{j,k}}\rc  +O\(t^{(2-1/\bb)n-\ep}\). 
\nn
\end{eqnarray}  
  It follows from (\ref{ord.1})   that the error term,   still,  only depends on $m$.

\medskip	The right-hand side of (\ref{m.25}) is precisely the desired expression in (\ref{m.2}). Therefore, to complete the proof  of Lemma \ref{lem-multiple}, we show that for all the other possible values of $a$, the integral in (\ref{4.25}) can be absorbed in the error term.

\subsection{ { $\bf a=e$}  but not all cycles are   of order two and ${\bf a\ne e}$ }\label{ss3.4}

We show that when  $  a=e$   but not all cycles are   of order two or when ${ a\ne e}$  
\begin{equation}
\Bigg | \int \mathcal{T}_{t}( x;\,\pi ,a )\prod_{j,k,i}\,dx_{j,k,i}\Bigg |=O\(\(t^{2}\psi^{-1} (1/t)\)^{\frac{m}{2} }t^{-\ep}\)\label{f9.50},
\end{equation}
 for some $\ep=\ep_{\bb}>0$.  In this subsection we do not assume that $m$ is even.

 Consider the basic formula  (\ref{m.9}).  Since we    only need an upper bound, we  take absolute values in the integrand and extend   all the time integrals to    $[0,t]$, as we have done several times above.   We refer to this integral as the extended integral. We take the time integrals and get an upper bound for (\ref{m.9}) involving the terms  $u$, $v$ and $w$. As we have done several times above, we choose $m$  of the $u$, $v$ and $w$ terms with arguments that span $R^{m}$. We then bound the remaining $u$, $v$ and $w$ terms  and then make a change of variables and integrate the $u$, $v$ and $w$ terms with the chosen arguments.
Since we want to find the smallest possible upper bound for the extended integral, it is obvious that we  first  integrate as many of the $w$ terms as possible, since such integrals are effectively bounded. (We continue to ignore slowly varying functions of $t$). We then try to integrate as many of the $v$ terms as possible.

  In order to do this efficiently, we divide the $v$ and $w$ terms into sets. As we construct the sets of    $v$ and $w$ terms, we  also choose a  subset $\II  $ of  the $v$ and $w$ terms with arguments  that are linearly independent.  The cardinality of this subset is a lower bound on the number of $v$ and $w$ terms that we can integrate.
 
This is  how we divide the   $v$ and $w$ terms into  sets.  
For each $\pi$ and $ a$ we define  a multigraph $G_{\pi ,a}$ with vertices   
$\{(j,k,i),\,0\leq j<k\leq K,\, 1\leq i\leq m_{j,k}\}$, and  an edge 
between the vertices $\pi_{l } (q-1)$ and $ \pi_{l}(q)$  whenever $( a_{1}(l,q), a_{2}(l,q))=(1,1)$, $l=0,\ldots K$, $2\le q\le m_{l}$, and an edge between the vertices  $\pi_{l } (1)$ and $ \pi_{l-1}(m_{l})$,   whenever $( a_{1}(l,1), a_{2}(l,1))=(1,1)$, $1\le l\le K$.

This graph divides the $w$ terms into cycles and chains.  Assume that there are 
  $S$ cycles.  We denote them by 
   $C_{s}=\{\phi_{s,1},\ldots, \phi_{s,l(s)}\}$, written in cyclic order,  where  the cycle length $l(s)=|C_{s} |\geq 1$ and $1\le s \le S$.  For each  $1\le s \le S$ we take the set of $l(s)$ terms
   \begin{equation}
  \mathcal{G}^{\mbox{\scriptsize cycle}} _{s}=\{w(x_{\phi_{s,2}}-x_{\phi_{s,1}}),\cdots,    w(x_{ \phi_{s,l(s)}}-x_{ \phi_{s,l(s)-1}}), w(x_{\phi_{s,1}}-x_{ \phi_{s,l(s)}}) \} \label{gr.1}.
   \end{equation}
Let 
\be
y_{\phi_{s,i}}=x_{\phi_{s,i}}-x_{\phi_{s,i-1}},\qquad i=2,\ldots,l(s).\label{4.71v}
\ee
  It is easy to see that $\{y_{\phi_{s,i}}\,|  i=2,\ldots, l(s)\}$, are linearly independent. We put the corresponding $w$ terms, $w(x_{\phi_{s,2}}-x_{\phi_{s,1}}),\cdots,    w(x_{ \phi_{s,l(s)}}-x_{ \phi_{s,l(s)-1}})$ into $\mathcal{I}$.
(On  the other hand, since
\be
\sum_{i=2}^{l(s)}y_{\phi_{s,i}}=-(x_{\phi_{s,1}}-x_{ \phi_{s,l(s)}}),\label{4.71s}
\ee   
we see that we can only extract $l(s)-1$ linearly independent variables from the $l(s) $ arguments of $w$ for a given $s$.)

  A cycle of length $1$ consists of a single point  $\phi_{s,1}=\phi_{l(s),1}$ in   the graph, so in this case
   \begin{equation}
  \mathcal{G}^{\mbox{\scriptsize cycle}} _{s}=\{w(0) \} \label{gr.1z}.
   \end{equation}
 We explain below how this can occur.  Obviously, $w(0)$ is not put into $\II$.

Next, suppose there are
  $S'$ chains. We denote them by 
   $C'_{s}=\{\phi'_{s,1},\ldots, \newline \phi'_{s,l'(s)}\}$, written in order,  where $l'(s)=|C'_{s} |\geq 2$ and $1\le s \le S'$.    Note that there are $l'(s)-1$, $w$ terms corresponding to  $C'_{s}$.  
Then for each  $1\le s \le S'$ we form the set of $l'(s)+1$ terms
   \bea
   &&
  \mathcal{G}^{\mbox{\scriptsize chain}}_{s}=\{v(x_{ \phi'_{s,1}}-x_{a(s)}), w(x_{\phi'_{s,2}}-x_{\phi'_{s,1}}),\cdots, \label{gr.2}\\
  &&\hspace{1 in}\cdots,    w(x_{ \phi'_{s,l(s)}}-x_{ \phi'_{s,l(s)-1}}), v(x_{b(s)}-x_{ \phi'_{s,l(s)}}) \}\nn
   \eea
   where $v(x_{ \phi'_{s,1}}-x_{a(s)})$ is the unique $v$ term associated with   $\De_{x_{ \phi'_{s,1}}}^{ 1}$, and similarly, $v(x_{b(s)}-x_{ \phi'_{s,l(s)}}) $ is the unique $v$ term associated with $\De _{ x_{ \phi'_{s,l(s)}}}^{1}$. 
  (This deserves further clarification. There may be other $v$ terms containing the variable  $x_{ \phi'_{s,1}}$  in the extended integral, but there is only one $v$ term   of the form  
 \begin{equation}
   \int_{0}^{t}\Big|\De_{x_{ \phi'_{s,1}}}^{ 1}p_{s}(x_{ \phi'_{s,1}}-u)\Big|\,ds,
   \end{equation}
 where $u$ is some other $x_{\cd}$ variable which we denote by $x_{a(s)}$. This is because one operator $\De_{x_{ \phi'_{s,1}}}^{ 1}$     is associated with $w(x_{\phi'_{s,2}}-x_{\phi'_{s,1}})$ and there are precisely two operators $\De_{x_{ \phi'_{s,1}}}^{ 1}$ in (\ref{m.9})).
   
It is easy to see that variables $ y_{\phi'_{s,i}}=x_{\phi'_{s,i}}-x_{\phi'_{s,i-1}},\,  i=2,\ldots, l(s) $, are linearly independent. 
 We put the  $w$ terms, $w(x_{\phi'_{s,2}}-x_{\phi'_{s,1}}),\cdots,    w(x_{ \phi'_{s,l(s)}}-x_{ \phi'_{s,l(s)-1}})$ into $\mathcal{I}$.  We leave the $v$ terms in $\mathcal{G}^{\mbox{\scriptsize chain}}_{s}$ out of $\mathcal{I}$.

 \medskip	
 At this stage we emphasize that the terms we have put in $\mathcal{I}$ from all cycles and chains have linearly independent arguments. If fact, the set of $x$'s appearing in the different chains and the cycles are disjoint. This is obvious for the cycles and the interior of the chains 
 since there are exactly two difference operators  
 $\De^{x}_{h}$ for each $x$. It also must be true for the endpoints of the chains, since if this is not the case they could be made into larger chains or cycles.
 
For the same reason,   if a $v$  term involving  $\De_{ x'}^{h}$ is not in any of the  sets of chains, then  $x'$ will not appear in the arguments of the terms that  are put in $\mathcal{I}$ from all the cycles and chains.

Suppose, after considering the $w$ terms and the $v$ terms associated with the chains of $w$ terms, that there are $p$ pairs of   $v$ terms left, each pair corresponding to difference operators    $\De^{1}_{z_{j}}$, $j=1,\ldots,p$; ($p$ may be $0$). Let
\be
\ZZ:=\{z_{1},\ldots,z_{p}\}
\ee
A typical $v$ term is of the form  
\begin{equation}
 v^{(j)}(z_{j}-u_{j'}):=v (z_{j}-u_{j'})=\int_{0}^{t}|\De^{h}_{z_{j} } p_{t}(z_{j}-u_{j'})|\,dt,\label{3.93qq}
   \end{equation}
where $u_{j'}$ is some  $x_{\cd}$ term. We use the superscript $(j)$ is  to keep track of the fact that this $v$ term is associated with the  difference operator $\De^{1}_{z_{j}}$. We distinguish between  the variables $z_{j} $ and $u_{j'}$ by referring to $z_{j}$ as a marked variable.   Note that if $u_{j'}$ is also in $\ZZ$, say $u_{j'}=z_{k}$,  then  $u_{j'}$ is also a marked variable but in a different $v$ term. (In this case, in $v^{(k)}(z_{k}-u_{k'})$, where $u_{k'}$ is some other $x_{\cd}$ variable.)

Thus $\ZZ$ is the collection of marked variables.
Consider the corresponding terms
\begin{equation}
  v^{(j)}(z_{j}-u_{j})\quad\mbox{and}\quad   v^{(j)}(z_{j}-v_{j}), \hspace{.2 in}j=1,\ldots,p\label{114.0}
   \end{equation}
  where $u_{j}$ and $v_{j}$ represent whatever terms $x_{\cd}$ and $x'_{\cd}$ are coupled with the two variables $z_{j}$. 
   
 There may be some $j$ for which $u_{j}$ and $v_{j}$ in (\ref{114.0}) are   both in $\ZZ$. Choose such a $j$.
  Suppose  $u_{j}=v_{j}=z_{k}$.  We set  
   \bea
   &&
  \mathcal{G}^{\mbox{\scriptsize $\ZZ,1$}}_{j}=\{  v^{(j)}(z_{j}-z_{k}),   v^{(j)}(z_{j}-z_{k}), \label{gr.3}\\
   &&\hspace{1 in}    v^{(k)}(z_{k}-u_{k}),  v^{(k)}(z_{k}-v_{k}) \}\nn
   \eea
  and put $v^{(j)}(z_{j}-z_{k})$ into $\mathcal{I}$.  Here $u_{k}$ and $v_{k}$ are whatever two variables appear with the two marked variables $z_{k}$.   
 
  On the other hand, suppose $u_{j}$ and $v_{j}$ are   both in $\ZZ$ but $u_{j}= z_{k}$ and $v_{j}= z_{l}$ with $k\neq l$. We set 
 \bea
   &&
  \mathcal{G}^{\mbox{\scriptsize $\ZZ,2$}}_{j}=\{  v^{(j)}(z_{j}-z_{k}),   v^{(j)}(z_{j}-z_{l}), \label{gr.4}\\
  &&\hspace{1 in}    v^{(k)}(z_{k}-u_{k}),  v^{(k)}(z_{k}-v_{k}),v^{(l)}(z_{l}-u_{l}),  v^{(l)}(z_{l}-v_{l}) \}\nn
   \eea
  and put both  $v^{(j)}(z_{j}-z_{k})$ and $v^{(j)}(z_{j}-z_{l})$ into $\mathcal{I}$.

  We then turn to the  elements in $\ZZ$  which have not yet appeared in the arguments of the terms that have  been put into  $\mathcal{I}$.  If there is another $j'$ for which $u_{j'}$ and $v_{j'}$ are   both in $\ZZ$, choose such a $j'$ and proceed as above. If there are no longer any such elements in $\ZZ$, choose some remaining element, say, $z_{i}$.   
 Set 
 \bea
   &&
  \mathcal{G}^{\mbox{\scriptsize $\ZZ,3$}}_{i}=\{  v^{(i)}(z_{i}-u_{i}),   v^{(i)}(z_{i}-v_{i})\} \label{gr.5} 
   \eea
   and if $u_{i}\not\in \ZZ$, place $ v^{(i)}(z_{i}-u_{i})$ into $\mathcal{I}$. If $u_{i} \in \ZZ$, so that 
 $v_{i}\not\in \ZZ$,   place $ v^{(i)}(z_{i}-v_{i})$ into $\mathcal{I}$.    
    
We    continue until we have exhausted $\ZZ$.  

The $v$ and $w$ terms in $\II  $  have  linearly independent arguments. We choose an additional $\II'=m-|\II|  $ terms from the remaining $u$, $v$ and $w$ terms so that the arguments of the $m $ terms are  a spanning set of $R^{m}$.  We  bound the remaining terms by their supremum. We then make a change of variables and   integrate separately  each of the $m$ terms in $\II\cup\II'$. Our goal is to integrate as few $u$ terms as necessary. 

Let $S_{1}$  denote the number of cycles of length $1$. The number of $w$ terms in $\II  $  from cycles is 
  \begin{equation}
 \sum_{s=1}^{S}(l(s) -1 )\geq   {1 \over 2}\sum_{s=1}^{S}l(s)-{S_{1} \over 2}.\label{cyc.1}
  \end{equation}
 
The number of $w$ terms in $\II  $  from  chains is
     \begin{equation}
 \sum_{s=1}^{S'}(l'(s) -1 )\geq   {1 \over 2}\sum_{s=1}^{S'}(l'(s) -1 )+{S' \over 2}.\label{cyc.2}
  \end{equation}

 The number of $w$ terms may be less than  $m$.   (In general it is, but we see below that   it is  possible that the number of $w$ terms may be equal to  $m$.)    Suppose there are $\rho $     terms of type $w$.
Then  the number of $v$ terms  must be 
 $2(m-\rho)$,  and consequently, the number of $u$ terms must be   $\rho$. 
  
  We note that 
  \begin{equation}
 \rho= \sum_{s=1}^{S}l(s)+\sum_{s=1}^{S'}(l'(s) -1 ).\label{cyc.4}
  \end{equation}
 Since the total number of $w$ terms is $\rho$, we see from (\ref{cyc.1})  and (\ref{cyc.2}) that the number of   $w$ terms in $\II$ is at least
 \begin{equation}
 { \rho \over 2}+{S' \over 2}-{S_{1} \over 2}. \label{cyc.5}
  \end{equation}
   This shows that for a given $\rho$ the  the  number of   $w$ terms in $\II$ is minimized when their are no chains.

  We now turn to the number of integrated $v$ terms.  Since the total number of $v$ terms is $2(m-\rho)$, and there  are also two $v$ terms in each set $  \mathcal{G}^{\mbox{\scriptsize chain}}_{s}$ we see that 
 \begin{equation}
2(m-\rho)=\sum_{i,j}| \mathcal{G}^{\mbox{\scriptsize $\mathcal{Z},i$}}_{j}|+2S'.\label{cyc.6}
 \end{equation}
 It is easily seen that we place in $\II  $ at least $1/4$ the number of $v$ terms in the sets $  \mathcal{G}^{\mbox{\scriptsize $\mathcal{Z},i$}}_{j}$
 for all $i,j$.
 Consequently, the number of   $v$ terms with arguments in $\II$ is at least 
 \begin{equation}
{1 \over 4}\sum_{i,j}| \mathcal{G}^{\mbox{\scriptsize $\mathcal{Z},i$}}_{j}|=\frac{m}{2}-{\rho \over 2}- {S' \over 2}.\label{cyc.7}
 \end{equation}
   Combined with (\ref{cyc.5}) we see that the number of  $w$ and $v$ terms in $\II$ is at least
 \begin{equation}
  \( \frac{m}{2}-{\rho \over 2}- {S' \over 2}\)+\({ \rho \over 2}+{S' \over 2}-{S_{1} \over 2}\)=\frac{m}{2}-{S_{1} \over 2}.\label{112}
   \end{equation}
   Since    $\frac{m}{2}-{S_{1} \over 2}$ is an integer, it is at least 
\be
\frac{m}{2}-{S_{1} \over 2}+{\bar 1\over 2},
\ee
where $\bar 1=0$,  if  $m-S_{1} $ is even, and $\bar 1=1$,  if  $m-S_{1} $ is  odd.

Suppose that $\rho \geq \frac{m}{2}+{S_{1} \over 2}-\bar 1/2$. Then, since there are $m$ terms that are integrated, the upper bound of the extended integral will be greatest   if we integrate  
 $\frac{m}{2} +{S_{1} \over 2}-\bar 1/2$ terms of the form  $u$, and bound the remaining $u$ terms  by their supremum.  (Note that  $\frac{m}{2} +{S_{1} \over 2}-\bar 1/2$ is also an integer.) This gives 
 a bound for the   $u$ terms    of
 \begin{equation}
t^{ \frac{m}{2}  +{S_{1} \over 2}-\bar 1/2}\(t \psi^{-1} (1/t)\)^{\rho- \frac{m}{2}  -{S_{1} \over 2}+\bar 1/2 }. \label{cyc.8}
 \end{equation}
   (We  ignore   slowly varying function of $ t$.)
 
Note that  when we integrate $ \frac{m}{2}+{S_{1} \over 2}-\bar 1/2$ terms of the form $u$,   we  only   integrate  $\frac{m}{2}-{S_{1} \over 2}+\bar 1/2$ terms of the form $v$ and $w$. By Lemma \ref{lem-vproprvtj} integrated $v$ terms are much larger than integrated $u$ terms.
What is the maximum number of $v$ terms that can be integrated?  

The  maximum number of $v$ terms that can be integrated  occurs when  all the $w$ terms are in cycles of length $1$ or $2$,  in which case    $(\rho-S_{1})/2$   terms of the form $w$ are integrated. This is easy to see, since in this case the right-hand side of (\ref{cyc.1}) is realized. (We point out in the paragraph containing (\ref{cyc.5}) that to minimize the number of $w$ terms that  are integrated there should be no chains.) 

We are left with  $\frac{m}{2}-{\rho \over 2}+\bar 1/2$ terms of the form $v$ that are integrated. Since, by (\ref{2.4j}), the supremum of the $v$ terms  are effectively bounded this gives a contribution from all $v$ terms of $\(t \psi^{-1} (1/t)\)^{\frac{m}{2}-\rho/2+\bar 1/2 } $. Combining the bounds for $u$ and $v$ terms we obtain 
  \bea
&& t^{\frac{m}{2}+{S_{1} \over 2}-\bar 1/2}\(t \psi^{-1} (1/t)\)^{\rho/2-{S_{1} \over 2}+\bar 1}
\label{erd.4z0} \\
&&\qquad=  \(t ^{2}\psi^{-1} (1/t)\)^{\frac{m}{2}}\(t \psi^{-1} (1/t)\)^{-(\frac{m}{2}-\rho/2)}
 \( \psi^{-1} (1/t)\)^{-{S_{1} \over 2}}\nonumber\\
 &&\hspace{2 in} \(t \(\psi^{-1} (1/t)\)^{2}\)^{\bar 1/2}. \nn
 \eea
   It follows from \cite[(4.77)]{book}   that 
 \be 
  t(\psi^{-1}(1/  t) )^{2}\le C\qquad\forall\,t\ge t_{0}.\label{4.77}
  \ee
  Therefore,  when $\rho \geq \frac{m}{2}+{S_{1} \over 2}-\bar 1/2$, (\ref{erd.4z0})  is bounded by
  \begin{equation}
  C\(t ^{2}\psi^{-1} (1/t)\)^{\frac{m}{2}}\(t \psi^{-1} (1/t)\)^{-(\frac{m}{2}-\rho/2)}
 \( \psi^{-1} (1/t)\)^{-{S_{1} \over 2}} \label{oz.1}.
  \end{equation}
 
  On the other hand,  when   $\rho < \frac{m}{2}+{S_{1} \over 2}-\bar 1/2$,   we  get the largest upper bound for the extended integral when   we integrate all
 $\rho$ of the  $u$ terms.  As above, to get the most $v$ terms integrated,    we only integrate $(\rho-S_{1})/2$ terms of the form $w$. Consequently, since  $m$  terms are integrated,    $m-{3\rho \over 2}+{S_{1} \over 2}$ of the  $v$ terms are integrated. (The remaining  $v$ terms  are  bounded by their supremum, which is effectively bounded.) 
 Combining the bounds for $u$ and $v$ terms we obtain  
 \bea
\lefteqn{ t^{ \rho }\(t \psi^{-1} (1/t)\)^{m-{3\rho \over 2}+{S_{1} \over 2}}\label{erd.5z0}}\\
&&   = \(t ^{2}\psi^{-1} (1/t)\)^{\frac{m}{2}}\(t \(\psi^{-1} (1/t)\)^{3}\)^{ -\rho/2}\(  \psi^{-1} (1/t)\)^{ \frac{m}{2} }\(t \psi^{-1} (1/t)\)^{ {  S_{1} \over 2}}. \nn
 \eea

\medskip	We now show that we  obtain   (\ref{f9.50}) when $S_{1}=0$. Consider the case when $\rho < \frac{m}{2}-\bar 1/2 $ and refer to (\ref{erd.5z0}). Note that  
\bea
   \(t \(\psi^{-1} (1/t)\)^{3}\)^{ -\rho/2}\(  \psi^{-1} (1/t)\)^{ \frac{m}{2}}&=&t^{(3/\bb-1)\rho/2-m/(2\bb)}L(t) \nn\\
   &<&t^{(3/\bb-1)m/4-m/(2\bb)}L(t)\nn\\
     &=&t^{-(\bb-1)m/(4\bb)}L(t),\label{116}
      \eea
where $L(t)$ is slowly varying at infinity and we use the facts that  $(3/\bb)-1>0$ and 
$\rho<\frac{m}{2}$.

Now we consider the case when $\rho \ge \frac{m}{2}-\bar 1/2 $. In dealing with (\ref{oz.1}) we also have that  $\rho \le m-1$, since we arrived at this inequality by assuming that all cycles are of order two, but are excluding the case when the graph $G_{\pi,a}$ consists solely of cycles of order two.  Therefore $\frac{m}{2}-(\rho/2)$ in   (\ref{oz.1}) is strictly positive. This observation and (\ref{116}) gives (\ref{f9.50}) when $S_{1}=0$.

 \medskip	 We now eliminate the restriction that $S_{1}=0$. This requires additional work since the estimates on the right-hand side of (\ref{erd.4z0}) and (\ref{erd.5z0}) are larger in this case. Actually
 we show that the bounds in  (\ref{erd.4z0})  and (\ref{erd.5z0}) that we obtained when $S_{1}=0$, remain the same  when $S_{1}\ne 0$.

  The only way  there can be cycles of length one     is in  terms of the type   
\begin{equation}
\De^{ 1}\De^{- 1}p_{(t-\sum_{q=1}^{m_{l-1}}r_{l-1,q})+r_{l,1}}(x_{\ga_{2\si_{l} (1)-1}}-x_{\ga_{2\si_{l-1} (n_{l-1})-1}})\label{ex.1}
\end{equation}
when $\ga_{2\si_{l} (1)-1}=\ga_{2\si_{l-1} (n_{l-1})-1}$.    In this case
\bea
   \int_{0}^{t}|\De^{ 1}\De^{- 1}p_{s }(x_{\ga_{2\si_{l} (1)-1}}-x_{\ga_{2\si_{l-1} (n_{l-1})-1}})|\,ds&=&w(0,t). \label{3.110}
   \eea
 Note    that  
 \begin{equation}
|\De^{ 1}\De^{- 1}p_{(t-\sum_{q=1}^{m_{l-1}}r_{l-1,q})+r_{l,1}}(0)|=2|\De^{ 1}p_{(t-\sum_{q=1}^{m_{l-1}}r_{l-1,q})+r_{l,1}}(0)|.\label{wo.1}
\end{equation}

This is how we bound the right-hand side of (\ref{f9.50}) when $G_{\pi,a}$ contains cycles of length one.  We return to the  basic formulas   (\ref{m.8}) and (\ref{m.9}). We obtain an upper bound for (\ref{m.9})  by taking the absolute value of the integrand. However, we do not, initially extend  the region of  integration with respect to time.  Instead  we proceed as follows: Let  $l'$ be the largest value of $l$   for which   $\ga_{2\si_{l} (1)-1}=\ga_{2\si_{l-1} (n_{l-1})-1}$.   We extend the integral with respect to $r_{l,q}$   to $[0,t]$  for all $l>l'$, and also  for $l=l'$ and $q>1$, and  bound these integrals with terms of the form $u(\cd,t)$, $v(\cd,t)$ and $w(\cd,t)$. We then consider  the integral of the term in (\ref{wo.1}) with respect to $r_{l',1}$.

Clearly
\be
  \int _{0}^{t}|\Delta^{1}p_{(t-\sum_{q=1}^{m_{l'-1}}r_{l'-1,q})+r_{l',1}}(0)| \,dr_{l',1}\label{ex.48} 
  \leq    \int _{t-\sum_{q=1}^{m_{l'-1}}r_{l'-1,q}}^{2t-\sum_{q=1}^{m_{l'-1}}r_{l'-1,q}}|\Delta^{1}p_{ s}(0)|   \,ds 
\ee
 If $\sum_{q=1}^{m_{l-1}}r_{l-1,q}\leq t/2$, we use (\ref{mac.1}) to bound the left-hand side of (\ref{ex.48})   
 by
 \begin{equation}
\int_{t/2}^{2t} |\De^{1}p_{r }(0)|\,dr\leq Ct\(\psi^{-1} (1/t)\)^{3}\le C\psi^{-1} (1/t),\label{wo.2}
 \end{equation}
 by (\ref{4.77}).
 
Suppose, on the other hand, that $\sum_{q=1}^{m_{l-1}}r_{l-1,q}\geq t/2$. Then for some $q$
we have $r_{l-1,q}\geq t/2m$. We do two things. We bound the contribution of 
 \begin{equation}
\Bigg | \(\(\De^{ 1}_{ x_{ \pi_{l-1}( q)}}\)^{a_{ 1}(l-1,q)}
\(\De^{ 1}_{ x_{ \pi_{l-1}( q-1)}}\)^{a_{ 2}(l-1,q)}\,p_{r_{l-1,q}}(x_{\pi_{l-1}(q)}-x_{\pi_{l-1}(q-1)})\)\Bigg |\label{wo.3}
 \end{equation}
by its supremum over $t/2m\leq r_{l-1,q}\leq t$. 

  To   express we use the notation  
\begin{equation}
\ov u(x,t)= \sup_{t/2m\leq r\leq t}u(x,r),\hspace{.2 in}\ov v(x,t)= \sup_{t/2m\leq r\leq t}v(x,r),\label{4.120}
\end{equation}
and
\begin{equation}
\hspace{.2 in}\ov w(x,t)= \sup_{t/2m\leq r\leq t}w(x,r),\hspace{.2 in}\label{}
\end{equation}
so that the bound of (\ref{wo.3}) may be 
\begin{equation}
   \ov u(x_{\pi_{l-1}(q)}-x_{\pi_{l-1}(q-1)},t),\qquad \ov v(x_{\pi_{l-1}(q)}-x_{\pi_{l-1}(q-1)},t)\label{3.28}
   \end{equation}
or
\[
   \ov w(x_{\pi_{l-1}(q)}-x_{\pi_{l-1}(q-1)},t)
   \]
according to whether there are no, one or two difference operators. 

The terms in (\ref{3.28}) no longer depend on $r_{l-1,q}$ therefore we can
 integrate (\ref{wo.1}) with respect to both 
  $ r_{l,1}$ and $ r_{l-1,q}$ and   use (\ref{wo.4a}) get
 \begin{equation}
\int_{0}^{2t}\int_{0}^{2t} |\De^{1}p_{r+s }(0)|\,dr\,ds \leq C\(t^{2}\(\psi^{-1}(1/t)\)^{3}+L(t)+ 1\),
\label{wo.4}
 \end{equation}
 where $L(t)$ is a slowly varying function at infinity.
 
  Consider how   (\ref{wo.3}) contributes to the bounds in (\ref{erd.4z0}) and (\ref{erd.5z0}). If there are no difference operators they would ultimately contribute either
  \begin{equation}
   \sup_{x}u(x,t)\qquad\mbox{or}\qquad \int u(x,t)\,dx\label{130}
   \end{equation}  
   Now because of the bound in (\ref{3.28}) we get a contribution of 
     \begin{equation}
   \sup_{x}\ov u(x,t)\qquad\mbox{or}\qquad \int \ov u(x,t)\,dx\label{131}
   \end{equation}
   The following table  summarizes results from Lemmas \ref{lem-vproprvtj} and \ref{lem-new3}. It shows that each term in (\ref{130}) is smaller than the corresponding term in (\ref{131}) by a factor of $Ct^{-1}$. Up to    factors  of $\log t$ the same diminution, or more, occurs when we compare   the two functions of $ v(x,t)$   with those of $ \ov v(x,t)$ and the two functions of $ w(x,t)$   with those of $ \ov w(x,t)$.   
   \begin{center}

\begin{tabular}{ |c|| c| c |}
\hline   &&\\
$f(x,t)$ &$\sup_{x}f(x,t)\leq $& $  \int   f(x,t)\,dx\leq $ \\  \hline  \hline &&\\
 $u(x,t)$&$C t\psi^{-1}(1/t)$ &  $t$\\
\hline &&\\
$ \ov u(x,t)$)&$C \psi^{-1}(1/t)$ &  $C$\\
\hline &&\\
$v(x,t)$&  $C\log t $  &  $Ct \psi^{-1}(1/  t)\log t  $\\
\hline &&\\
$ \ov v(x,t)$&  $C  \(\psi^{-1}(1/  t)\)^{2} \leq C/t$
&$C   \psi^{-1}(1/  t)\log t  $\\
\hline && \\
$w(x,t)$&    $C$&$C(\log t)^{2}$\\
   \hline &&\\
$ \ov w(x,t)$&$C  \(\psi^{-1}(1/  t)\)^{3}$& $C\(\psi^{-1}(1/  t)\)^{2}\leq C/t$
\\\hline
\end{tabular}
\end{center}

  To read the table note that the second line  states that   $\sup_{x}u(x,t)\le C \psi^{-1}(1/t)$ and $\int u(x,t)\,dx\le t$, and similarly for the remaining lines.

Combined with (\ref{wo.4}) we see that we have reduced the bounds in (\ref{erd.4z0}) and (\ref{erd.5z0}) by a factor of 
  \begin{equation}
   C\(t \(\psi^{-1}(1/t)\)^{3}+{L(t)+ 1\over t}\)\le C\psi^{-1}(1/t),
   \end{equation}
  where for the last inequality we use   (\ref{4.77}),   as we do in (\ref{wo.2}).

 We  apply a similar procedure for each $l$ in decreasing order, with one exception. If $r_{l-1,q}=r_{l-1,1}$, i.e., $q=1$ in this case and we are also in the 
(\ref{wo.1}) with $l$ is replaced by $l-1$, we skip this term  because this it has already been modified.    We then proceed to deal with remaining terms as we did when we assumed that there were no cycles of length one.  
 
 Consequently,
  if there are $S_{1}$ cycles of length $1$ we have diminished the bounds in  (\ref{erd.4z0}) and (\ref{erd.5z0}) by a factor of at least 
  $ C\(\psi^{-1} (1/t)\)^{{S_{1} \over 2}}$, if $S_{1}$ is even and  by a factor of at least 
  $ C\(\psi^{-1} (1/t)\)^{{S_{1} \over 2}+\frac{1}{2}}$, if $S_{1}$ is odd.
  
In the case of   (\ref{erd.4z0}) we are precisely in the case we considered when $S_{1}=0$, which gives (\ref{f9.50}). In the case of   (\ref{erd.5z0}) the final factor is now $\(t(\psi^{-1}(1/t))^{2}\)^{\frac{S_{1}}{2}}$, which is bounded by a constant by (\ref{4.77}). Thus we are again in the case we considered when $S_{1}=0$, which also gives (\ref{f9.50}).

  \medskip	
 It follows from (\ref{m.25}) and (\ref{f9.50}) that  when $m$ is even
 \begin{eqnarray} && 
\sum_{a}\sum_{\pi_{0},\ldots, \pi_{K}}\int \mathcal{T}_{t}( x;\,\pi,a)\prod_{j,k,i}\,dx_{j,k,i}
\label{m.25haa}\\ 
&&\quad
=  \prod_{\stackrel{j,k=0}{j< k}}^{K}{( 2n_{j,k})!\over 2^{ n_{j,k}}n_{j,k}!}\(4c_{\psi,1}\)^{ n_{j,k}} E\lc\prod_{\stackrel{j,k=0}{j< k}}^{K}\(\al_{j,k,t}\)^{
n_{j,k}}\rc  +O\(t^{(2-1/\bb)n-\ep}\). 
\nn
\end{eqnarray}  
  	 We  now show  that we get the same estimates when  $ \mathcal{T}_{t}( x;\,\pi ,a )$   is replaced by $ \mathcal{T}'_{t}( x;\,\pi ,a )$; (see   (\ref{m.7}) and (\ref{m.9})).

We point out, in the paragraph containing (\ref{3.24q}) that terms of the form  $\De^{1}\De^{-1}p_{\cdot}^{\sharp}$ in (\ref{m.7}) are always of the form $\De^{1}\De^{-1}p_{\cdot}$. Therefore, in showing that (\ref{4.25}) and (\ref{m.8}) have the same asymptotic behavior as $t\to \ff$ we need only consider how the proof of (\ref{m.25haa}) must be modified when the arguments of the density functions with one or no difference operators applied is effected by adding $\pm 1$.

 \label{page24}It is easy to see that the presence of these terms   has no effect on the  integrals that are $O\(\(t^{2}\psi^{-1}(1/t) \)^{n} t^{-\ep}\) $ as $t\to 0$.    This is because in evaluating these expressions we either integrate over all of $R^{1}$ or else use bounds that hold on all of  $R^{1}$.  Since terms with one difference operator only occur in these estimations, we no longer need to be concerned with them.

   Consider the terms with no difference operators applied to them, now denoted by    $p^{\sharp}$. So, for example,  instead of $F (\wt\si, 0,\ldots,0) $ on  the right-hand side of (\ref{m.15sw}), we now have
 \begin{eqnarray} 
&& \int\(\int_{\wt\RR_{0}(0)\times\cdots \times\wt\RR_{K}(0)} \prod_{ l=0}^{ K}\, p^{\sharp}_{(1-\sum_{q=1}^{n_{l-1}}r_{l-1,2q-1}-s_{l-1})+r_{l,1}}\right.
\label{m.14k}\\
&&\left.
\hspace{.2 in}(x_{\wt  \si_{l} (1) }-x_{\wt  \si_{l-1} (n_{l-1}) }) \prod_{q=2}^{n_{l}} p^{\sharp}_{r_{l,2q-1}}(x_{\wt  \si_{l} (q) }-x_{\wt  \si_{l} (q-1) })
\prod_{q=1}^{n_{l}}\,dr_{l,2q-1}\)  \,d x    \nonumber.
\end{eqnarray} 
  Suppose that   $p^{\sharp}_{r}(y_{\si (i)}-y_{\si (i-1)})=p_{r} (y_{\si (i)}-y_{\si (i-1)}\pm 1)$.  We write this term as
\begin{equation}\qquad
p^{\sharp}_{r}(y_{\si (i)}-y_{\si (i-1)})=p_{r} (y_{\si (i)}-y_{\si (i-1)} ) +\De^{\pm  1}p_{r}(y_{\si (i)}-y_{\si (i-1)}). \label{4.52}
   \end{equation}
   Substituting all such terms  into  (\ref{m.14k}) and expanding we get (\ref{m.25haa}) and many other terms with at least one $p_{r} (y_{\si (i)}-y_{\si (i-1)} )$ replaced by $\De^{\pm  1}p_{r} (y_{\si (i)} -y_{\si (i-1)})$.

    Substitute (\ref{4.52}) into (\ref{m.14k}) and write it as the sum of $2^{m }$ terms. One term, which contains no difference operator, is the term we analyzed when we replaced $p^{\sharp}$   by $p$.  All the other terms contain at least one difference operator.  It is easy to see that all these other terms are $O\(\(t^{2}\psi^{-1}\p1/t\)^{n}t^{-\ep}\)$, for some $\ep>0$.

 By  (\ref{4.60}) the term with no difference operators is bounded by $C((t^{2}\newline \psi^{-1}(1/t))^{n} )$. This bound is obtained by extending the integrals to $[0,t]$ in (\ref{m.14av}) and integrating or bounding the resulting terms   $u(\cd,t)$. We are in a situation similar to the one considered in the paragraph containing (\ref{130}). Each difference operator in the  other  terms   replaces  a $u(\cd,t)$ term by  a $v(\cd,t)$ term. By Lemma \ref{lem-vproprvtj} each replacement reduces $C((t^{2}\newline \psi^{-1}(1/t))^{n} )$ by a factor of at least $(t\psi^{-1}(1/t))^{-1}$. Therefore, the replacement of $p$ by $p^{\sharp}$ doesn't change (\ref{m.25haa}) when $m$ is even.
 
\medskip	 We now obtain (\ref{m.2a}).    In     Subsection \ref{ss3.4}  we do not require that  $m$ is even.   Therefore, (\ref{m.2a}) follows from (\ref{f9.50}) unless $G_{\pi,a}$ consists solely of cycles of order two and there are no terms with a single  difference operator.    Therefore, (\ref{m.2a}) follows from (\ref{f9.50}) when $m$ is odd  unless we are in the situation covered in Subsection \ref{assinged}.    This also holds when when $p_{\cdot}$     is replaced  by $p_{\cdot}^{\sharp}$
for the reasons given in the case when $m$ is even. 
 
  However, if any of the  $m_{j,k}$  are odd we can not be in the situation covered in Subsection \ref{assinged}.  Consider  the multigraph $G_{\pi }$ described in the paragraph following  (\ref{multigraph}), with vertices   
$\{(j,k,i),\,0\leq j<k\leq K,\, 1\leq i\leq m_{j,k}\}$, and  an edge 
between the vertices $\pi_{l } (2q-1)$ and $ \pi_{l}(2q)$  for each $ 0\leq l\leq K$ and $1\leq q\leq n_{l}$.  Each vertex is connected to two edges.  Suppose that $  \{(j,k,i)\}=\pi_{l}(2q)$,  with $j=l$ and    $ k=l'\neq l$. Then there is a unique $q'$ such that   $  \pi_{l'}(2q')$ or $\pi_{l'}(2q'-1)$    is also equal to  $ \{(j,k,i)\}$.  

 Suppose    $  \pi_{l'}(2q')=\{(j,k,i)\}$  
and consider $  \pi_{l}(2q-1)$ and $ \pi_{l'}(2q'-1)$   Suppose that $  \pi_{l}(2q-1)=\{(j,k',i')\}$ for some $k'$
and $  \pi_{l'}(2q'-1)=\{(j',k,i'')\}$ for some $j'$. In order that $G_{\pi }$ consist of cycles or order two, we must have   $(j,k',i')=(j',k,i'')$, in particular, $j'=j,k'=k$, (but, of course, $i\neq i'$). This shows that for $G_{\pi }$ to consist of cycles or order two    $m_{j,k}$ must be even for each $j,k$.
 \qed

 \subsection{Proof of Lemma \ref{lem-m4}}\label{sec-m4}

For any $A\subseteq [0,3t]^{n}$ we set 
\begin{eqnarray}
&&F_{A}
 =\int\lc \int_{A} \prod_{ l=0}^{ K}\, p_{r_{l,1}}(x_{\wt  \si_{l} (1) }-x_{\wt  \si_{l-1} (n_{l-1}) }) \right.
\label{vm4.1}\\
&&
\left. \hspace{1 in}\prod_{q=2}^{n_{l}} p_{r_{l,2q-1}}(x_{\wt  \si_{l} (q) }-x_{\wt  \si_{l} (q-1) })
\prod_{ l=0}^{ K}\prod_{q=1}^{n_{l}}\,dr_{l,2q-1}\rc\prod_{q=1}^{n_{l}} \,dx_{\wt  \si_{l} (q) }.    \nn
\end{eqnarray} 
Rather than bound the time integral by that over $[0,3t]^{n}$  as we have in the past, we have to be more careful.

It follows from   (\ref{m.14}), paying special attention to the time variable of    $p_{\cd}$ in the second line,  that 
\begin{equation}
F (\si,s_{0},\ldots, s_{K})=F_{A_{s_{0},\ldots, s_{K}}}\label{vm4.5}
\end{equation}
where
\begin{eqnarray}
 A_{s_{0},\ldots, s_{K}}&=&\lc r\in R_{+}^{n}\,\Bigg| \sum_{\la=0}^{l-1}(t-\sum_{q=1}^{n_{\la}}r_{\la,2q-1}-s_{\la})\leq \sum_{q=1}^{n_{l}}\, r_{l,2q-1}\right.    \label{vm4.6}\\
&&   \left. \hspace{.2 in}\leq \sum_{\la=0}^{l-1}(t-\sum_{q=1}^{n_{\la}}r_{\la,2q-1}-s_{\la})+(t-s_{l});\,l=0,1,\ldots,K\rc.
\nn
\end{eqnarray}
  In particular  
\begin{eqnarray}
 A_{0,\ldots, 0} &=&\lc r\in [0,3t]^{n}\,\Bigg| \sum_{\la=0}^{l-1}(t-\sum_{q=1}^{n_{\la}}r_{\la,2q-1})\leq \sum_{q=1}^{n_{l}}\, r_{l,2q-1}\right.    \label{vm4.6a}\\
&&   \left. \hspace{.2in}\leq \sum_{\la=0}^{l-1}(t-\sum_{q=1}^{n_{\la}}r_{\la,2q-1})+t);\,l=0,1,\ldots,K\rc.
\nn
\end{eqnarray}
Let  $\phi_{l}(r)=\sum_{\la=0}^{l}(t-\sum_{q=1}^{n_{\la}}r_{\la,2q-1})$. We have
\begin{eqnarray}
\lefteqn{ A_{s_{0},\ldots, s_{K}}\De\,A_{0,\ldots, 0}
\label{4.124}}\\
&& \subseteq \bigcup_{l=1}^{K}   \lc r\in [0,3t]^{n}\,\Bigg| \phi_{l-1}(r)-\sum_{\la=0}^{l-1} s_{\la}\leq \sum_{q=1}^{n_{l}}\, r_{l,2q-1}\leq  \phi_{l-1}(r)\rc \nonumber\\
&&\qquad\bigcup_{l=0}^{K}  \lc r\in [0,3t]^{n}\,\Bigg| \phi_{l-1}(r)+t-\sum_{\la=0}^{l} s_{\la}\leq \sum_{q=1}^{n_{l}}\, r_{l,2q-1}\leq  \phi_{l-1}(r)+t\rc. \nonumber\\
&& =\bigcup_{l=1}^{K} \mathcal{A}_{l}\cup \mathcal{B}_{l}\nonumber 
\end{eqnarray}
where, setting $\bar \phi_{l-1}(r)=\phi_{l-1}(r)-\sum_{q=1}^{n_{l}-1}\, r_{l,2q-1}$ we can write 
\begin{equation}
 \mathcal{A}_{l}=\lc r\in [0,3t]^{n}\,\Bigg| \bar \phi_{l-1}(r)-\sum_{\la=0}^{l-1} s_{\la}\leq  \, r_{l,2n_{l}-1}\leq  \bar \phi_{l-1}(r) \rc\label{defa}
\end{equation}
and
\begin{equation}
 \mathcal{B}_{l}=\lc r\in [0,3t]^{n}\,\Bigg| \bar \phi_{l-1}(r)+t-\sum_{\la=0}^{l-1} s_{\la}\leq  \, r_{l,2n_{l}-1}\leq  \bar \phi_{l-1}(r)+t \rc.\label{defb}
\end{equation}
(The first union in (\ref{4.124}) are the points in $A_{s_{0},\ldots, s_{K}}$ that are not in $A_{0,\ldots, 0}$ and the second union are the points in  $A_{0,\ldots, 0}$ that are not in $A_{s_{0},\ldots, s_{K}}$.)

  Note that each time $r_{l,2n_{l}-1}$ is contained in an interval of length $2(K+1)n\sqrt t  $.

We bound each $F_{ \mathcal{A}_{l}}$ and  $F_{ \mathcal{B}_{l}}$ as in (\ref{4.60}) except that we only  integrate with respect to $r_{l,2n_{l}-1}$ over $\AA_{l}$ or $\BB_{l}$. Therefore, instead of getting a bound of $u(x,t)$  or $\int u(x,t)\,dx$ from this term we get a smaller bound.
 
To see this, for fixed $a, b\geq 0$,  let
\begin{equation}
u_{a,b}(x)=\int_{a}^{a+b}p_{s}(x)\,ds.\label{nj.1}
\end{equation}
Clearly
\begin{equation}
\int u_{a,b}(x)\,dx=\int_{a}^{a+b}1 \,ds=b.\label{nj.2}
\end{equation}
In addition, by (\ref{jr.1}),
\bea
&&
\sup_{x} u_{a,b}(x)=\sup_{x} \int_{a}^{a+b}\int e^{ipx} e^{-s\psi(p)}\,dp\,ds \label{nj.3}\\
&&=\int_{a}^{a+b}\int   e^{-s\psi(p)}\,dp\,ds \leq \int_{0}^{b}\int   e^{-s\psi(p)}\,dp\,ds\leq Cb\psi^{-1}(b). 
\nn
\eea

Using (\ref{nj.1})--(\ref{nj.3}) with $b=C\sqrt t$,  and Lemma \ref{lem-vproprvtj} we see that  the bound in (\ref{4.60}) is reduced by a factor of   at least $(t\psi^{-1}(1/t))^{-(\frac{1}{2}-\ep')}$ for any $\ep'>0$. \qed

 \section{Proof of Lemmas \ref{lem-hue.1}--\ref{lem-hue.3}}\label{sec-5}

{\bf  Proof of Lemma \ref{lem-hue.1}  } 
Using the multinomial theorem   on the sum in (\ref{e1.6}) we have  
\begin{eqnarray}
&& E\(\( \wt I_{l,t} \)^{m}\)=\sum_{\wt m\in \MM}\({m! \over \prod_{\stackrel{j,k=0}{j< k}}^{l-1}( m_{j,k}! )}\)E\(\prod_{\stackrel{j,k=0}{j< k}}^{l-1} \(I_{j,k,t/l} \)^{m_{j,k}}\),\,\,
\label{mnt.1}
\end{eqnarray}
 where 
 \[\MM=\lc \wt m=\{m_{j,k}, 0\leq j<k\leq l-1\}\,\Bigg |\,
 \sum_{\stackrel{j,k=0}{j< k}}^{l-1}m_{j,k}=m \rc.\]

We now use Lemma \ref{lem-multiple}, with $t$ replaced by $t/l$ to compute the expectation on the right-hand side of (\ref{mnt.1}).      We get that when all the $m_{j,k}$ are even, there exists an $\ep>0$   such that  
\begin{eqnarray} \lefteqn{ 
E\(\( \wt I_{l,t} \)^{m}\)
\label{mn.2}}\\ 
&&\hspace{-.1in}
= \sum_{\wt m\in \MM}\({m! \over \prod_{\stackrel{j,k=0}{j< k}}^{l-1}( m_{j,k}!) }\) \prod_{\stackrel{j,k=0}{j< k}}^{l-1}{( 2n_{j,k})!\over 2^{ n_{j,k}}(n_{j,k}!)}\( 4c_{\psi ,1}\)^{ n_{j,k}} E \prod_{\stackrel{j,k=0}{j< k}}^{l-1}\(\al_{j,k,t/l}\)^{
n_{j,k}}  \nonumber\\
&&\hspace{2 in} +O(l^{m}(t^{2}\psi^{-1} (1/t))^{n}t^{-\ep}).
\nn
\end{eqnarray}
  (Recall   that when all the $m_{j,k}$ are even, $m_{j,k}=2n_{j,k}$ for all $j$ and $k$ and $n=m/2$.) Here we   use  the fact that
\begin{equation}
\sum_{\wt m\in \MM}\({m! \over \prod_{\stackrel{j,k=0}{j< k}}^{l-1} (m_{j,k}!) }\)=l^{m}\label{mn.3}
\end{equation}
 to compute the error term. It also follows from Lemma \ref{lem-multiple} that 
 \begin{equation}
    E\(\( \wt I_{l,t} \)^{m}\)=O(l^{m}(t^{2}\psi^{-1} (1/t))^{m/2}t^{-\ep})  \label{mn.2x}
   \end{equation}
    if any of the $m_{j,k}$
are odd. 
  (Lemma \ref{lem-multiple} is for a fixed partition  of $m$. In (\ref{mn.2}) and (\ref{mn.2x}) we include the factor $l^{m}$, to account for the number of possible partitions. ) Recall that $l=[\log t]^{q}$ for some  $q>0$.

When $m_{j,k}=2n_{j,k}$ for all $j$ and $k$,  
 \begin{equation}
\({m! \over \prod_{\stackrel{j,k=0}{j< k}}^{l-1} (m_{j,k}! )}\) \prod_{\stackrel{j,k=0}{j< k}}^{l-1}{( 2n_{j,k})!\over 2^{ n_{j,k}}(n_{j,k}!)}={(2n)! \over 2^{n}n!}{n! \over \prod_{\stackrel{j,k=0}{j< k}}^{l-1}(n_{j,k}!)} \label{mn.4}.
 \end{equation}
Using this in (\ref{mn.2}) we get  
 \begin{eqnarray} \lefteqn{
E\(\( \wt I_{l,t} \)^{m}\)
\label{mn.5}}\\ 
&&
= {( 2n)!\over 2^{ n}n!}\(4c_{\psi ,1}\)^{ n}\sum_{\NN}\({n! \over \prod_{\stackrel{j,k=0}{j< k}}^{l-1} n_{j,k}! }\)  E\lc\prod_{\stackrel{j,k=0}{j< k}}^{l-1}\(\al_{j,k,t/l}\)^{
n_{j,k}}\rc \nonumber\\
&&\hspace{2 in} +O(l^{m}(t^{2}\psi^{-1} (1/t))^{n}t^{-\ep}).,
\nn
\end{eqnarray}
where $\NN$ is defined similarly to $\MM$. Using the multinomial theorem as in (\ref{mnt.1}) we see that the sum in (\ref{mn.5}) is equal to $E\lc\(\wt\al_{l,t}\)^{
n}\rc $, which completes the proof of (\ref{e1.8ax}). 
  \qed

\noindent {\bf  Proof of Lemma \ref{lem-hue.2}  }
  By Kac's moment formula
  \begin{eqnarray}
 E\lc\( \al_{t}\)^{
n}\rc&=&E\(\(\int ( L^{ x}_{ t})^{ 2}\,dx\)^{n}\)
\label{e1.60}\\
&=&  2^{n}\sum_{\pi }  \int \int_{\{\sum_{i=1}^{2n}r_{i}\leq t\}}\prod_{i=1}^{2n}
p_{r_{i}}(x_{\pi(i)}-x_{\pi(i-1)})\,\,    \prod_{i=1}^{2n}\,dr_{i}
\,\,    \prod_{i=1}^{n}\,dx_{i} \nonumber,
\end{eqnarray}
where the sum runs over all maps $\pi:\,[1,2n]\mapsto [1,n]$ with $|\pi^{-1}(i)|=2$ for each $i$. The factor $2^{n}$ comes from the fact that we can interchange each  $x_{\pi(i)}$ and $x_{\pi(i-1)}$, $i=1,\ldots,2n$. 
  
 It is not difficult to see that we can find a subset  $J=\{i_{1},\ldots, i_{n}\}\subseteq [1,2n]$, such that each of  $x_{1},\ldots, x_{n}$ can be written as a linear combination of    $y_{j}:=x_{\pi(i_{j})}-x_{\pi(i_{j}-1)}$, $j=1,\ldots ,n$. For  $i\in J^{c}$ we use the bound  $p_{r_{i}}(x_{\pi(i)}-x_{\pi(i-1)})\leq p_{r_{i}}(0)$, then change variables and integrate  out the $y_{j}$, to see that
\begin{eqnarray}
&& \int \(\prod_{i=1}^{2n}
\int_{0}^{t}p_{r_{i}}(x_{\pi(i)}-x_{\pi(i-1)})\,dr_{i}\)
\,\,    \prod_{i=1}^{n}\,dx_{i}
\label{e1.61a}\\
&& \qquad\leq  C\(\int_{0}^{t}p_{r}(0)\,dr \)^{n} \int \(\prod_{i\in J} 
\int_{0}^{t}p_{r_{i}}(x_{\pi(i)}-x_{\pi(i-1)})\,dr_{i}\)
\,\,    \prod_{i=1}^{n}\,dx_{i}\nonumber\\
 &&\qquad =  C u^{n}(0,t)   \(\prod_{i\in J} 
\int \int_{0}^{t}p_{r_{i}}(y_{i})\,dr_{i}\,dy_{i}\)\nn\\
&&\qquad=Cu^{n}(0,t)  \(\int  u(x,t)\,dx\)^{n} \le C\(t^{2}\psi^{-1}(1/t) \)^{n},
\nonumber
\end{eqnarray}
 where we use  (\ref{jr.1j}) and  (\ref{ee1})  for the last line. This shows that
    \begin{equation}
 \| \ \al_{t}\|_{n}\leq C t^{2}\psi^{-1}(1/t)  \label{e1.59z},
  \end{equation}
for all $t$ sufficiently large,   where $C $ depends only on $n$,   and where $\|\cd\|_{n}:=(E(\cd)^{n})^{1/n}$.  
 
It follows from (\ref{e1.59z})  that 
 for  
$l$ sufficiently  large,
  \begin{eqnarray}
  \Big | \| 2\wt\al_{l,t}\|_{n}-\| \al_{t}\|_{n}\Big |&\leq &\| 2\wt\al_{l,t}- \al_{t}\|_{n}=  \| \sum_{ j=0 }^{l-1}\al_{j,j,t/l}\|_{n}\label{4.9aq}\\
  &\le&  l\| \ \al_{0,0,t/l}\|_{n}
  =l\| \ \al_{t/l}\|_{n}\nonumber\\
  &\le&  C t^{2}\frac{\psi^{-1}(l/t)}{l }.\nn 
  \end{eqnarray} 
  We next show that
   when $l=l(t)=[\log t]^{q}$ for any $q>0$,
   \begin{equation}
 \lim_{t\rar \ff} {\Big | \| 2\wt\al_{l,t}\|_{n}-\| \al_{t}\|_{n}\Big | \over t^{2}\psi^{-1} (1/t)}=0.\label{4.9aqr}
   \end{equation} 
  This follows from  (\ref{4.9aq})   since 
  \begin{equation}
  \lim_{t\rar \ff} {\psi^{-1}(l/t) \over l\psi^{-1} (1/t)} =0.  \label{e1.57b}
  \end{equation}
To obtain (\ref{e1.57b}) we use \cite[Theorem 1.5.6, (iii)]{BGT} to see that for all $\de>0$, there exists a $t_{0}$, such that for all $t\ge t_{0}$
  \begin{equation}
   {\psi^{-1}(l/t) \over \psi^{-1} (1/t)}\le l^{(1/\bb)+\de}.\label{pot.1}
   \end{equation}
  Obviously, we pick $\de$ such that $(1/\bb)+\de<1$. 
  
  The statement in (\ref{2.11}) follows from (\ref{4.9aqr}).
\qed

For the  proof of Lemma \ref{lem-hue.3} 
we need the following   lemma which is proved in Section \ref{sec-9}.

\begin{lemma}\label{lem-varep} Under the hypotheses of Theorem \ref{theo-clt2r}\begin{equation}
E\(\int ( L^{ x+1}_{t}- L^{ x}_{ t})^{ 2}\,dx\)=4c_{\psi,0}t+O\(g( t)\)\label{expep}
\end{equation}
as $t\to \ff$, where  
\be
g( t)=     \left \{\begin{array}{ll}   t^{2}\(\psi^{-1}(1/t)\)^{3}&\qquad 3/2<\bb\le2\\\\
L(t)&\qquad\bb=3/2\\\\
C&\qquad1<\bb<3/2 \end{array}  \right.\label{4.9}
  \ee
and   $ L(\cd )$ is  slowly varying at infinity.
Also  
\be \mbox{Var}\(\int ( L^{ x+1}_{t}- L^{ x}_{ t})^{ 2}\,dx\)   \leq  C   t^{2}\psi^{-1}(1/t)  \log  t   .\label{9.3}
\ee 
\end{lemma}

\medskip	 \noindent {\bf  Proof of Lemma \ref{lem-hue.3}  } We prove this lemma by showing that 
\begin{equation}
 \lim_{t\to\ff}{  \sum_{ j=0 }^{l-1}    E( I_{j,j,t/l})-  E\(\int ( L^{ x+1}_{t}- L^{ x}_{ t})^{ 2}\,dx\)\over t(\psi^{-1}(1/t))^{1/2}}=0
   \end{equation}
and 
\begin{equation}
    \lim_{t\to\ff}{  \sum_{ j=0 }^{l-1} \(  I_{j,j,t/l}-  E( I_{j,j,t/l})  \)\over t(\psi^{-1}(1/t))^{1/2}}=0\qquad  \label{4.37a}
   \end{equation}
 in $L^{2}$, where $l=l(t)=[\log t]^{q}$, for some  $q$ sufficiently large.

Set  
 \be
 \phi(t)=t^{2}\psi^{-1}(1/t).\label{5.4}
 \ee
It follows from  (\ref{pot.1}) that 
\begin{eqnarray}
&&{l \phi(t/l)\over \phi(t)}= {\psi^{-1}(l/t) \over l \psi^{-1} (1/t)}\le l^{(1/\bb)+\de-1},
 \\
&&   \nonumber
\end{eqnarray}
  for all  $\de>0$.
Recall $l =[ \log t]^{q}$.  We   choose  a $\de$ and $q<\ff$ such that   
 \begin{equation}
{l(t)\phi(t/l(t))\over \phi(t)}=O\(\frac{1}{\log^{2}t}\),\label{5.5}
   \end{equation}
as $t\to\ff$.

We see from (\ref{expep}) that 
\begin{eqnarray} 
&&  \lim_{t\to\ff}{  \sum_{ j=0 }^{l-1}    E( I_{j,j,t/l})-  E\(\int ( L^{ x+1}_{t}- L^{ x}_{ t})^{ 2}\,dx\)\over t(\psi^{-1}(1/t))^{1/2}}\label{5.12}                                                      \\
&&\qquad  =\lim_{t\to\ff}{l(t)O\(g( t/l(t))\)+O\(g( t)\)\over t(\psi^{-1}(1/t))^{1/2}}=0.\nn
\end{eqnarray}
The last equality follows from the fact that $t(\psi^{-1}(1/t))^{1/2}$ is regularly varying as $t\rar\ff$   with index $1-1/(2\bb)>1/2$, since $\bb>1$, whereas $g( t)$ is regularly varying as $t\rar\ff$ with index $(2-3/\bb)^{+} \leq 1/2$ since, $\bb\leq 2$,   and $l(t)$ is slowly varying.

 Since $I_{j,j,t/l}$ are   independent and identically distributed, we  obtain (\ref{4.37a}) by showing   that 
\begin{equation}
  \lim_{t\to\ff}l(t)\mbox{Var} \,\({I_{j,j,t/l} \over t(\psi^{-1}(1/t))^{1/2}}\)=0. \label{5.16}
  \end{equation}
 Using (\ref{9.3}) and (\ref{5.4}) we see that 
\begin{equation}
 l(t)\mbox{Var} \,\({I_{j,j,t/l} \over t(\psi^{-1}(1/t))^{1/2}}\)=O\({l(t)\phi(t/l(t))\over \phi(t)}\log t\) \label{5.16x}.
  \end{equation}
as $t\to\ff$. Thus (\ref{5.16}) follows from (\ref{5.5}).   \qed

    \section{Proofs of Lemmas \ref{lem-vproprvtj}--\ref{lem-2.5} and \ref{lem-3.2} }\label{sec-lemproofs}

  Since the  L\'evy processes, $X$, that we are concerned with satisfy (\ref{regcond2}), it follows from the Riemann Lebesgue Lemma that they have transition probability density functions, which we designate as $p_{s}(\cd)$. Taking the inverse Fourier transform of the characteristic function $X_{s}$, and using the symmetry of $\psi$, we see that 
 \bea
p_{s }(x)&=&{1 \over    2 \pi}\int e^{ipx} \, e^{-s\psi (p) }\,dp\label{69dd} \\
&=&{1 \over     \pi}\int_{0}^{\ff}\cos (px)\, e^{-s\psi (p) }\,dp.\nn 
\eea

Our basic hypothesis is that $\psi(\la)$  is regularly varying at 0 with   index $1<\bb\leq 2$. Therefore $\psi(\cd)$ is    asymptotic  to an increasing function near zero. Considering the way we use $\psi(\cd)$ in the estimates below, we can assume that $\psi(\la)$ is strictly increasing for $0\le \la\le \la_{0}$, for some $\la_{0}>0$, and that $\psi^{-1}(\la)$ is well defined for $0\le \la\le \la_{0}$. Actually, we are really interested in  $\psi^{-1}(1/s)$ as $s\to\ff$. Therefore, there exists an $s_{0}$ such that $\psi^{-1}(1/s)$, as a function of $s$ is regularly varying with index $- 1/\bb$ for   $s\ge s_{0}$.
 
\medskip	 The next two lemmas give fundamental estimates that are used in the proofs of the lemmas in Section \ref{sec-indf}.
 
 \begin{lemma}\label{lem-new1}  Let    $X$ be a symmetric L\'{e}vy process 
with L\'{e}vy exponent
 $\psi(\la)$ that is regularly varying at $0$ with index $1<\bb\leq 2$ and satisfies (\ref{regcond2})--(\ref{1.12}).
   Then for all  $\ga\ge 1$ and for all $s$ sufficiently large and all $x\in R^{1}$, 
    \bea
   p_{s}(x)&\le &C\( \psi^{-1}(1/s) \wedge\frac{1}{\psi^{-1}(1/s)x^{2}} \) \label{8.15a};\\
   |\De^{\ga}p_{s}(x)|&\le&  C\ga^{2}\( \(\psi^{-1}(1/  s)\)^{2}\wedge  \frac{1+\log^{+} |x|}{ x^{2}}\);\label{151}\\
 |\De^{\ga} \De^{-\ga}p_{s}(x)|&\le& C \ga^{2}\( \(\psi^{-1}(1/  s)\)^{3}\wedge \frac{\psi^{-1}(1/  s)}{ x^{2}}\)\label{132}. 
\eea
 \end{lemma}
 
 \begin{lemma}\label{lem-vproprvt} Let    $X$ be a symmetric L\'{e}vy process 
with L\'{e}vy exponent
 $\psi(\la)$ that is regularly varying at $0$ with index $1<\bb\leq 2$ and satisfies (\ref{regcond2})--(\ref{1.12}).
   Then for all $t$ sufficiently large and all $x\in R^{1}$  
   \begin{equation}
u(x,t):= \int_{0}^{t}\, \,p_{s }(x)\,ds\le C\(t\psi^{-1} (1/t)\wedge \frac{t(1+\log^{+} |x|)}{|x|}\)\label{jr.1},
 \ee
   \begin{equation}
 v_{\ga}(x,t):=\int_{0}^{t}\, |\De ^{ \ga}\,p_{s }(x)|\,ds\leq C\ga^{2} \(\log t\wedge{ t\psi^{-1}(1/t)\over |x|}\wedge t  {1+  \log^{+} |x|\over x^{2}}\)  \label{2.4}
 \ee
and
 \begin{equation}
w_{\ga}(x,t):=\int_{0}^{t}\,|\De^{ \ga}\De^{ -\ga} \,p_{s }(x)|\,ds \leq C\ga^{2} \(1\wedge { \log t\over |x| }\wedge {t\psi^{-1} (1/t)\over |x|^{2}}\).\label{jrst1.3y}
\end{equation}
 \end{lemma}

 We use the following lemma repeatedly.
 
 \begin{lemma} \label{lem-6.3}For all   $p\in R^{1}$ and $s ,q>0 $, 
\be
e^{-s\psi(p)}\le \frac{C}{s^{q}\psi^{q}(p)}.\label{bound}
\ee
 \end{lemma} 
 
 \Proof This is elementary since   for all $q>0$
\begin{equation}
e^{-s\psi(p)}\le \frac{\sup _{s\ge 0}s^{q}\psi^{q}(p)e^{-s\psi(p)} }{s^{q}\psi^{q}(p)}.
   \end{equation}\qed

    \noindent{\bf  Proof of Lemma \ref{lem-new1} } We first note that by (\ref{bound})  with $q=1$, and (\ref{regcond2})    \bea
\lefteqn  {\int_{0}^{\ff} e^{-s\psi(p)}\,dp \label{8.15new}}\\
 && \le \(\int_{0}^{\psi^{-1}(1/s)} e^{-s\psi(p)}\,dp+ \int_{\psi^{-1}(1/s)}^{1} e^{-s\psi(p)}\,dp+ \int_{1}^{\ff} e^{-s\psi(p)}\,dp\)\nn\\
& &\le   \(\psi^{-1}(1/s) +{1 \over   s  }\int_ {\psi^{-1}(1/s)}^{1} \frac{1}{\psi  (p)}\,  \,dp+{1 \over   s  }\int_  {1}^{\ff} \frac{1}{\psi  (p)}\,  \,dp\nn\)\\
&& \le C\(\psi^{-1}(1/s)+\frac{1}{s}\) \le C \psi^{-1}(1/s),    \nn
\eea
for all $s$ sufficiently large. 
Therefore, it follows from (\ref{69dd}),    that for all $s$ sufficiently large
 \begin{equation}
   p_{s}(x)\le C\(  \psi^{-1}(1/s)  \)\label{8.15}.
   \end{equation}

  By integration by parts  
 \bea 
  p_{s }(x) 
  & =&{1 \over     \pi x}   \int_{0}^{\ff}e^{-s\psi (p) } \,d(\sin px) \label{8.19w} \\
  & = &-{ 1 \over     \pi x}  \int_{0}^{\ff}\sin px\,\left (\frac{d}{dp}e^{-s\psi (p) }\right )\,dp \nn\\
    & =&-\frac{ 1}{\pi x^{2}}  \int_{0}^{\ff}\cos px\left (\frac{d^{2}}{dp^{2}}e^{-s\psi (p) }\right ) \,dp.  \nn
 \eea
 where the last line uses   the  fact  that $\psi'(0)=0$,  which follows from  (\ref{regcond})  and the first inequality in (\ref{88.m}).
 
We have
 \bea
    \frac{d^{2}}{dp^{2}} e^{-s\psi (p) } &=& \(s^{2}(\psi'(p))^{2}-s\psi''(p)\)e^{-s\psi (p) }  \label{118}
   \eea
 By (\ref{bound}) and (\ref{88.m}) for $p\le 1$
   \begin{equation}
	   s^{2}(\psi'(p))^{2}e^{-s\psi (p) }\le C{s(\psi'(p))^{2}\over \psi(p) }\le C{s\,\psi(p) \over p^{2} }
   \end{equation}
   and 
   \begin{equation}
   s|\psi''(p) |e^{-s\psi (p)} \le  C{s\,\psi(p) \over p^{2} }.
   \end{equation} 
Therefore, for all $s$ sufficiently large 
 \bea
    \bigg| \int_{0}^{\psi^{-1}(1/s)}\cos px\left (\frac{d^{2}}{dp^{2}}  e^{-s\psi (p) } \right ) \,dp   \bigg|&\le& Cs \int_{0}^{\psi^{-1}(1/s)} { \psi(p) \over p^{2} }\,dp\nn\\
    &\le&{ C \over \psi^{-1}(1/s)}\label{7.17} .
   \eea
   By (\ref{bound}), (\ref{118}) and (\ref{88.m}), for $p\le 1$  
     \bea
    \bigg| \frac{d^{2}}{dp^{2}} e^{-s\psi (p) }  \bigg|    &\le&C \left\{\(\frac{\psi'(p)}{\psi(p) }\)^{2}+\frac{|\psi''(p)|}{\psi (p)}\right\}\le \frac{C}{p^{2}}.\label{7.17f}
   \eea
    Therefore, for all $s$ sufficiently large
 \bea
    \bigg| \int_ {\psi^{-1}(1/s)}^{1}\cos px\left (\frac{d^{2}}{dp^{2}}  e^{-s\psi (p) } \right ) \,dp   \bigg|&\le& C   \bigg| \int_ {\psi^{-1}(1/s)}^{1}{ 1 \over p^{2} }\,dp\nn\\
    &\le&{ C \over \psi^{-1}(1/s)}\label{7.17w} .
   \eea
 By (\ref{7.17f})  and (\ref{1.12}) 
  \be \bigg| \int_  {1}^{\ff}\cos px\left (\frac{d^{2}}{dp^{2}}  e^{-s\psi (p) } \right ) \,dp   \bigg| \le C  \int_  {1}^{\ff}  \left\{\(\frac{\psi'(p)}{\psi(p) }\)^{2}+\frac{|\psi''(p)|}{\psi (p)}\right\}  \,dp\le C. \label{124}
    \ee 
   Using (\ref{8.15}), (\ref{8.19w}), (\ref{7.17}), (\ref{7.17w})  and  (\ref{124}) we get (\ref{8.15a}). 
 
  \medskip	We now obtain (\ref{132}).  
\bea
 \De^{ \ga}\De^{ - \ga} p_{s }(x) &= &2p_{s }(x)-p_{s }(x+ \ga)-p_{s }(x- \ga)\label{72dd}\\
  &=&{4 \over  \pi}\int_{0}^{\ff}  \cos (px) \sin^{2}(p\ga   /2) \,e^{-s\psi (p) }\,dp.\nn
   \eea  
 Therefore, by (\ref{bound})
 \bea
\lefteqn{| \De^{\ga}\De^{ - \ga} p_{s }(x) |}\label{134}\\
&&\le  \nn C\int_{0}^{\ff}    \sin^{2}(p\ga   /2) \,e^{-s\psi (p) }\,dp\\
&&\le C\ga^{2}\(\int_{0}^{\psi^{-1}(1/  s)}  p^{2} \,dp\,+\frac{1}{s^{3}}\int_ {\psi^{-1}(1/  s)} ^{1} {p^{2}\over \psi^{3}(p) }\,dp+\frac{1}{s^{3}}\int_  {1}^{\ff} {p^{2}\over \psi^{3}(p) }\,dp\)\nn\\
&&\le C \ga^{2}\(\psi^{-1}(1/  s)\)^{3}.\nn
   \eea
 We next show that 
    \begin{equation}
\De^{\ga}\De^{ -\ga} p_{s }(x)=\frac{8}{\pi}{K_{\ga}(s,x) \over x^{2} }\label{89.1dd}
\end{equation}
where  
\begin{equation}
 K_{\ga}(s,x):= \int_{0}^{\ff}\sin^{2}( px/2)\( \sin^{2}(p\ga/2)\, e^{-s\psi (p) }\)''\,dp.\label{89.2dd}
\end{equation}
 To get this we integrate by parts in (\ref{72dd}),  
\bea
  &&   \int_{0}^{\ff}  \cos p x \sin^{2}(p\ga /2) \, e^{-s\psi (p) }\,dp \label{89dd}\\
  &&\qquad= \frac{1}{x}   \int_{0}^{\ff} \sin^{2} (p\ga /2) \, e^{-s\psi (p) }\,d(\sin px)\nn \\
   &&\qquad= -\frac{1}{x}   \int_{0}^{\ff}\sin px\( \sin^{2}(p \ga/2)\, e^{-s\psi (p) }\)'\,dp\nn\\
   &&\qquad=  -\frac{1}{x}   \int_{0}^{\ff} \( \sin^{2}(p\ga /2)\, e^{-s\psi (p)} \)'\,d\( \int_{0}^{p}\sin  rx\,dr\)\nn\\
       &&\qquad=  -\frac{1}{x^{2}}   \int_{0}^{\ff} \( \sin^{2}(p\ga /2) \, e^{-s\psi (p) }\)'\,d\(1-\cos  px \)\nn\\
         &&\qquad= \frac{2}{x^{2}}   \int_{0}^{\ff}\sin^{2}( px/2)\( \sin^{2}(p\ga /2)\, e^{-s\psi (p) }\)''\,dp\nn.
   \eea
which gives  (\ref{89.1dd}).
  
 Let $g(p)= e^{-s\psi (p) }$ and note that 
 \begin{equation}
 \( 2\sin^{2} (p \ga/2)\, e^{-s\psi (p) }\)'=  \ga g(p)\sin (p \ga)   +2g'(p)\sin^{2} (p\ga /2)\label{90}
   \end{equation}
and  
 \be 
 \( 2\sin^{2} (p/2)\, e^{-s\psi (p) }\)''\label{91dd} =  \ga^{2}g(p)\cos (p \ga)  +2\ga g'(p)\sin (p \ga)  +2g''(p)\sin^{2} (p\ga /2). 
   \ee
 Substituting (\ref{91dd}) in (\ref{89.2dd})
 we    write $K_{\ga}(s,x)= I+ II + III$. 
 
 Note that  
      \bea 
  | I| &= &\ga^{2}  \bigg| \int_{0}^{\ff}\cos (p\ga)  \sin^{2}( px/2)e^{-s\psi (p) }\,dp \bigg|  \label{88.12}\\
    &\le& \ga^{2}  \int_{0}^{\ff}  e^{-s\psi (p) }\,dp\le C\ga^{2}\psi^{-1}(1/s)  \nn
    \eea
    by (\ref{8.15new}).
    
By (\ref{91dd})  
    \bea
   | II|  &=& 2\ga   \Big |\int_{0}^{\ff}\sin (p\ga)  \sin^{2}( px/2)g'(p)\,dp\Big | \label{88.13dd} \\
    &\le& C \ga^{2} \int_{0}^{\psi^{-1} (1/s)}ps| \psi '(p)|\,e^{-s\psi (p)  }\,dp \nn \\
     & & \quad +C\ga^{2} \int_{\psi^{-1} (1/s)}^{1}ps| \psi '(p)|\,e^{-s\psi (p)  }\,dp\nn \\
     &&\qquad+C\ga \nn\int_{1}^{\ff}s| \psi '(p)|\,e^{-s\psi (p)  }\,dp\nn
     \eea
 By (\ref{88.m}) and (\ref{bound}) the first of these last three integrals
 \begin{equation}
   \le  C\int_{0}^{\psi^{-1} (1/s)} s \psi (p)\,e^{-s\psi (p)  }\,dp\le  C\int_{0}^{\psi^{-1} (1/s)}\,dp\le   C \psi^{-1} (1/s).
   \end{equation}    
   By (\ref{88.m}) and (\ref{bound}) the second of the  last three integrals in (\ref{88.13dd})
 \begin{equation}
   \le   \frac{C}{s} \int_{\psi^{-1} (1/s)}^{1} {s^{2} \psi^{2} (p)\over \psi(p)}\,e^{-s\psi (p)  }\,dp \le   \frac{C}{s} \int_{\psi^{-1} (1/s)}^{1}{dp\over \psi(p)}\le   C \psi^{-1} (1/s).
   \end{equation}   
     By (\ref{1.12}) and (\ref{bound}) the third of the  last three integrals in (\ref{88.13dd})
   \begin{equation}
 \frac{1}{s}\int_{1}^{\ff}s^{2}\psi ^{2}(p){| \psi '(p)|\over \psi ^{2}(p)} \,e^{-s\psi (p)  }\,dp\le  \frac{1}{s}\int_{1}^{\ff} {| \psi '(p)|\over \psi ^{2}(p)} \,dp\le \frac{C}{s}.
 \end{equation} 
  Since $1/s<\psi^{-1} (1/s)$ for all $s$   sufficiently large, and $\ga\ge1$, we see that
      \be 
   |II|\le C\ga^{2}\psi^{-1}(1/s)\qquad\forall\, x\in R^{1}\label{143aa}.
 \ee 
 
  Similarly,   
        \bea 
 | III | &=& 2  \Big |\int_{0}^{\ff}\sin^{2} (p \ga) \sin^{2}( px/2)g''(p)\,dp  \Big|    \label{151dd}\\
     &\le &C \ga^{2} \int_{0}^{\psi^{-1}(1/s)}p^{2}  \(s |\psi''(p)|   + s^{2} |\psi'(p)|^{2}  \)   e^{-s \psi (p) }  \,dp  
   \nn\\ 
   &  &\quad+C\ga^{2} \int_{\psi^{-1}(1/s)}^{1}  p^{2}\(  s |\psi''(p)|   + s^{2} |\psi'(p)|^{2}  \)e^{-s \psi (p) }  \,dp  
   \nn\\ 
  &  &\qquad+C \int_1^{\ff}  \(  s |\psi''(p)|   + s^{2} |\psi'(p)|^{2}  \)e^{-s \psi (p) }   \,dp\nn  \\
   &\le &C  \ga^{2}\int_{0}^{\psi^{-1}(1/s)}   \(s  \psi (p)   + s^{2} \psi^{2} (p)   \)   e^{-s \psi (p) }  \,dp  
   \nn\\ 
   &  &\quad+{C\ga^{2}\over s} \int_{\psi^{-1}(1/s)}^{1} \frac{1}{\psi(p)}  \(s^{2}  \psi^{2} (p)   + s^{3} \psi ^{3}(p)   \) e^{-s \psi (p) }  \,dp  
   \nn\\ 
   &  &\qquad+\frac{C}{s} \int_1^{\ff}  \(  s^{2}\psi^{2}(p) {|\psi''(p)| \over \psi^{2}(p)}  + s^{3}\psi^{3}(p) {|\psi'(p)|^{2}\over \psi^{3}(p)}  \)e^{-s \psi (p) }   \,dp\nn  \\
   &\le& C\ga^{2}\psi^{-1}(1/s).\nn
 \eea
   Note that  $\lim_{\la\to\ff} \psi(\la)=\ff$, see e.g. \cite[Lemma 4.2.2]{book},  so that (\ref{1.12}) implies that
\begin{equation}
   \int_{1}^{\ff}\frac{|\psi'(\la)| ^{2}}{\psi^{3}(\la)}\,d\la<\ff,\quad  \int_{1}^{\ff}\frac{|\psi''(\la)|}{\psi^{2} (\la)}\,d\la<\ff.\label{1.12aa}
   \end{equation}
   We use this to bound the next to last line in (\ref{151dd}).

Combining (\ref{134}),    (\ref{89.1dd}),   (\ref{88.12}),  (\ref{143aa}), and (\ref{151dd}) we get (\ref{132}). 

\medskip	 
 We now obtain (\ref{151}). Note that  
\bea
  \De^{\ga} p_{s}(x)&= & p_{s }(x+\ga)-p_{s }(x )\label{70}
\\
  &=&{1\over  \pi}\int_{0}^{\ff} \(\cos p(x+ \ga)- \cos px\)\, e^{-s\psi (p) }\,dp \nn\\
   &=&-{2 \over  \pi}\int_{0}^{\ff}   \cos (px) \sin^{2}(p\ga /2)e^{-s\psi (p) } \nn\\
   &&\hspace{.5 in}-\frac{ 1}{\pi}\int_{0}^{\ff}   \sin (px) \sin (p\ga )\,  e^{-s\psi (p) }\,dp \nn   
    \eea
 Thus
\be 
 \De^{  \ga }p_{s }(x)  =  -{1 \over 2}\De^{  \ga}\De^{ - \ga} p_{s }(x)-\frac{ 1} {\pi}\int_{0}^{\ff}   \sin (px ) \sin   (p\ga)  \, e^{-s\psi (p) }\,dp. \label{6.10}
 \ee
 The second order difference is bounded in (\ref{132}). We deal with the 
 second integral which is bounded by  
 \bea
\lefteqn{  \int_{0}^{\ff}  |  \sin (p\ga) | \, e^{-s\psi (p) }\,dp\label{612}}\\
 &&\le\ga  \(\int_{0}^{\psi^{-1}(1/s)} p \,dp+ \int_{\psi^{-1}(1/s)}^{1} pe^{-s\psi(p)}\,dp+ \int_{1}^{\ff}  e^{-s\psi(p)}\,dp\)\nn\\
 &&\le   C \ga\(\(\psi^{-1}(1/s)\)^{2}+{1 \over   s^{2}  }\int_ {\psi^{-1}(1/s)}^{1} \frac{p}{\psi^{2}  (p)}\,  \,dp+{1 \over   s^{2}  }\int_  {1}^{\ff} \frac{1}{\psi^{2}  (p)}\,  \,dp\nn\)\\
 &&\le  C\ga\(\psi^{-1}(1/s)\)^{2}  \nn.
\eea
 This gives us the first bound in (\ref{151}). To obtain the second bound we
integrate by parts twice to get 
  \bea
  &&   \int_{0}^{\ff} \sin (p x) \sin  (p\ga)  \, e^{-s\psi (p) }\,dp \label{89s}\\
  &&\qquad=- \frac{1}{x}   \int_{0}^{\ff} \sin  (p\ga)   \, e^{-s\psi (p) }\,d(\cos px)\nn \\
   &&\qquad = \frac{1}{x}   \int_{0}^{\ff}\cos (px)\(\sin  (p\ga)   \, e^{-s\psi (p) }\)'\,dp\nn\\
   &&\qquad=  \frac{1}{x^{2}}   \int_{0}^{\ff} \(\sin  (p\ga)  \, e^{-s\psi (p) }\)'\,d\(  \sin  px \)\nn\\
         &&\qquad= -\frac{1}{x^{2}}   \int_{0}^{\ff}\sin( px)\(\sin  (p\ga) \, e^{-s\psi (p) }\)''\,dp.\nn\\
          &&\qquad :=   \frac{G}{x^{2}}  .\nn
   \eea
  
 Since
 \bea
&& \(\sin  (p\ga)   \, e^{-s\psi (p) }\)'' 
 =\( -\ga^{2} \sin  (p\ga)  -2s \ga \cos (p\ga)  \,\psi'(p)\right.\\
 &&\qquad\hspace{1.2in}\left.-\sin (p\ga)  (s \,\psi''(p)-s^{2}(\psi'(p))^{2}\)e^{-s\psi (p) }\nn,
   \eea
 we can write 
 \be G= G_{1}+G_{2}+G_{3}, \label{550} \ee          
  where 
  \begin{eqnarray}
 |G_{1}|&=& \ga^{2}  \Big |\int_{0}^{\ff}\sin (px )\sin (p\ga)  \, e^{-s\psi (p) } \,dp\Big | \label{89.2edd}\\
  &  
  \leq &C \ga^{2} \(\psi^{-1} (1/s)\)   \nonumber,
  \end{eqnarray}   
  for all $s$ sufficiently large, by (\ref{8.15new}).
  
   Using   (\ref{bound}),   (\ref{88.m})  and (\ref{1.12}), we see that  
     \begin{eqnarray}
  |G_{2}|&=&2  \ga\Big |\int_{0}^{\ff}\sin px \cos(p\ga)  \( \psi'(p)\,s e^{-s\psi (p)}\) \,dp 
  \label{89.2fdd}\\
  &\le&   C\ga  \int_{0}^{1 }  |\sin px|\, |\psi'(p)| s e^{-s\psi (p)}  \,dp + C\ga   \int_  {1}^{\ff}  |\psi'(p)|se^{-s\psi (p)}  \,dp \nn\\
   &\le&   C\ga  \int_{0}^{1} {|\sin px| \over p}   \,dp 
     + {C\ga\over s}   \int_  {1}^{\ff}  {|\psi'(p)|\over \psi^{2}(p)}   \,dp \nn\\
  &\le& C\ga\(1+\log^{+} x+(1/s)\) \nn,
  \end{eqnarray} 
 where we use  
\begin{equation}
 \int_{0}^{1} {|\sin px| \over p}   \,dp= \int_{0}^{|x|} {|\sin p| \over p}   \,dp \le    C\ga \(1+\log^{+} |x| \).
\label{8.15d}
\end{equation}  
Therefore, for $s$ sufficiently large  
\begin{equation}
  |G_{2}|\le    C\ga \(1+\log^{+} |x| \).
   \end{equation}
       Similarly,  
   \begin{eqnarray}
 |G_{3}|&=&  \Big |\int_{0}^{\ff}\sin px\,\, \sin p \ga \( s\psi''(p)   - s^{2} (\psi'(p))^{2}    \)e^{-s\psi (p) } \,dp\Big | 
  \nn\\
    &\leq&  \label{89.2gdd}C\ga    \int_{0}^{ 1}|\sin px|\,p \( s|\psi''(p)| + s^{2} (\psi'(p))^{2}   \)e^{-s\psi (p) } \,dp
 \\
    &&\qquad+ C \int_{1}^{ \ff}  \( s|\psi''(p)| + s^{2} (\psi'(p))^{2}   \)e^{-s\psi (p) }    \,dp
  \nn\\
  &\le&  C  \ga   \int_{0}^{1} |\sin px| \(  {s\psi(p)|\over p} + {s^{2}  \psi^{2}(p)\over p}    \)  e^{-s\psi (p) }  \,dp\nn\\
  &&\qquad +C \int_{1}^{ \ff} \({|\psi''(p)|\over \psi(p) } +  {(\psi'(p))^{2} \over \psi^{2}(p) }\)\,dp \nn\\
  &\le&   C \ga \(\int_{0}^{1} {|\sin p x| \over p}   \,dp  +C\)\le C\ga \(1+\log^{+} |x|\). 
  \nonumber  \end{eqnarray}
Combining (\ref{89.2edd})--(\ref{89.2gdd}) we get the second bound in (\ref{151}).
\qed

   We use the next two lemmas in the proof of Lemma \ref{lem-vproprvt}.

\begin{lemma}  Let    $X $ be a symmetric L\'{e}vy process 
with  L\'{e}vy exponent
 $\psi(\la)$  that is regularly varying at zero with index $1<\bb\leq 2$ and satisfies (\ref{regcond2}).
 Then for any $r\geq 0$   and $t>0$, \begin{equation} 
 \int_{0}^{t} s^{r}e^{-s\psi (p)} \,ds\leq C_{k}\(t\wedge {1 \over \psi  (p)}\)^{r+1} ; \label{el.1}
\end{equation}
for all $t\ge0$, where $C_{k}<\ff$, is a constant depending on $k$. Furthermore, for any $r\geq 0$   and all $t$ sufficiently large,  
\be 
 \int_{0}^{\ff} \psi^{r} (p) \int_{0}^{t} s^{r}e^{-s\psi (p)} \,ds\,dp 
 \le Ct\psi^{-1} (1/t).\label{jr.0}
  \ee 
\end{lemma}

\Proof
 The first part of the bound in  the first inequality in  (\ref{el.1}) comes from taking  $e^{-s\psi (p)}\leq 1$;  the second from letting $t=\ff$.  
  
  Since
 \begin{equation}
  \psi^{r} (p)  s^{r}e^{-s\psi (p)}= 2^{r}   \psi^{r} (p)  \({s\over 2}\)^{r} e^{-s\psi (p)/2}e^{-s\psi (p)/2},
   \end{equation}it follows from  (\ref{bound}) and (\ref{8.15new}) that
 \begin{equation}
\int_{0}^{\ff} \psi^{r} (p)  s^{r}e^{-s\psi (p)} \,dp\leq C\int_{0}^{\ff}   e^{-s\psi (p)/2} \,dp \leq C\psi^{-1} (1/s)\label{sr.j1}
 \end{equation}
  for all $s$ sufficiently large. On the other hand
 for any fixed $t_{0}$,  
  \begin{equation}
  \int_{0}^{\ff}   \int_{0}^{t_{0}}  e^{-s\psi (p)/2} \,ds\,dp= 2\int_{0}^{\ff}    {1-e^{-t_{0}\psi (p)/2} \over \psi (p)} \,dp<\ff,\label{1415}
  \end{equation}
  by  (\ref{regcond2}).  Putting these two together, and using the fact that  $\psi^{-1} (1/t)$ is regularly varying at infinity, gives (\ref{jr.0}).   \qed

\begin{lemma} Under the hypotheses of Theorem \ref{theo-clt2r}, for $r=0,1,\dots$
  \bea 
  \int_{0}^{1}{|\sin px| \over p}\psi^{r}(p)\(t\wedge {1 \over \psi (p)}\)^{r+1} \,dp\label{6.36} &\le&   Ct(1+\log^{+}|x|) ; \\
    \int_{0}^{1} \psi^{r}(p)\(t\wedge {1 \over \psi (p)}\)^{r+1} \,dp  &\le&  Ct\psi^{-1} (1/t)\label{6.7};\\
     \int_{0}^{t}\int_{0}^{\ff}  |    \sin  p\ga |\,  e^{-s\psi (p) }\,dp\,ds&\le& C \ga \log t,\label{6.7a}
   \eea
for all $t$ sufficiently large.

 \end{lemma}

\Proof
We first note that  for $r=0,1,\dots$
   \begin{equation}
   \psi^{r+1}(p)\(t\wedge {1 \over \psi (p)}\)^{r+2} \le \psi^{r}(p)\(t\wedge {1 \over \psi (p)}\)^{r+1} \label{6.9}
   \end{equation} 
So we need only prove  (\ref{6.36})  and (\ref{6.7}) for $r=0$.    In this case, (\ref{6.36}) follows immediately from (\ref{8.15d}).

For (\ref{6.7}) we have 
  \begin{eqnarray}
&  &   \int_{0}^{1}   \(t \wedge {1 \over \psi(p)}\)  \,dp\nn \\
& &\qquad \leq  t \int_{0}^{\psi^{-1} (1/t)}  \,dp +  \int_{\psi^{-1} (1/t)}^{1} {1 \over  \psi (p)} \,dp     \le  Ct\psi^{-1} (1/t)\nn,
 \end{eqnarray}
 for all $t$ sufficiently large.

 By (\ref{1415}) and (\ref{612}), there exists a $t_{0}$ such that for all $t\geq t_{0}$,
\begin{eqnarray}
\lefteqn{\int_{0}^{t}\int_{0}^{\ff}  |  \sin  p\ga |\,  e^{-s\psi (p) }\,dp\,ds
\label{1416}}\\
&& \leq \int_{0}^{t_{0}}\int_{0}^{\ff}     e^{-s\psi (p) } \,dp\,ds+   \int_{t_{0}}^{t}\int_{0}^{\ff} |  \sin  p\ga |\,  e^{-s\psi (p) }\,dp\,ds  \nonumber\\
&& \leq C+C\ga \int_{t_{0}}^{t}\(\psi^{-1}(1/s)\)^{2}\,ds\leq C+C\ga \log t,  \nonumber
\end{eqnarray}
where for the last bound we   use (\ref{4.77}).
(This bound can not be smaller since  we may have  $\psi(p)=p^{2}$.)  
  \qed

 \medskip	
 
\noindent  {\bf  \Proof   of Lemma \ref{lem-vproprvt} }   For   the first bound in  (\ref{jr.1})
 we use (\ref{69dd}) and (\ref{jr.0}) with $r=0$ to get 
\begin{eqnarray}
&&\int_{0}^{t}p_{s }(x)\,ds \leq {1 \over \pi}\int_{0}^{\ff}  \int_{0}^{t} e^{-s\psi (p)} \,ds\,dp= O(t\psi^{-1} (1/t)),\label{h.10}
\end{eqnarray}
as $t\to\ff$. For the second bound in  (\ref{jr.1}) we use (\ref{88.m}), (\ref{el.1})  and (\ref{6.36}), 
to see that 
 \bea 
   \int_{0}^{t}p_{s }(x)\,ds& = &{1 \over     \pi}  \int_{0}^{\ff}\cos px\,  \int_{0}^{t}e^{-s\psi (p) }\,ds\,dp \label{h.101}\\
  & =& \Bigg | {1 \over     \pi x} \int_{0}^{\ff} \int_{0}^{t}e^{-s\psi (p) }\,ds \,d(\sin px)\Bigg |  \nn \\
  &\le &{1 \over     \pi |x|}  \int_{0}^{\ff}|\sin px|\,\left |\frac{d}{dp} \int_{0}^{t}e^{-s\psi (p) }\,ds\right |\,dp \nn\\
    &\le&\frac{C}{|x|}  \int_{0}^{\ff}|\sin px|  |\psi'(p)|\(t\wedge\frac{1}{\psi(p)}\)^{2} \,dp  \nn\\
        &\le&\frac{C}{|x|} \( \int_{0}^{1}{|\sin px|\over p}  \psi(p)\(t\wedge\frac{1}{\psi(p)}\)^{2} \,dp  \right.\nn\\
     &&\qquad\left.+\nn\int_{1}^{\ff}\frac{ |\psi'(p)|}{ \psi^{2}(p)}\,dp\)\\
     &\le &  C   \frac{t(1+\log^{+} |x|)}{|x|} \nn.
 \eea
Thus we get (\ref{jr.1}).

We next obtain (\ref{jrst1.3y}). Consider (\ref{72dd}).
For    $\ga> 1$
 \bea \lefteqn{ \int_{0}^{\ff}  {\sin^{2}(p\ga/2)\over    \psi( p) }\,dp\label{72j.2}}\\
 & &\le C\ga^{2}  \int_{0}^{1/\ga}  {p^{2}\over    \psi( p) }\,dp +  \int_{1/\ga}^{1}  {1\over    \psi( p) }\,dp+ \int_ {1} ^{\ff} {1\over    \psi( p) }\,dp\nn\\
  & &\le C\(\frac{1}{\ga\psi(1/\ga)}+1\)  \nn,
 \eea
   and for $\ga=1$ the integral is a constant. It follows from this   and (\ref{4.77})
    that 
    \be\sup_{x\in R^{1}}\int_{0}^{\ff}|\De^{ \ga}\De^{ -\ga} p_{s }(x)|
 \,ds\leq C \ga^{2}  .\label{72j.3}
 \end{equation}
    This gives the first bound in  (\ref{jrst1.3y}).   
     
 To obtain the third bound in (\ref{jrst1.3y}), consider (\ref{89.1dd})--(\ref{91dd}).      By  (\ref{jr.0}) with $r=0$,   we have 
      \bea 
  \int_{0}^{t} | I|\,ds &= &\ga^{2} \int_{0}^{t} \Big |\int_{0}^{\ff}\cos p\ga \sin^{2}( px/2)e^{-s\psi (p) }\,dp\Big |\,ds \label{88.12aa}\\
    &\le&\ga^{2}   \int_{0}^{\ff} \(\int_{0}^{t}e^{-s\psi (p) }\,ds\)\,dp \le C\ga^{2}t\psi^{-1} (1/t),\nn
     \eea
 for all $t$ sufficiently large.        
     Using (\ref{88.m}). (\ref{1.12})  and  (\ref{jr.0}) with $r=1$ we get 
    \bea
      \int_{0}^{t} | II|\,ds &=&2\ga\int_{0}^{t} \Big |\int_{0}^{\ff}\sin p\ga \sin^{2}( px/2)g'(p)\,dp\Big |\,ds \label{88.13} \\
    &\le& 2\ga \int_{0}^{1}|\sin (p\ga) \, \psi '(p)|\,\(\int_{0}^{t} s e^{-s\psi (p) }\,ds \)\,dp \nn \\
     & & \qquad +2\ga\int_{1}^{\ff}|\sin (p\ga) \, \psi '(p)|\,\(\int_{0}^{t} s e^{-s\psi (p) }\,ds \)\,dp\nn \\
       &\le& C\ga^{2}  \int_{0}^{ 1}  \, \psi (p)\,\(\int_{0}^{t} s e^{-s\psi (p) }\,ds \)\,dp\nn\\
    & & \qquad + C\ga\int_{1}^{\ff}{|  \psi '(p)|\over  \psi^{2}(p)}\,dp\nn \\
       &\le&C\ga^{2}\(t\psi^{-1} (1/t)\)\nn+C\ga\int_{1}^{\ff}{|  \psi '(p)|\over  \psi^{2}(p)}  \,dp\nn\\
    &\le&C\ga^{2}\(t\psi^{-1} (1/t)\)+C\ga\nn.   
          \eea
Similarly,  
        \bea 
    \int_{0}^{t} | III|\,ds &=& 2\int_{0}^{t} \Big |\int_{0}^{\ff}\sin^{2} (p\ga/2) \sin^{2}( px/2)g''(p)\,dp \Big |\,ds 
    \label{fgf}\\
     &\le &C\ga^{2} \int_{0}^{1}p^{2}\(\int_{0}^{t} \(  s |\psi''(p)|   + s^{2} |\psi'(p)|^{2}  \)e^{-s \psi (p) }\,ds \)  \,dp  
   \nn\\ 
   &  &\qquad+C \int_{1}^{\ff} \(\int_{0}^{\ff} \(  s |\psi''(p)|   + s^{2} |\psi'(p)|^{2}  \)e^{-s \psi (p) }\,ds \)  \,dp  
   \nn\\ 
   &\le &C\ga^{2} \int_{0}^{1} \(\int_{0}^{t} \(  s \psi(p)   + s^{2} \psi^{2}(p)   \)e^{-s \psi (p) }\,ds \)  \,dp  
    \nn\\
     & &\qquad+C \int_{1}^{\ff}    \( { \psi''(p)\over \psi^{2}(p)}  + { (\psi'(p))^{2}\over \psi^{3}(p)} \)  \,dp  
    \nn\\
&\le&C\ga^{2}t\psi^{-1} (1/t)+C\nn.
 \eea
Combining (\ref{88.12aa})-(\ref{fgf}) with (\ref{89.1dd})  we get the third bound in  (\ref{jrst1.3y}).

To get the second  bound in (\ref{jrst1.3y}) we use  the third integral in (\ref{89dd})  to see that 
\begin{equation}
\De^{ \ga}\De^{ -\ga} p_{s }(x)=-\frac{4}{\pi} {L(s,x) \over x  }
\label{89.1j}
\end{equation}
where
\begin{equation}
 L(s,x)=     \int_{0}^{\ff}\sin px\( \sin^{2}(p\ga/2) \, e^{-s\psi (p) }\)'\,dp.\label{89.2k}
\end{equation}   
  Using (\ref{90}) and (\ref{88.m}) we see that  
        \bea 
       \lefteqn{
 \int_{0}^{t}  |L | \,ds\label{90.mm}}\\&&\le C  \ga \int_{0}^{t}\int_{0}^{\ff } |  \sin p\ga|  g(p)\,dp\,ds+   C\int_{0}^{t}\int_{0}^{\ff }  \sin^2{(p\ga/2)}|  g'(p)|\,dp \,ds\nn.
   \eea
   By (\ref{6.7a})  the first term on the right-hand side is bounded by $C\ga^{2}\log t.$
  For the second term we note that  
   \begin{eqnarray}
   && \int_{0}^{t}\int_{0}^{\ff }  \sin^2{(p\ga/2)}|  g'(p)|\,dp \,ds
   \label{5.81}\\
  &&\qquad\le C\ga^{2}\int_{0}^{1}  p^{2}  |\psi'(p)|   
\int_{0}^{t}   s e^{-s\psi (p) }\,ds \,dp\nn\\
&&\hspace{1in} +C \int_{1}^{\ff}     |\psi'(p)|   
\int_{0}^{t}   s e^{-s\psi (p) }\,ds \,dp \nn\\
&&\qquad=IV+V.\nn
   \end{eqnarray}
By (\ref{el.1})
   \bea 
  IV&\le&C\ga^{2} \int_{0}^{1}  p \psi(p)
\(t^{2}\wedge \frac{1}{\psi^{2}(p)}\) \,dp\\
&\le& Ct^{2}\ga^{2}\int_{0}^{\psi^{-1}(1/t)}p \psi(p)\,dp+C\ga^{2}\int_{\psi^{-1}(1/t)}^{1}{p\over  \psi(p)}\,dp\nn\\
&\le& C  \ga^{2}t( \psi^{-1}(1/t))^{2} +C\ga^{2}\int_{\psi^{-1}(1/t)}^{1}{1\over p}\,dp\nn\\
&\le& C \ga^{2} \log t \nn ,
 \eea
where we use (\ref{4.77})  which   implies  that  $p/\psi(p)\le C/p$  for $p\in [0,1]$. The integral $V\le C$ by (\ref{el.1}) and (\ref{1.12}). Using all the material from (\ref{89.1j}) to this point we get the second bound in  (\ref{jrst1.3y}).  This completes the proof of (\ref{jrst1.3y}).

\medskip	 
 Using (\ref{6.10}), (\ref{72j.3}), (\ref{jrst1.3y}) and (\ref{6.7a}) we get the first bound in (\ref{2.4}).
 
  We now obtain the third bound in (\ref{2.4}). Considering (\ref {6.10}) and 
  (\ref{jrst1.3y}), it suffices to show that  
  \begin{equation}
  \int_{0}^{t} \Big |\int_{0}^{\ff}  \sin (px ) \sin (p\ga) \, e^{-s\psi (p) }\,dp \Big |\,ds\le C t \ga^{2} {1+  \log^{+}  |x|\over x^{2}}.  \label{70x.5}
\end{equation}
Consider (\ref{89s})--(\ref{550}). We have     \begin{eqnarray}
 \int_{0}^{t}|G_{1}|\,ds&=&\ga^{2}\int_{0}^{t} \Big |\int_{0}^{\ff}\sin px \sin p   \ga\, e^{-s\psi (p) } \,dp\Big |\,ds\label{89.2e}\\
  &  
  \leq &  \ga^{2}\int_{0}^{t}  \int_{0}^{\ff} \, e^{-s\psi (p) } \,dp\leq C\ga^{2}t\psi^{-1} (1/t) \nonumber,
  \end{eqnarray}   
by (\ref{jr.0}).

Using   (\ref{6.36})  and (\ref{1.12}), we see that 
    \begin{eqnarray}
 \int_{0}^{t}|G_{2}|\,ds&=&2 \ga\int_{0}^{t} \Big |\int_{0}^{\ff}\sin px \cos p \ga\( \psi'(p)\,s e^{-s\psi (p)}\) \,dp\Big |\,ds
  \label{89.2f}\\
  &\le&   2  \ga\int_{0}^{1}|\sin px||\psi'(p)|\(\,\int_{0}^{t}s e^{-s\psi (p)}\,ds\) \,dp + 2 \ga \int_{1}^{\ff} {|\psi'(p)|\over \psi^{2} (p)}   \,dp \nonumber\\
  &\le&   C \ga\int_{0}^{1}{|\sin px|\over p}\psi(p)\(t\wedge {1 \over \psi (p)}\)^{2} \,dp+C \ga \nonumber\\
  &\le&   C\ga t\( 1+  \log^{+} |x| \)\nonumber .
  \end{eqnarray} 
Similarly,  
   \begin{eqnarray}
\lefteqn{  \int_{0}^{t}|G_{3}|\,ds\label{89.2g}}\\
  &&=\int_{0}^{t} \Big |\int_{0}^{\ff}\sin px\,\, \sin p \ga \( s\psi''(p)   - s^{2} (\psi'(p))^{2}    \)e^{-s\psi (p) } \,dp\Big |\,ds
  \nn\\
    &&\leq C  \ga \int_{0}^{ 1}|\sin px|\,p\(\int_{0}^{t}\( s|\psi''(p)| + s^{2} (\psi'(p))^{2}   \)e^{-s\psi (p) } \,ds\) \,dp
  \nn\\
    &&\quad+ C  \int_{1}^{ \ff} \(\int_{0}^{\ff}\( s|\psi''(p)| + s^{2} (\psi'(p))^{2}   \)e^{-s\psi (p) } \,ds\) \,dp
  \nn\\
  &&\le C \ga  \int_{0}^{1}{|\sin px|\over p}   \psi(p)\(t\wedge {1 \over \psi (p)}\)^{2}\,dp\nn\\
  && \hspace{.2in}  + C \ga  \int_{0}^{1}{|\sin px|\over p}   \psi^{2}(p)\(t\wedge {1 \over \psi (p)}\)^{3}\,dp    
  \nn\\
  &&\qquad +C \int_{1}^{ \ff} {|\psi''(p)|\over \psi^{2}(p) }\,dp+ C \int_{1}^{ \ff} {(\psi'(p))^{2} \over \psi^{3}(p) }\,dp \nn\\
  &&  \leq C \ga \( 1+  \log^{+} |x| \)t.\nonumber
  \end{eqnarray}
  
This completes the proof of  (\ref{70x.5}) and gives us the third bound in (\ref{2.4})

The second bound in (\ref{2.4}) follows from the third line of (\ref{89s}) and the observation that  
\bea
 && \Bigg |\int_{0}^{t} \int_{0}^{\ff}\cos px\(\sin (p \ga) \, e^{-s\psi (p) })\)'\,dp \Bigg | \\
 &&\qquad\le\nn\ga\int_{0}^{\ff}\int_{0}^{t}e^{-s\psi (p) }\,ds\,dp+\int_{0}^{\ff}\int_{0}^{t}|\sin p \ga||\psi'(p)|se^{-s\psi(p)}\,ds\,dp\\
 &&\qquad\le C\ga(t\psi^{-1}(1/t))+\ga\int_{0}^{1} \psi(p)\(t\wedge\frac{1}{\psi(p)}\)^{2}\,dp\nn\\
 &&\hspace{1in}+\int_{1}^{\ff}\frac{|\psi'(p)|}{\psi^{2}(p)}\,dp\nn\\
  &&\qquad\le C(1\vee\ga)(t\psi^{-1}(1/t))\nn.
   \eea 
   In this chain of inequalities we use (\ref{jr.0}), (\ref{88.m}), (\ref{1.12}) and (\ref{6.7}). 
\qed

\noindent {\bf Proof  of Lemma \ref{lem-vproprvtj} }
The inequalities in (\ref{jr.1j})--(\ref{jrst1.3yj}) follow immediately from Lemma \ref{lem-vproprvt}. 

The inequality in (\ref{ee1}) is trivial, since $p_{s}(x)$ is a probability density for all $s>0$. 
 
To obtain    (\ref{151ee1}) we use (\ref{2.4}) with $\ga=1$  to  get \label{page65}
  \begin{eqnarray}
   &&\int _{0}^{\ff} v (x,t)\,dx
  \label{72j.3f} \\
  &  &\qquad\le C \( \log t\int_{0}^{t\psi^{-1}(1/t)}\,dx+t\psi^{-1}(1/t)\int_{t\psi^{-1}(1/t)}^{2c/\psi^{-1}(1/t)}\frac{1}{x}\,dx\right.\nn\\
  &&\qquad\qquad\qquad\left.+Ct\int_{{2c/\psi^{-1}(1/t)}}^{\ff}     {\log  x\over x^{2}}\,dx\) \nonumber\\
  &&\qquad\le  C(t\psi^{-1}(1/t) \log  t),\nn\end{eqnarray}
for all $t$ sufficiently large. Note that when $\bb<2$ it is clear that $t\psi^{-1}(1/t)<1/\psi^{-1}(1/t)$ for all $t$ sufficiently large. In general we use (\ref{4.77})  with $c$ representing the constant.

  For (\ref{132ee2}) we use  (\ref{jrst1.3y})    to see that
  \begin{eqnarray}
  &&\int_{0}^{\ff}  \(\,\int_{0}^{t}\,|\De^{ 1}\De^{ -1} \,p_{s }(x)|\,ds\)^{2}\,dx
  \label{72j.3n}\\
  &&\qquad \le C \int_{0}^{\log t}  \(\,\int_{0}^{t}\,|\De^{ 1}\De^{ -1} \,p_{s }(x)|\,ds\)^{2}\,dx\nn\\
  &&\qquad\qquad \qquad+2\int_{\log t}^{\ff}  \(\,\int_{0}^{t}\,|\De^{ 1}\De^{ -1} \,p_{s }(x)|\,ds\)^{2}\,dx  \nonumber\\
  &&\qquad \leq C  \log t +C\int_{\log t}^{\ff}  {(\log t)^{2}\over x^{2}} \,dx \le C\log t \nn.
  \end{eqnarray}
A similar argument gives
(\ref{132ee3}) since
  \be 
   \int_{u}^{\ff}  \(\,\int_{0}^{t}\,|\De^{ 1}\De^{ -1} \,p_{s }(x)|\,ds\)^{2}\,dx
 \leq   C\int_{u}^{\ff}  {(\log t)^{2}\over x^{2}} \,dx =C{(\log t)^{2}\over u}.  \label{72j.3nv}
  \ee
Finally, to obtain (\ref{132ee1}) we use  (\ref{jrst1.3y})  to get
\begin{eqnarray}
&&\int_{0}^{\ff} \int_{0}^{t} \Big |\De^{ 1}\De^{ -1}\,p_{s }(x)\Big |\,ds\,dx
\label{st.30}\\
&&= \int_{0}^{1}\int_{0}^{t} \Big |\De^{ 1}\De^{ -1}\,p_{s }(x)\Big |\,ds\,dx
+\int_{1}^{t\psi^{-1} (1/t)}  \int_{0}^{t} \Big |\De^{ 1}\De^{ -1}\,p_{s }(x)\Big |\,ds\,dx\nn\\
&&\hspace{1 in}
+\int_{t\psi^{-1} (1/t)}^{\ff} \int_{0}^{t} \Big |\De^{ 1}\De^{ -1}\,p_{s }(x)\Big |\,ds\,dx  \nonumber\\
&&\leq C \int_{0}^{1}1\,dx+C\log t\int_{1}^{t\psi^{-1} (1/t)} {1 \over |x|} \,dx\nn\\
&&\hspace{1 in}
+C\int_{t\psi^{-1} (1/t)}^{\ff}{t\psi^{-1} (1/t) \over |x|^{2}}\,dx  \nonumber\\
&&\leq C   +C  (\log  t)^{2}+C \nn. 
\end{eqnarray}
 \qed

 \begin{lemma}\label{lem-new3}   Under the hypotheses of 
Theorem  \ref{theo-clt2r},  for all $t$ sufficiently large   and all  $x\in R^{1}$   
    \bea
  \ov u(x,t) :=\sup_{\de t\le s\le t} p_{s}(x)&\le &C\( \psi^{-1}(1/t) \wedge\frac{1}{\psi^{-1}(1/t)x^{2}} \) \label{8.15aq};\qquad\\
 \ov v(x,t):=    \sup_{\de t\le s\le t}  |\De^{1}p_{s}(x)|&\le& C\( \(\psi^{-1}(1/  t)\)^{2}\wedge  \frac{1+\log^{+} |x|}{ x^{2}}\);\qquad\label{151q}\\
 \ov w(x,t):=\sup_{\de t\le s\le t}  |\De^{1} \De^{-1}p_{s}(x)|&\le& C \( \(\psi^{-1}(1/  t)\)^{3}\wedge \frac{\psi^{-1}(1/  t)}{ x^{2}}\)\label{132q}. \,
\eea
In addition  
 \bea
   \int   \ov u(x,t)\,dx&\le &C  \label{8.15aqw};\\
  \int  \ov v(x,t)\,dx&\le& C   \psi^{-1}(1/  t)\log t ;\label{151qw}\\
 \int  \ov w(x,t)\,dx&\le& C   \(\psi^{-1}(1/  t)\)^{2}\leq \frac{C}{t}. \label{132qw}
\eea
 \end{lemma}

\Proof By (\ref{8.15a})
 \begin{equation}
     \sup_{\de t\le s\le t} p_{s}(x)\le C\( \psi^{-1}(1/\de t) \wedge\frac{1}{\psi^{-1}(1/t)x^{2}} \) \label{8.15aqqq},
   \end{equation}
and by the regular variation property $\psi^{-1}(1/\de t)\le C   \psi^{-1}(1/  t)$. (The constant depends on $\de$ but that doesn't matter.)   The inequalities in (\ref{151q}) and (\ref{132q}) follow similarly from (\ref{151}) and (\ref{132}). 
 
 The inequalities in (\ref{8.15aqw})--(\ref{132qw}) follow easily from (\ref{8.15aq})--(\ref{132q}). For (\ref{8.15aqw}) we write 
 \be
\int\ov u(x,t)\,dx\le C\int _{0}^{a} \psi^{-1}(1/t)\,dx+ \int _{a}^{\ff}\frac{1}{\psi^{-1}(1/t)x^{2}}  \,dx,
 \ee
 where $a=1/\psi^{-1}(1/t)$. For (\ref{151qw}) and (\ref{132qw}) we proceed similarly with $a=1/\psi^{-1}(1/t)$ in both cases. \qed
 
  {\bf \Proof  of Lemma \ref{lem-2.4} }   The inequality in (\ref{mac.1}) follows from (\ref{6.10}) and (\ref{132}).
  
  To obtain (\ref{wo.4a}) we write
  \bea
 &&  \int_{0}^{2t}\int_{{0}}^{2t} |\De^{1}p_{r+s }(0)|\,dr\,ds\\
  &&\qquad = \int_{{0}}^{2t}u |\De^{1}p_{u }(0)|\,du+  \int_ {2t}^{4t}(4t-u) |\De^{1}p_{u }(0)|\,du\nn.
        \eea
  By (\ref{6.10}) and (\ref{72dd})   
\bea
  \int_{{0}}^{2t}u |\De^{1}p_{u }(0)|\,du &\le &\frac{2}{\pi}\int_{0}^{2t}u \int _{0}^{\ff} \sin^{2}(p   /2) \,e^{-u \psi (p) }\,dp\,du \\
   &=&\frac{2}{\pi}   \int _{0}^{\ff} \sin^{2}(p   /2) \int_{0}^{2t} u \,e^{-u\psi (p) }\,du\,dp. \nn
         \eea
  In addition  
  \begin{equation}
   \int_{0}^{2t} u \,e^{-u\psi (p) }\,du\le\( \frac{1}{\psi^{2}(p)}\(1-e^{-2t\psi (p)}\)\wedge \frac{Ct}{\psi(p)}\),
   \end{equation}       
  where, for the final inequality we use Lemma \ref{lem-6.3}. Consequently, for all $t$ sufficiently large,
  \bea 
   && \int _{0}^{\ff} \sin^{2}(p   /2) \int_{0}^{2t} u \,e^{-u\psi (p) }\,du\,dp\label{6.97} \\
   &&\qquad\le Ct \int_{0}^{\psi^{-1}(1/t)}\frac{p^{2}}{\psi (p)}\,dp+ \frac{1}{4}\int _ {\psi^{-1}(1/t)}^{1}\frac{p^{2}}{\psi^{2} (p)}+\nn\int_{1}^{\ff}\frac{1}{\psi^{2}(p)}\,dp\\
  &&\qquad\le C\(t^{2}\(\psi^{-1}(1/t)\)^{3}+ 1\)+\frac{1}{4}\int _ {\psi^{-1}(1/t)}^{1}\frac{p^{2}}{\psi^{2} (p)}\,dp.\nn 
 \eea 
  Note that 
   \begin{eqnarray} &&
\int _ {\psi^{-1}(1/t)}^{1}\frac{p^{2}}{\psi^{2} (p)}\,dp\nn  \le\left\{\begin{array}{ll} \displaystyle
 Ct^{2}\(\psi^{-1}(1/t)\)^{3}  &\mbox{ if }\bb> 3/2\\
L(t)&\mbox{ if }\bb=3/2\\
C&\mbox{ if }\bb<3/2,
\end{array}
\right.
\label{e1.8aj}
\end{eqnarray}
 where $L(t)$ is a slowly varying function at infinity. Therefore
 \begin{equation}
  \int_{{0}}^{2t}u |\De^{1}p_{u }(0)|\,du\le   C\(t^{2}\(\psi^{-1}(1/t)\)^{3}+L(t)+ 1\).
   \end{equation}        
 In addition    
 \bea 
    \int_ {2t}^{4t}(4t-u) |\De^{1}p_{u }(0)|\,du&\le& Ct \int_ {2t}^{4t}  |\De^{1}p_{u }(0)|\,du\\
    &=&Ct\int_ {0}^{2t}  |\De^{1}p_{v+2t}(0)|\,dv\nn
    \eea
  Note that by (\ref{6.10}), (\ref{72dd}) and (\ref{132q})
 \bea
  |\De^{1}p_{v+2t}(0)|&=&\frac{1}{2}  |\De^{1}\De^{-1}p_{v+2t}(0)|\\
  &\le&\frac{1}{2} |\De^{1}\De^{-1}p_{ 2t}(0)|\le\nn C\(\psi^{-1}(1/t)\)^{3}.
 \eea
    Here we also use the fact that $\De^{1}\De^{-1}p_{s}(0)$ is decreasing in $s$, and the regular variation of $\psi$. Consequently
 \begin{equation}
       \int_ {2t}^{4t}(4t-u) |\De^{1}p_{u }(0)|\,du \le Ct^{2}\(\psi^{-1}(1/t)\)^{3}.
   \end{equation}
   Thus we obtain (\ref{wo.4a}). \qed

   {\bf \Proof  of Lemma \ref{lem-2.5} }   The equality in  (\ref{2.13}) follows easily from (\ref{69dd}).

   The equality in (\ref{h3.1}) follows from (\ref{f9.2}) integrated with respect to $r$ and $r'$.
   
   For (\ref{h3.2}) we use Parseval's Theorem, (see (\ref{f9.2}))   to get
\be 
\int  \(\int_{0}^{t }\De^{ 1}\De^{ -1}\,\,p_{s}(x)\,ds\)^{2}\,dx\label{rst.28} ={16 \over  \pi}\int_{0}^{\ff} { \sin ^{4}(p/2) \over \psi^{2}(p)  }\(1-e^{-t \psi (p) }\)^{2}\,dp\nn.
\ee
 To complete the proof of (\ref{h3.2}) we   note that by (\ref{bound})  
\be 
\int_{0}^{\ff} {\sin ^{4}(p/2) \over \psi^{2}(p)  } e^{-t \psi (p) } \,dp\le \frac C{t^{1/3}} \int_{0}^{\ff} {\sin ^{4}(p/2) \over \psi^{7/3}(p)  }   \,dp \le \frac C{t^{1/3}}.
\label{rst.28m} 
\ee
  \qed

\noindent{\bf Proof of Lemma  \ref{lem-3.2} } The first inequality is given in (\ref{2.4}). The second inequality follows from  the definition of $v_{*}$ in (\ref{vstar}). For (\ref{v_{***}}) we note that in the proof of (\ref{151ee1}),   on page \pageref{page65}, we  are actually integrating $v_{*}(x,t)$. \qed

\section{Proof of Lemma \ref{lem-ilt} }\label{sec-ilt}

Set 
\begin{equation}
h=\psi^{-1}(1/t) \label{h},\hspace{.2 in}\mbox{so that }\hspace{.2 in}\psi (h)=1/t. \label{h.0}
 \end{equation}

\noindent {\bf  Proof of Lemma \ref{lem-ilt} } By Kac's moment formula, (see (\ref{e1.60})),
 \begin{eqnarray} 
\lefteqn{ E\(\(\int  (L^{ x}_{t})^{2} \,dx\)^{n}\) \label{ilt.2}}\\
&& =\sum_{ \si  }\int  \(\int_{\{\sum_{i=1}^{2n} s_{i}\leq t\}}\prod_{ i=1}^{2 n}\, p_{s_{i} }(y_{\si (i)}-y_{\si (i-1)})\,ds_{i}\) \prod_{i=1}^{n}\,dy_{i} \nn\\
&& =\sum_{ \si  }\frac{1}{(2\pi)^{2n}}\int  \(\int_{\{\sum_{i=1}^{2n} s_{i}\leq t\}}\prod_{ i=1}^{2 n}\, \int e^{ip_{i}(y_{\si (i)}-y_{\si (i-1)})}e^{-s_{i}\psi  (p_{i}) }\,dp_{i}\,ds_{i}\) \nn\\
&&\hspace{4.2in}\prod_{i=1}^{n}\,dy_{i} \nn \\
&& =\({t\over 2\pi}\)^{2n} \sum_{ \si  }\int  \(\int_{\{\sum_{i=1}^{2n} s_{i}\leq 1\}}\prod_{ i=1}^{ 2n}\, \int e^{ip_{i}(y_{\si (i)}-y_{\si (i-1)})}e^{-s_{i}t\psi  (p_{i}) }\,dp_{i}\,ds_{i}\)\nn\\
&&\hspace{4.2in}\prod_{i=1}^{n}\,dy_{i} \nn.
 \end{eqnarray}

  Here the  sum  in the second line runs over all  maps $\si $ of $\{1,\ldots, 2n\}$ into $\{1,\ldots, n\}$, such that $|\si ^{-1}(j)|=2$ for each $1\leq j\leq n$, and we set $\si(0)=0$.  Thus, by   (\ref{h.0})  and many changes of variables  
 \begin{eqnarray}
\lefteqn{ (2\pi)^{2n}(t^{2}\psi^{-1} (1/t))^{-n} E\(\(\int  (L^{ x}_{t})^{2} \,dx\)^{n}\) \label{ilt.4}}  \\
&& =(2\pi)^{2n}t^{-2n}h^{-n} E\(\(\int  (L^{ x}_{t})^{2} \,dx\)^{n}\)
     \nn\\
&& =h^{-n}\sum_{ \si  }\int  \(\int_{\{\sum_{i=1}^{2n} s_{i}\leq 1\}}\prod_{ i=1}^{2 n}\, \int e^{ip_{i}(y_{\si (i)}-y_{\si (i-1)})}e^{-s_{i}t\psi  (p_{i}) }\,dp_{i}\,ds_{i}\) \prod_{i=1}^{n}\,dy_{i} \nn\\
&& =h^{n}\sum_{ \si  }\int  \(\int_{\{\sum_{i=1}^{2n} s_{i}\leq 1\}}\prod_{ i=1}^{ 2n}\, \int e^{ip_{i}h(y_{\si (i)}-y_{\si (i-1)})}e^{-s_{i}t\psi  (p_{i}h) }\,dp_{i}\,ds_{i}\) \prod_{i=1}^{n}\,dy_{i} \nn\\
&& =\sum_{ \si  }\int  \(\int_{\{\sum_{i=1}^{2n} s_{i}\leq 1\}}\prod_{ i=1}^{2 n}\, \int e^{ip_{i}(y_{\si (i)}-y_{\si (i-1)})}e^{-s_{i}\psi  (p_{i}h)/\psi (h)}\,dp_{i}\,ds_{i}\) \prod_{i=1}^{n}\,dy_{i}. \nn
 \end{eqnarray}
Using the regular variation of $\psi$ at zero 
 the proof follows once we   justify interchanging the limit and the integrals.

 For $\si $ fixed let   
 \bea
f_{h}(y)&=& \int_{\{\sum_{i=1}^{2n} s_{i}\leq 1\}}\prod_{ i=1}^{2 n}\, \int e^{ip_{i}y_{i}}e^{-s_{i}\psi  (p_{i}h)/ \psi (h)}\,dp_{i}\,ds_{i}\label{ilt.6j}\\
 &=& 2^{2n}\int_{\{\sum_{i=1}^{2n} s_{i}\leq 1\}}\prod_{ i=1}^{2 n}\, \int_{0}^{\ff} \cos p_{i}y_{i}\,e^{-s_{i}\psi  (p_{i}h)/ \psi (h)}\,dp_{i}\,ds_{i}\nn
 \eea
Considering  (\ref{ilt.4})  it suffices to show that  for each fixed $y=(y_{1},\ldots, y_{n})$   
 \begin{equation}
 \lim_{h\rar 0} f_{h}(y)=2^{2n} \int_{\{\sum_{i=1}^{2n} s_{i}\leq 1\}}\prod_{ i=1}^{2 n}\, \int_{0}^{\ff} \cos p_{i}y_{i}\,e^{-s_{i} p_{i}^{\bb}}\,dp_{i}\,ds_{i}.\label{ilt.7}
 \end{equation}
and   $f_{h}(y)$ is bounded and integrable in $y$, uniformly in $h\leq h_{0}$, for some $ h_{0}>0$, sufficiently small. In fact we show that     
\begin{equation}
   \sup_{h\le h_{0}}|f_{h}(y)|\le C\prod_{i=1}^{2n}\(1\wedge\frac{1}{y_{i}^{2}}\).\label{8.6}
   \end{equation}

We first obtain (\ref{ilt.7}). For $M$ large, write 
\begin{equation}
1=\prod_{i=1}^{2n}\(1_{\{0\leq p_{i}\leq M\}}
+1_{\{p_{i}\geq M\}}\)\label{ilt.8}
\end{equation}
and  
  \begin{equation}
f_{h}(y)=2^{2n} \int_{\{\sum_{i=1}^{2n} s_{i}\leq 1\}}\prod_{ i=1}^{ 2n}\, \int_{0}^{M} \cos p_{i}y_{i}\,e^{-s_{i}\psi  (p_{i}h)/ \psi (h)}\,dp_{i}\,ds_{i}+G_{h}.\label{ilt.6a}
 \end{equation}
Here  $G_{h}$ is a sum of many terms, in each of which  $p_{i}\geq M$,  for at least one $1\le i\le 2n$. Suppose  there are $k$ terms with $p_{i}\geq M$.  We   bound these terms by
 \be 
2^{2n} \(\int_{0}^{1}\int_{0}^{M}  e^{-s \psi  (p h)/ \psi (h)}\,dp \,ds \)^{2n-k}
 \(\int_{0}^{1}\int_{M}^{\ff}  e^{-s \psi  (p h)/ \psi (h)}\,dp \,ds \)^{k}.
 \label{ilt.9}
 \ee
By (\ref{regcond}),  for any $\ep>0$, (see also \cite[Theorem 1.5.6]{BGT}),
 \begin{eqnarray}
 &&\int_{0}^{1}\int_{0}^{M}  e^{-s \psi  (p h)/ \psi (h)}\,dp \,ds
 \label{pba}\\
  &&\qquad\leq 1+\int_{0}^{1}\int_{1}^{M}  e^{-s \psi  (p h)/ \psi (h)}\,dp \,ds   \nonumber\\
  &&\qquad\leq 1+\int_{0}^{1}\int_{1}^{M}  e^{-s Cp^{\bb-\ep}}\,dp \,ds,   \nonumber
 \end{eqnarray}
 which is bounded by a constant independent of $M$. 
Using the regular variation of $\psi$  at zero, we have
 \begin{eqnarray}
\lefteqn{\int_{0}^{1}\int_{M}^{\ff}   e^{-s \psi  (p h)/ \psi (h)}\,dp \,ds\label{ilt.10}}\\
 && \le \psi (h)\int_{M}^{\ff} {1\over \psi (hp) } \,dp\nn =\frac{\psi(h)}{h} \int_{hM}^{\ff}\frac{1}{\psi(s)}\,ds\nn\\
 &&=\frac{\psi(h)}{h} \int_{hM}^{1}\frac{1}{\psi(s)}\,ds+\frac{\psi(h)}{h} \int_{1}^{\ff}\frac{1}{\psi(s)}\,ds\nn\\
 &&\le C \frac{\psi(h)M}{\psi(hM)}+C\frac{\psi(h)}{h}\nn,
   \eea    
   for all $h$ sufficiently small. 
   Therefore,   as in (\ref{pba}),
   \begin{equation}
   \limsup_{h\to 0}\int_{0}^{1}\int_{M}^{\ff}   e^{-s \psi  (p h)/ \psi (h)}\,dp \,ds\le \frac{ C}{M^{\bb-1}}.
   \end{equation}
Thus 
      \begin{equation}
   \limsup_{h\to 0} |G_{h}|\le \frac{ C}{M^{\bb-1}}.\label{8.11}
   \end{equation}
   
   Now consider the integral in (\ref{ilt.6a}). By the regular variation of $\psi$  at zero
  and the Dominated Convergence Theorem, 
 \begin{eqnarray}
&& \lim_{h\rar 0} \int_{\{\sum_{i=1}^{2n} s_{i}\leq 1\}}\prod_{ i=1}^{2 n}\, \int_{0}^{M} \cos p_{i}y_{i}\,e^{-s_{i}\psi  (p_{i}h)/ \psi (h)}\,dp_{i}\,ds_{i}
 \label{ilt.11} \\
 &&\qquad = \int_{\{\sum_{i=1}^{2n} s_{i}\leq 1\}}\prod_{ i=1}^{2 n}\, \int_{0}^{M} e^{ip_{i}y_{i}}\cos p_{i}y_{i} \,e^{-s_{i} p_{i}^{\bb}}    \,dp_{i}\,ds_{i}.  \nonumber 
 \end{eqnarray} 
Thus we get (\ref{ilt.7}). 

 We show below that for any $J\subseteq \{1,\ldots,2n \}$ we have
\begin{equation}
   \sup_{h\le h_{0}}|f_{h}(y)|\le C\prod_{i\in J} \(\frac{1}{y_{i}^{2}}\).\label{8.6j}
   \end{equation}
In particular, (\ref{8.6j}) also holds when $J $ is the empty set, so that $   \sup_{h\le h_{0}}\newline |f_{h}(y)|\le C$. Using this it is easy to see that   (\ref{8.6}) holds.

It  follows from integrating by parts twice that 
   \be 
    \int_{0}^{\ff} \cos (py)\,\,  e^{-s \psi  (p h)/ \psi (h)}\,dp\label{112a} =-\frac{1}{y^{2}}  \int_{0}^{\ff} \cos (py)\(e^{-s \psi  (p h)/ \psi (h)}\)^{''}\,dp.\nn
    \ee 
 where we use      the  fact  that $\psi'(0)=0$,  which follows from  (\ref{regcond})  and the first inequality in (\ref{88.m}).
    Applying this for all  $i\in J$ we see that  
     \bea
f_{h}(y)&=&   \prod_{i\in J} \(\frac{-1}{y_{i}^{2}}\) \int_{\{\sum_{i=1}^{2n} s_{i}\leq 1\}} \prod_{i\in J}\, \int_{0}^{\ff} \cos p_{i}y_{i}\,\(e^{-s_{i}\psi  (p_{i}h)/ \psi (h)}\)^{''}\,dp_{i}\,ds_{i}\nonumber\\
&&\hspace{1 in}  \prod_{i\in J^{c}}\, \int_{0}^{\ff} \cos p_{i}y_{i}\,e^{-s_{i}\psi  (p_{i}h)/ \psi (h)}\,dp_{i}\,ds_{i}.\label{ilt.6}
 \eea
 Therefore
     \bea
|f_{h}(y)|&\leq &  
 \prod_{i\in J} \(\frac{1}{y_{i}^{2}}\) \(\,\int_{0}^{1} \int_{0}^{\ff}  \,\bigg|\(e^{-s \psi  (p h)/ \psi (h)}\)^{''}\bigg|\,dp \,ds \)^{|J|}\nonumber\\
&&\hspace{.7 in}  \( \int_{0}^{1}\int_{0}^{\ff}  \,e^{-s \psi  (p h)/ \psi (h)}\,dp \,ds\)^{|J^{c}|}.\label{ilt.6k}
 \eea

\medskip 	It is easily seen that
\begin{eqnarray}
 \int_{0}^{1}\int_{0}^{\ff}  \,e^{-s \psi  (p h)/ \psi (h)}\,dp \,ds&\leq &C\sup_{h}\int \(1\wedge {\psi (h) \over  \psi(hp)}\)\,dp  .\label{rvc.8}
\end{eqnarray}
  Therefore,  for $h\le 1$,  
\begin{eqnarray}
&&\int \(1\wedge {\psi (h) \over  \psi(hp)}\)\,dp
\label{rvc.8h}\\
&&\qquad\leq  \int_{0}^{1} 1 \,dp+\int_{1}^{1/h}   {\psi (h) \over  \psi(hp)}\,dp +\int_{1/h}^{\ff}   {\psi (h) \over  \psi(hp)}\,dp  \nonumber\\
&&\qquad= 1+{\psi (h) \over h}\(\int_{h}^{1}   {1 \over  \psi(p)}\,dp+\int_ {1}^{\ff}   {1 \over  \psi(p)}\,dp\)\leq C.  \nonumber
\end{eqnarray}
 Consequently, to obtain
(\ref{8.6j})   we need only show    that,   for $h\le 1$,
   \begin{equation}
  \int_{0}^{1} \int_{0}^{\ff}  \bigg|\(e^{-s \psi  (p h)/ \psi (h)}\)^{''} \bigg|\,dp\,ds<\ff.\label{8.21}
   \end{equation}
 We have  
  \bea
 &&   \bigg|\(e^{-s \psi  (p h)/ \psi (h)}\)^{''} \bigg|\\
 &&\qquad\le      \bigg| { h^{2}(\psi'(hp))^{2}\over\psi^{2}(h)}  \bigg|s^{2}e^{-s \psi  (p h)/ \psi (h)}+\bigg| {h^{2} \psi''(hp)) \over\psi (h)}  \bigg|se^{-s \psi  (p h)/ \psi (h)}\nn.
   \eea
 Using  (\ref{88.m}) and  (\ref{1.12})
  we see that,   for $h\le 1$,
\begin{eqnarray}
\lefteqn{  \int_{0}^{1} \int_{0}^{\ff} \bigg| {h^{2}(\psi'(hp))^{2}\over\psi^{2}(h)}  \bigg|s^{2}e^{-s \psi  (p h)/ \psi (h)}\,dp\,ds
\label{rvc.8i}}\\
&&\leq C   \int_{0}^{\ff} \bigg| {h^{2}(\psi'(hp))^{2}\over\psi^{2}(h)}  \bigg|\(1\wedge {\psi (h) \over  \psi(hp)}\)^{3}\,dp\nonumber\\
&&\leq \frac{Ch^{2}}{\psi^{2}(h)} \int_{0}^{1}(\psi'(hp))^{2}\,dp+  Ch^{2}  \psi(h) \int_{1}^{\ff}\frac{(\psi'(hp))^{2}}{\psi^{3}(hp)}\,dp\nn\\
&&\leq \frac{C h }{\psi^{2}(h)} \int_{0}^{h}(\psi'(s))^{2}\,ds+ C  h  \psi (h) \int_{h}^{\ff}\frac{(\psi'(s))^{2}}{\psi^{3}(s)}\,ds\nn\\
&&\leq C\( \frac{h }{\psi^{2}(h)} \int_{0}^{h}{ \psi^{2}(s)\over s^{2}}\,ds+  h  \psi (h) \int_{h}^{1}\frac{1}{s^{2}\psi (s)}\,ds+C\)\nn\\
     &&\leq C'\nn.
\end{eqnarray}
Similarly, 
\begin{eqnarray}
\lefteqn{  \int_{0}^{1} \int_{0}^{\ff} \bigg| {h^{2} \psi''(hp)) \over\psi (h)}  \bigg|se^{-s \psi  (p h)/ \psi (h)}\,dp\,ds
\label{rvc.8ia}}\\
&&\leq C   \int_{0}^{\ff}\bigg| {h^{2} \psi''(hp)) \over\psi (h)}  \bigg|\(1\wedge {\psi (h) \over  \psi(hp)}\)^{2}\,dp\nonumber\\
&&\leq  \frac{ C h^{2}}{\psi (h)} \int_{0}^{1} \psi''(hp)  \,dp+   Ch^{2} \psi (h) \int_{1}^{\ff}\frac{\psi''(hp)}{\psi^{2}(hp)}\,dp\nn\\
&&\leq  \frac{ C h }{\psi (h)} \int_{0}^{h} \psi''(s) \,ds+  C h  \psi (h) \int_{h}^{\ff}\frac{\psi''(s)}{\psi^{2}(s)}\,ds\nn\\
&&\leq C\( \frac{ h }{\psi (h)} \int_{0}^{h}{ \psi (s)\over s^{2}}\,ds+  h  \psi (h) \int_{h}^{1}\frac{1}{s^{2}\psi (s)}\,ds+C\)\nn\\
&&\leq  C'.\nn
   \eea
   Thus we obtain (\ref{8.21}).
\qed

\section{Estimates for the mean  and variance   }\label{sec-9}

\noindent{\bf Proof of Lemma  \ref{lem-varep} }
By the     Kac  moment formula
\begin{eqnarray}
&&E\(\int ( L^{ x+1}_{t}- L^{ x}_{ t})^{ 2}\,dx\)
\label{kacv.3}\\
&&\qquad =2 \int \int_{\{\sum_{i=1}^{2}r_{i}\leq t\}}\De^1p_{r_{1}}(x)\De^1 p_{r_{2}}( 0)\,dr_{1}\,dr_{2}\,dx \nonumber\\
&&\quad \qquad+2 \int \int_{\{\sum_{i=1}^{2}r_{i}\leq t\}}p_{r_{1}}(x)\De^1 \De^{-1} p_{r_{2}}( 0)\,dr_{1}\,dr_{2}\,dx. \nonumber
\end{eqnarray}
When we  integrate with respect to  $x$ we get zero in the first integral and one in the second. Consequently, by (\ref{69dd})
\begin{eqnarray}
E\(\int ( L^{ x+1}_{t}- L^{ x}_{ t})^{ 2}\,dx\) &=& 2   \int_{\{\sum_{i=1}^{2}r_{i}\leq t\}} \De^1 \De^{-1} p_{r_{2}}( 0)\,dr_{1}\,dr_{2}\label{kacv.4}\\
&=& 4  \int_{0}^{t} (t-r)\(p_{r }( 0)-p_{r }( 1)\) \,dr \nn\\
 &=&  \frac{8}{\pi}\int_{0}^{\ff} \sin^{2} p/2 \int_{0}^{t} (t-r)e^{-r\psi(p)} \,dr\,dp.\nn
\end{eqnarray}
Note that  
\begin{equation}
   \int_{0}^{t} (t-r)e^{-r\psi(p)}\,dr=\frac{t}{\psi(p)}-\frac{1-e^{-t\psi(p)}}{\psi^{2}(p)}.
   \end{equation}
 By (\ref{1.16a})  
\begin{equation}
   {8 t\over \pi}\int_{0}^{\ff}   {\sin^{2} (p/2) \over \psi(p)}  \,dp  =4c_{\psi,0}t.\label{7.6}
   \end{equation}
  Therefore the absolute value of the  error term in (\ref{expep}) is  
\begin{equation}
  {8   \over \pi}\int_{0}^{\ff}   {\sin^{2} (p /2) \over \psi^{2}(p)} \(1-e^{-t\psi(p)}\) \,dp \le    {8   \over \pi} \int _{0}^{\ff}  {\sin^{2} (p /2) \over \psi^{2}(p)} \(1\wedge t\psi(p)\) \,dp .\label{7.7}
   \end{equation}
   
   We break this last integral into three parts and see that it is bounded by
   \be C\(   t \int_{0}^{\psi^{-1}(1/t)} {p^{2}\over \psi (p)} \,dp +      \int_ {\psi^{-1}(1/t)} ^{1 }{p^{2} \over \psi^{2} (p)} \,dp +\int_{1 }^\ff  {1 \over \psi^{2} (p)} \,dp\)\label{7.8}
 \ee 
 We have
   \begin{equation}
   t \int_{0}^{\psi^{-1}(1/t)} {p^{2}\over \psi (p)} \,dp\le C t^{2}\(\psi^{-1}(1/t)\)^{3}\label{7.9}.
   \end{equation}
   In addition 
   \begin{equation}
  \int_{1 }^\ff  {1 \over \psi^{2} (p)} \,dp\le C .\label{7.10}
   \end{equation}
If $\bb>3/2$
\begin{equation}
      \int_ {\psi^{-1}(1/t)} ^{ 1}{p^{2} \over \psi^{2} (p)} \,dp \le C t^{2}\(\psi^{-1}(1/t)\)^{3}.\label{7.11}
   \end{equation}
If $\bb=3/2$
\begin{equation}
     \int_ {\psi^{-1}(1/t)} ^{1 }{p^{2} \over \psi^{2} (p)} \,dp \le C L(t) \label{7.12}  \end{equation}
for some function $L $ that is slowly varying at infinity. 
If $\bb<3/2$
\begin{equation}
     \int_ {\psi^{-1}(1/t)} ^{1}{p^{2} \over \psi^{2} (p)} \,dp \le C . \label{7.13}\end{equation}
Using (\ref{7.7})--(\ref{7.13}) we get (\ref{4.9}). 

\medskip	Let 
\be
Z=\int ( L^{ x+1}_{t}- L^{ x}_{ t})^{ 2}\,dx.\label{7.14}
\ee
 We get an upper bound for the variance of  $Z$ by finding an upper bound for $EZ^{2}$ and using  (\ref{expep}) to estimate $(EZ)^{2}$.  We proceed as in the beginning of the proof of Lemma \ref{lem-multiple}, however there are enough differences that it is better to repeat some of the arguments.
 
 By the   Kac  Moment Theorem  
\begin{eqnarray}
 & &
E\(\prod_{i=1}^{2}\(\De^{1}_{x_{i}}L_{t}^{x_{i}}\)\(\De^{1}_{y_{i}}L_{t}^{y_{i}}\)\)\label{kacv}\\
& &\quad= \prod_{i=1}^{2}\(\De^{1}_{x_{i}}\De^{1}_{y_{i}}\)\sum_{\si}\int_{\{\sum_{i=1}^{4}r_{i}\leq t\}}\prod_{i=1}^{4}p_{r_{i}}( \si(i) - \si(i-1))\,\,    \prod_{i=1}^{4}\,dr_{i} \nn
\end{eqnarray}
where the sum runs over all  bijections    $\si:\,[1,4]\mapsto \{x_{i},y_{i},\,1\leq i\leq 2\}$   and we take $\si(0)=0$.   
We rewrite (\ref{kacv}) 
 so that each   $\De^{1}_{\cd}$ applies to a single $p_{\cd}$ factor and then set  $y_{i}=x_{i}$ and then integrate with respect to $x_{1},\ldots,x_{m}$ to get 
\begin{eqnarray}
&&E\(\(\int ( L^{ x+1}_{t}- L^{ x}_{ t})^{ 2}\,dx\)^{2}\)
\label{kacv2}\\
&&\qquad =4\sum_{\pi,a}  \int \int_{\{\sum_{i=1}^{4}r_{i}\leq t\}}\prod_{i=1}^{4}
\(\De^1_{x_{\pi(i)}}\)^{{a_{1}(i)}}\(\De^1_{x_{\pi(i-1)}}\)^{{a_{2}(i)}}\nn\\
&&\hspace{1 in}p^{\sharp}_{r_{i}}(x_{\pi(i)}-x_{\pi(i-1)})\,\,    \prod_{i=1}^{4}\,dr_{i}
\,\,    \prod_{i=1}^{2}\,dx_{i}. \nonumber
\end{eqnarray}

In (\ref{kacv2}) the sum runs over all maps $\pi:\,[1,4]\mapsto [1,2]$ with $|\pi^{-1}(i)|=2$ for each $i$ and over all    
$a=(a_{ 1},a_{ 2})\,:\,[1,\ldots, 4]\mapsto \{ 0,1\}\times \{ 0,1\}$ with the
property that for each $i$ there   are   exactly two factors of the form $\De^{
1}_{ x_{i}}$. The factor 4 comes from the fact that we can interchange each $y_{i}$ and $x_{i}$, $i=1,2$.   As usual we take $\pi(0)=0$. 

  As we did in Section \ref{sec-3}, we continue the analysis with $p^{\sharp}$ replaced by $p$.  

Note  that in (\ref{kacv2}) it is possible to have  \label{page79} `bound states', that is values of $i$ for which $\pi (i)=\pi (i-1)$. 
 We first consider the terms in (\ref{kacv2}) with two bound states. There are  two possible maps. They are  $(\pi(1),\pi(2),\pi(3),\pi(4))=(1,1,2,2)$ and $(\pi(1), \pi(2),  \pi(3), 	\pi(4))=(2,2,1, 1)$. The terms in (\ref{kacv2}) for the  map   $(\pi(1),\pi(2),\pi(3),\pi(4))=(1,1,2, 	2)$ are of the form
 
 \begin{equation}
   \prod_{i=1}^{4}
\(\De^1_{x_{\pi(i)}}\)^{{a_{1}(i)}}\(\De^1_{x_{\pi(i-1)}}\)^{{a_{2}(i)}} p _{r_{i}}(x_{\pi(i)}-x_{\pi(i-1)})\label{7.16},
   \end{equation}
 where the  density terms have the form
 \begin{equation}
   p_{r_{1}}(x_{1})p_{r_{2}}(y_{1}-x_{1} )p_{r_{3}}(x_{2}-y_{1} )p_{r_{4}}(y_{2}-x_{2} ),
   \end{equation}
 and where $y_{i}-x_{i}=0$,  $1=1,2$. 
 
 The value of the integrals of the terms in (\ref{7.16}) depend upon how the difference operators are distributed. In many cases the integrals are equal to zero. For example suppose we have
 \begin{equation}
  \De^{1}_{x_{1}}   p_{r_{1}}(x_{1})  \De^{1}_{x_{1}}p_{r_{2}}(0 )  \De^{1}_{x_{2}}p_{r_{3}}(x_{2}-x_{1} )\De^{1}_{x_{2}}p_{r_{4}}( 0),
   \end{equation}
   which we obtain by setting  $y_{1}=x_{1}$.   Written out this term is 
 \bea
&&    \(p_{r_{1}}(x_{1}+1)-p_{r_{1}}(x_{1} ) \) \De^{1}_{x_{1}}p_{r_{2}}(0 ) \\
&&\qquad \(p_{r_{3}}(x_{2}-x_{1}+1 )-p_{r_{3}}(x_{2}-x_{1} )\)\De^{1}_{x_{2}}p_{r_{4}}( 0)\nn
   \eea
  By a change of variables one sees that the integral of this term with respect to $x_{1}$ and $x_{2}$ is zero. 
  
  The only non-zero integrals  in  (\ref{7.16}) comes from 
  \begin{equation}
   p_{r_{1}}(x_{1})  \De^{1}\De^{-1}p_{r_{2}}(0 )   p_{r_{3}}(x_{2}-x_{1} )\De^{1}\De^{-1}  p_{r_{4}}( 0).
   \end{equation}
  The integral of this term with respect to $x_{1}$ and $x_{2}$ is
    \begin{equation}
   \De^{1}\De^{-1}p_{r_{2}}(0 )   \De^{1}\De^{-1}  p_{r_{4}}( 0).
   \end{equation}
 We get the same contribution when $(\pi(1),\pi(2),\pi(3),\pi(4))=(2,2,1,1)$.
 Consequently, the contribution to  (\ref{kacv2}) of maps with two bound states is  
\begin{eqnarray}
&&
  8\int_{\{\sum_{i=1}^{4}r_{i}\leq t\}}\De^1 \De^{-1} p_{r_{2}}( 0)\, \De^{1 }\De^{-1} p_{r_{4}}( 0)\,\prod_{i=1}^{4}\,dr_{i}\label{kacv.8} \\
  &&\qquad=32\int_{\{\sum_{i=1}^{4}r_{i}\leq t\}}\(p_{r_{2} }( 0)-p_{r_{2} }( 1)\)\, \(p_{r_{4} }( 0)-p_{r_{4} }( 1)\)\,\prod_{i=1}^{4}\,dr_{i}\nn\\
  &&\qquad=16\int_{\{u+v\leq t\}}(t-u-v)^{2}\(p_{u }( 0)-p_{u }( 1)\)\, \((p_{v}( 0)-p_{v }( 1)\)\,du\,dv.\nn\\
     &&\qquad\le 16t^{2}\(\int_{0}^{\ff} \(p_{u }( 0)-p_{u }( 1)\)\, du\)^{2}=(4c_{\psi,0}t)^{2}  ,\nn
\end{eqnarray}
(see (\ref{2.13})).

\medskip	
We next consider the contribution from terms with exactly one bound state. These   come from maps of the form $(\pi(1),\pi(2),\pi(3),\pi(4))=(1,2,2,1)$ or $(\pi(1),\pi(2),\pi(3),\pi(4))  =(2,1,1,2)$. These terms  give   non-zero contributions of the form
\begin{eqnarray}
\lefteqn{Q_{2}:=\int\int_{\{\sum_{i=1}^{4}r_{i}\leq t\}}p_{r_{1}}(x) \De^{1}_{x}p_{r_{2}}(y-x)\, \De^1_{y}\De^{-1}_{y}p_{r_{3}}( 0)\,\De^1_{x} p_{r_{4}}(x-y)
\nn}\\
&&  \hspace{3 in} \,\prod_{i=1}^{4}\,dr_{i} \,dx\,dy\label{kacv.14}\\
&&\qquad 
=\int\int_{\{\sum_{i=1}^{4}r_{i}\leq t\}} \De^{-1} p_{r_{2}}(y )\, \De^1 \De^{-1} p_{r_{3}}( 0)\,\De^{-1}  p_{r_{4}}(y) \,\prod_{i=1}^{4}\,dr_{i} \,dy\nonumber;\\
\lefteqn{Q_{3}:=\int\int_{\{\sum_{i=1}^{4}r_{i}\leq t\}}p_{r_{1}}(x) \De^{1}_{x} \De^1_{y}p_{r_{2}}(y-x)\,  p_{r_{3}}( 0)\,\De^1_{x} \De^h_{y} p_{r_{4}}(x-y)
\nn}\\
&&  \hspace{3 in} \,\prod_{i=1}^{4}\,dr_{i} \,dx\,dy\label{pkacv.17aa}\\
&&
=\int\int_{\{\sum_{i=1}^{4}r_{i}\leq t\}} \De^1\De^{-1}p_{r_{2}}(y )\,   p_{r_{3}}( 0)\,\De^1\De^{-1} p_{r_{4}}(y) \,\prod_{i=1}^{4}\,dr_{i} \,dy;\nonumber
\end{eqnarray} 
and 
\begin{eqnarray}
\lefteqn{Q_{4}:=\int\int_{\{\sum_{i=1}^{4}r_{i}\leq t\}}p_{r_{1}}(x) \De^{1}_{x} \De^1_{y}p_{r_{2}}(y-x)\, \De^1_{y} p_{r_{3}}( 0)\,\De^h_{x}  p_{r_{4}}(x-y)
\nn}\\
&&  \hspace{3 in} \,\prod_{i=1}^{4}\,dr_{i} \,dx\,dy\label{pkacv.17bb}\\
&&
=\int\int_{\{\sum_{i=1}^{4}r_{i}\leq t\}} \De^1\De^{-1}p_{r_{2}}(y )\, \De^1  p_{r_{3}}( 0)\,\De^{-1} p_{r_{4}}(y) \,\prod_{i=1}^{4}\,dr_{i} \,dy.\nonumber
\end{eqnarray} 

For further explanation consider  $Q_{2}$. This arrangement comes from the sequence $(x_{1},y_{2},x_{2},y_{1} )$. The expression it is equal to comes by making the change of variables, $y-x\to y$ and then integrating with respect to $x$. 

Integrating and using    (\ref{2.4j}),  (\ref{jrst1.3yj}) and (\ref{151ee1}) we see that
   \bea
  | Q_{2}|&\le  &t \(\int_{0}^{t}|\De^{1}\De^{-1} p_{s }(0 )|\,\,ds\)\,\,\int \( \int_{0}^{t}|\De^{-1} p_{r }(y )|\,\,dr \)^{2}\,dy  
\label{pkacv.15}\\
& \leq &  t\,w(0,t)\sup_{x}v(x,t)   \int v(y,t) \,dy  \nonumber\\
&\leq &Ct^{2}\psi^{-1}(1/t) (\log t)^{2}.\nonumber
      \eea

 To obtain a bound for  $Q_{3}$ we use   (\ref{jr.1j}) and (\ref{132ee2}) to see that  it is bounded in absolute value by 
\bea
&&
t \(\int_{0}^{t} p_{s }(0 )\,\,ds\)\,\,\int \( \int_{0}^{t}|\De^{1} \De^{-1} p_{r }(y )|\,\,dr \)^{2}\,dy\nonumber\\
&&\qquad= tu(0,t)   \int w^{2}(y,t) \,dy \nonumber\\
&& \qquad\label{pkacv.17a}\leq C t^{2}\psi^{-1}(1/t)\log t.
\eea
Integrating $Q_{4}$    we see that  it is bounded in absolute value by  
\bea
&&
  t\int_{0}^{t} \bigg|\De^1  p_{r }( 0)\bigg|\,dr  \int\(\int_{0 }^t \left| \De^1\De^{-1}p_{r }(y )\right|\, dr       \int_{0 }^t  \left| \De^{-1}p_{r }(y )\right|\, dr  \)\,dy\nonumber\\
  &&  \qquad\leq    tv(0,t)\sup_{x}v(x,t)   \int w(y,t) \,dy      \nonumber\\
   &&\qquad \le Ct(\log t)^{3},
   \eea
by (\ref{2.4j})  and (\ref{132ee1}).    \medskip

Finally, we consider  the contribution from terms in (\ref{kacv2}) with no bound states. These have to be from $\pi$ of the form $(\pi(1),\pi(2),\pi(3),\pi(4))\newline =(1,2,1,2)$ or of the form $(\pi(1),\pi(2),\pi(3),\pi(4))=(2,1,2,1)$. 
They give contributions of the form
\begin{eqnarray}
\lefteqn{Q_{5}\label{pkacv.20}}\\
&&:=\int\int_{\{\sum_{i=1}^{4}r_{i}\leq t\}}p_{r_{1}}(x) \De^{1}_{x}p_{r_{2}}(y-x)\, \De^1_{y}\De^{1}_{x}p_{r_{3}}( x-y)\,\De^1_{y} p_{r_{4}}(y-x)\nn 
\\
&&  \hspace{3.3 in} \,\prod_{i=1}^{4}\,dr_{i} \,dx\,dy\nn\\
&&
=\int\int_{\{\sum_{i=1}^{4}r_{i}\leq t\}} \De^{-1}p_{r_{2}}(y )\, \De^1 \De^{-1}p_{r_{3}}( y)\,\De^{1} p_{r_{4}}(y) \,\prod_{i=1}^{4}\,dr_{i} \,dy\nonumber
\end{eqnarray} 
and  
\begin{eqnarray}
\lefteqn{Q_{6}\label{pkacv.21}}\\
&&:=\int\int_{\{\sum_{i=1}^{4}r_{i}\leq t\}}p_{r_{1}}(x) \De^{1}_{x} \De^1_{y}p_{r_{2}}(y-x)\,  p_{r_{3}}( x-y)\,\De^h_{x} \De^1_{y} p_{r_{4}}(x-y)
\nn\\
&&  \hspace{3.3 in} \,\prod_{i=1}^{4}\,dr_{i} \,dx\,dy\nn\\
&&
=\int\int_{\{\sum_{i=1}^{4}r_{i}\leq t\}} \De^1\De^{-1}p_{r_{2}}(y )\,   p_{r_{3}}( y)\,\De^1\De^{-1} p_{r_{4}}(y) \,\prod_{i=1}^{4}\,dr_{i} \,dy\nonumber.
\end{eqnarray} 

Clearly   
\bea
&&  |Q_{5}|   
 \le t \int   \(\int_{0}^{t}|  \De^{-1}p_{r }(y ) |\,dr \)\\
&&\qquad\qquad\qquad  \(\int_{0}^{t}|  \De^{1}p_{r }(y ) |\,dr \) \(\int_{0}^{t}|  \De^{1}\De^{-1}p_{r }(y ) |\,dr \)\,dy\nonumber\\
&&\qquad\leq t \sup_{x}v^{2}(x,t)   \int w(y,t) \,dy.\nn\\
    &&\qquad\le Ct(\log t )^{3},\nn \eea   
by    (\ref{2.4j})  and (\ref{132ee1}).  

 The term  $Q_{6}$ is bounded the same way we bounded $Q_{3}$ and has the same bound.

 \medskip	 
   We can now  obtain an upper bound for the variance. Note that by (\ref{expep})
 \be 
 \(EZ\)^{2}=   \(E\(\int ( L^{ x+1}_{t}- L^{ x}_{ t})^{ 2}\,dx\)\)^{2}=\(4c_{\psi,0}t\)^{2}+O\( tg( t) \).\label{7.29}
 \ee  
  Therefore, it follows from  (\ref{kacv2}) and (\ref{kacv.8}) 
 that  
 \begin{equation}
   \mbox{Var }Z\le EZ^{2}- \(EZ\)^{2}=\sum_{j=2}^{6}|Q_{j}|+Ctg( t)
   \end{equation}
as $t\to \ff$.  Thus we see that
\be \mbox{Var }\, Z   \label{varep} \leq  C\(tg(t)  +  t^{2}\psi^{-1}(1/t)  \log  t\)  .
\ee 
Note that for all $t$ sufficiently large
\begin{equation}
 t  g(t)\le  (t \psi^{-1}(1/t))^{3}\le C  t^{2}\psi^{-1}(1/t) ,
   \end{equation}
  where we use (\ref{4.77}).  Thus we get (\ref{9.3}).
\qed

\section{Appendix: Kac Moment Formula}\label{sec-Kac}

Let    $X=\{X_{t},t\in R_{+} \}$    denote a symmetric L\'evy process with continuous local time $L=\{L_{t}^{x}\,;\,(x,t)\in R^{1}\times R_{+}\}$.   Since $L$ is continuous   we have the occupation density formula,  
 \begin{equation}
\int^{t}_{0} g(X_{s})\,ds=\int g(x) L_{t}^{x}\,dx\label{odf},
 \end{equation}
 for all continuous functions  $g $ with compact support. (See, e.g.  \cite[Theorem 3.7.1]{book}.)

 Let $f(x)$ be a  continuous  function on $R^{1}$ with   compact  support with $\int f(x)\,dx=1$. Let $f_{\ep,y}(x):={1 \over \ep}f\({x-y \over \ep}\)$.   I.e., $f_{\ep,y} (x)$ is an approximate $\de$-function at $x$. Set
\begin{equation}
L^{ x }_{ t,\ep}=\int_{0}^{t} f_{\ep,x }\(X_{s}\)\,ds. \label{ka.3}
\end{equation}
It follows from  (\ref{odf}) that
\begin{equation}
L^{ x }_{ t }=\lim_{\ep\rar 0}L^{ x }_{ t,\ep}\hspace{.2 in}\mbox{a. s. }\label{ka.0}
\end{equation}
Let $p_{t}(x,y)$ denote the probability density of $X_{t}$.

\begin{theorem} [Kac Moment Formula]\label{KMT}   Let    $X=\{X_{t},t\in R_{+} \}$    denote a symmetric L\'evy process with continuous local time $L=\{L_{t}^{x}\,;\,(x,t)\in R^{1}\times R_{+}\}$. 
For any fixed $0<t<\ff$, bounded continuous $g$, and any $x_{1},\ldots, x_{m},z\in R^{1}$,
\begin{eqnarray}
&& 
E^{z}\(\prod_{i=1}^{m}      L^{ x_{i}}_{ t} \,g(X_{t})\)= \sum_{\pi} \int_{\{\sum_{j=1}^{m} r_{j}\leq t\} } \prod_{j=1}^{m}p_{r_{j}}(x_{\pi(j-1)},x_{\pi(j)})\label{ka.1}\\
&& \hspace{2 in}\(\int p_{t-r_{m}}(x_{\pi(m)},y)g(y)\,dy\)\prod_{j=1}^{m}\,dr_{j},
\nn
\end{eqnarray}
where the sums run over all permutations  $\pi$ of $\{1,\ldots, m\}$ and    $\pi(0) : =0$
and $x_{0}: =z$.

 \end{theorem}
 
\Proof Let
\bea
&&
F_{t}(x_{1},\ldots, x_{m})=\int_{\{\sum_{j=1}^{m} r_{j}\leq t\} } \prod_{j=1}^{m}p_{r_{j}}(x_{ j-1 },x_{ j}) \label{ka.10}\\
&& \hspace{2 in}\(\int p_{t-r_{m}}(x_{m},y)g(y)\,dy\)\prod_{j=1}^{m}\,dr_{j}\nn
\eea
 Then 
\begin{eqnarray}
\lefteqn{ 
E^{z}\(\prod_{i=1}^{m}      L^{ x_{i}}_{ t,\ep}\,g(X_{t}) \)\label{ka.4j} }\\
&&= \sum_{\pi}\int_{\{0\leq t_{\pi(1)}\leq \ldots \leq t_{\pi(m)}\leq t\} }E^{z}\( \prod_{j=1}^{m}f_{\ep,x_{j} }\(X_{t_{\pi(j)}}\)\,g(X_{t}) \)\prod_{j=1}^{m}\,dt_{\pi(j)}
\nn\\
&&= \sum_{\pi}\int_{\{0\leq t_{1}\leq \ldots \leq t_{m}\leq t\} }E^{z}\( \prod_{j=1}^{m}f_{\ep,x_{\pi(j)} }\(X_{t_{j}}\) \,g(X_{t})\)\prod_{j=1}^{m}\,dt_{j}
\nn\\
&&= \sum_{\pi}\int \int_{\{\sum_{j=1}^{m} r_{j}\leq t\} } \prod_{j=1}^{m}f_{\ep,x_{\pi(j)} }(y_{j})p_{r_{j}}(y_{j-1},y_{j})\\
&& \hspace{1.5 in}\(\int p_{t-r_{m}}(y_{m},y)g(y)\,dy\)\prod_{j=1}^{m}\,dr_{j}\nn\,dy_{j}
\nn\\
&&= \sum_{\pi}\int F_{t}( y_{0},\ldots,  y_{m})\prod_{j=1}^{m}f_{\ep,x_{\pi(j)} }(y_{j})\,dy_{j}
\nn
\end{eqnarray}
where $ y_{0} :=z$.

 Since the integrand in (\ref{ka.10}) is dominated by $(2\pi)^{-m/2}\prod_{j=1}^{m }r_{j}^{-1/2}$
it follows from the Dominated Convergence Theorem  that $F_{t}(x_{1},\ldots, x_{m})$
is a continuous function of  $(x_{1},\ldots, x_{m}) $ for all $0\leq t<\ff$ and all $m$. 
It then follows immediately from (\ref{ka.4j}) and the fact that  $\prod_{j=1}^{m}f_{\ep,x_{\pi(j)} }(y_{j})$ has compact support that
\begin{eqnarray}
&& 
\lim_{\ep\rar 0}E\(\prod_{i=1}^{m}      L^{ x_{i}}_{ t,\ep}\,g(X_{t}) \) = \sum_{\pi}F_{t}( x_{\pi(0)},x_{\pi(1)},\ldots,  x_{\pi(m)}).
\label{ka.11}
\end{eqnarray}
 A repetition of the above proof  shows that $E\(\lc \prod_{i=1}^{m}      L^{ x_{i}}_{ t,\ep}\rc^{2} \)$
is bounded uniformly in $\ep>0$.  This fact and  (\ref{ka.0})   show that 
\begin{equation}
\lim_{\ep\rar 0}E\(\prod_{i=1}^{m}      L^{ x_{i}}_{ t,\ep}\,g(X_{t}) \)=E\(\prod_{i=1}^{m}      L^{ x_{i}}_{ t }\,g(X_{t}) \).\label{ka.11j}
\end{equation}
Obviously (\ref{ka.11}) and (\ref{ka.11j}) imply (\ref{ka.1}). \qed

   \def\noopsort#1{} \def\printfirst#1#2{#1}
\def\singleletter#1{#1}
      \def\switchargs#1#2{#2#1}
\def\bibsameauth{\leavevmode\vrule height .1ex
      depth 0pt width 2.3em\relax\,}
\makeatletter
\renewcommand{\@biblabel}[1]{\hfill#1.}\makeatother
\newcommand{\bysame}{\leavevmode\hbox to3em{\hrulefill}\,}

\bigskip
\noindent
\begin{tabular}{lll} 
      & Jay Rosen & Michael Marcus\\
      & Department of Mathematics& Department of Mathematics\\
     &College of Staten Island, CUNY& City College, CUNY\\
     &Staten Island, NY 10314& New York, NY 10031\\ &jrosen30@optimum.net &mbmarcus@optonline.net
\end{tabular}

\end{document}